\documentclass[a4paper,11pt]{amsart}

\newfont{\cyr}{wncyr10 scaled 1100}

\usepackage[left=2.7cm,right=2.7cm,top=3.5cm,bottom=3cm]{geometry}
\usepackage{amsthm,amssymb,amsmath,amsfonts,mathrsfs,amscd}
\usepackage{verbatim}
\usepackage[utf8]{inputenc}
\usepackage[all]{xy}
\usepackage{latexsym}
\usepackage{array,longtable,booktabs}
\usepackage{multirow}
\usepackage[usenames,dvipsnames]{color}
\usepackage{listings}
\usepackage{mathtools}
\usepackage{xcolor}
\usepackage{tikz}
\usepackage{tikz-cd}
\usetikzlibrary{matrix,arrows,decorations.pathmorphing}
\usepackage{bm}
\usepackage[shortlabels]{enumitem}
\usepackage[pagebackref]{hyperref}
\hypersetup{hidelinks,pdftitle={Irregular weight-one Eisenstein series and Euler systems},pdfauthor={Javier Polo and Oscar Rivero}}

\theoremstyle{plain}
\newtheorem*{thmA}{Theorem~A}
\newtheorem*{thmB}{Corollary~B}
\newtheorem*{thmC}{Theorem~C}
\newtheorem*{thmD}{Theorem~D}

\newtheorem{theorem}{Theorem}[section]
\newtheorem{hyp}[theorem]{Hypothesis}
\newtheorem{corollary}[theorem]{Corollary}
\newtheorem{lemma}[theorem]{Lemma}
\newtheorem{proposition}[theorem]{Proposition}

\theoremstyle{definition}
\newtheorem{definition}[theorem]{Definition}
\newtheorem{defi}[theorem]{Definition}
\newtheorem{question}[theorem]{Question}
\newtheorem{examplewr}[theorem]{Example}
\newtheorem{ass}[theorem]{Assumption}

\theoremstyle{remark}
\newtheorem{obswr}[theorem]{Observation}
\newtheorem{remarkwr}[theorem]{Remark}

\newenvironment{remark}[1][]{\ifx\relax#1\relax\begin{remarkwr}\else\begin{remarkwr}[#1]\fi\begin{upshape}}{\end{upshape}\end{remarkwr}}
\newenvironment{example}[1][]{\ifx\relax#1\relax\begin{examplewr}\else\begin{examplewr}[#1]\fi\begin{upshape}}{\end{upshape}\end{examplewr}}

\DeclareMathOperator{\et}{et}

\newcommand{\Q}{\mathbb{Q}}
\newcommand{\Z}{\mathbb{Z}}

\newcommand{\C}{\mathbb{C}}

\newcommand{\Frob}{\mathrm{Fr}}

\newcommand{\ord}{{\mathrm{ord}}}
\newfont{\gotip}{eufb10 at 12pt}

\newcommand{\cO}{{\mathcal O}}

\newcommand{\hf}{{\mathbf{f}}}

\numberwithin{equation}{section}

\begin{document}

\title[Irregular weight-one Eisenstein series]{Irregular weight-one Eisenstein series and Euler systems}

\author{Javier Polo and \'Oscar Rivero}

\begin{abstract}
Let $f=E_1(\chi_1,\chi_2)$ be an irregular Eisenstein series of weight one, that is, with $\chi_1(p)=\chi_2(p)$, and consider the generically cuspidal ordinary family passing through it. When $\chi_1$, $\chi_2$ and $\chi_1\chi_2^{-1}$ are primitive of the same conductor, we show that the specialization at $f$ of the modular-symbol $p$-adic $L$-function of this family, with periods taken in the cohomology of the compactified modular curves, factors as a product of two Kubota--Leopoldt $p$-adic $L$-functions attached to $\chi_1$ and $\chi_2$. Since there are no complex critical values in weight one, the factorization is obtained directly from the cohomological construction: we compute the tame-twisted modular measure by means of the localized Sharifi map and explicit reciprocity, and show that tame twisting identifies the compact period modules at weight one. This gives a weight-one analogue of the factorization of Bella\"iche--Dasgupta for critical-slope Eisenstein series. As applications, and under explicit specialization hypotheses, we compare Kato classes with circular units, Heegner classes with elliptic units, and diagonal cycles with Kato classes through two successive Eisenstein degenerations. We also discuss the passage from compact to open periods, and include an exposition of the congruence and degeneration methods underlying these comparisons.
\end{abstract}

\date{\today}

\address{J. P.: Departamento de Matem\'aticas (Universidade de Santiago de Compostela), Santiago de Compostela, Spain}
\email{javier.polo.noche@rai.usc.es}

\address{O. R.: CITMAga and Departamento de Matem\'aticas, Universidade de Santiago de Compostela, Santiago de Compostela, Spain}
\email{oscar.rivero@usc.es}


\maketitle

\setcounter{tocdepth}{1}
\begingroup\small\tableofcontents\endgroup

\section{Introduction}

The decomposition of an Eisenstein Galois representation into characters suggests a corresponding factorization of its $p$-adic $L$-function into Dirichlet $p$-adic $L$-functions. At an Eisenstein point lying on a generically cuspidal family, the function we consider is obtained by specializing a modular-symbol measure in that family, with periods chosen in a cohomological realization; its value at the Eisenstein point therefore depends on these periods, and the expected factorization need not follow from interpolation of complex values. In a primitive common-conductor setting, we prove that, when the periods are taken in the cohomology of the compactified modular curves, the specialization of this measure at an irregular weight-one Eisenstein point $E_1(\chi_1,\chi_2)$ is nonzero and factors as a product of two Kubota--Leopoldt $p$-adic $L$-functions attached to $\chi_1$ and $\chi_2$, and we then study the consequences for Euler systems.

Our main result may be seen as a weight-one analogue of the factorization established by Bella\"iche--Dasgupta \cite[Theorem~1.1]{BD15}. For a critical-slope Eisenstein series of weight at least two, their formula expresses the $p$-adic $L$-function, on the appropriate sign component, as the product of two Dirichlet $p$-adic $L$-functions and an explicit cyclotomic logarithm. Here we consider an ordinary Eisenstein series of weight one whose two Hecke roots at $p$ coincide, that is, an irregular one; since weight one has no complex critical integers, we must obtain the factorization from the modular-symbol construction in families.

Our starting point is the work of Shih--Wang \cite{SW21,SW22} on the exceptional case of Sharifi's conjecture, whose factorization concerns a weight-one Eisenstein series with one trivial constituent and uses an integral Gorenstein hypothesis to choose periods. We work instead with the cohomology of the compactified modular curves in characteristic zero, where local freeness follows from Betina--Dimitrov--Pozzi, and pass to two nontrivial constituents by tame twisting. This requires control of the twisting correspondence on the period modules, including its nonvanishing at the Eisenstein point, which we establish by evaluating the twisted modular symbols through local reciprocity. The same calculation thus gives the two Dirichlet factors and the nonvanishing needed to identify the compact periods at weight one.

A related factorization occurs for Katz $p$-adic $L$-functions restricted to the cyclotomic line: in the CM case, \cite[Corollary~2.3]{SW22} connects the modular $p$-adic $L$-function with the Katz $p$-adic $L$-function and Gross's formula. Our common-conductor result treats a different range, where both constituents are nontrivial. More broadly, the comparison of Euler systems as a way of proving instances of $p$-adic Artin formalism is also central in Dasgupta's factorization of the Rankin $p$-adic $L$-function of a tensor square \cite{Das}.

We then compare Kato classes with circular units, and consider analogous questions for Heegner classes and for two successive Eisenstein degenerations of diagonal cycles; in each case, the analytic factorizations suggest which Euler systems should occur after specialization. To compare the classes themselves, we also need regular specialization in the chosen cohomological realization and a certain compatibility with the local conditions, which we state as explicit hypotheses.

\subsection{The weight-one factorization}\label{subsec:intro-factorization}

Let $f=E_1(\chi_1,\chi_2)$, with $\chi_1$ even and $\chi_2$ odd, and assume that $f$ is irregular at $p$, that is, $\chi_1(p)=\chi_2(p)$; we write $\alpha$ for this common value. Since the two roots of the Hecke polynomial at $p$ coincide, $f$ has a unique ordinary $p$-stabilization, which we still denote by $f$. The Galois representation of $f$ is the sum of the two characters $\chi_1$ and $\chi_2$, and this is why one expects its $p$-adic $L$-function to involve two Dirichlet factors. The geometry of \cite{BeDiPo} gives a generically cuspidal ordinary branch through $f$, \'etale over weight space; we write $\mathbf f$ for the family and $A=L[[X]]$ for its completed local ring, with $X=0$ at weight one.

The periods in our formula come from ordinary cohomology along the tower of compactified modular curves: its signed torsion-free quotient on which the Hecke operators act through the eigenvalues of $\mathbf f$ is a free $A$-module of rank one, denoted by $\mathscr P_{\mathbf f,\mathrm{cpt}}$. Ordinary modular symbols define a measure with values in this module, and choosing a basis $\gamma_{\mathbf f}$ gives a scalar measure in the family, whose specialization at weight one we denote by $\mathscr L_p^{\mathrm{cpt}}(f)$. We refer to this as the \emph{compact period} normalization, where $\mathrm{cpt}$ refers to the compactified modular curve; the precise quotient, measure and interpolation formula at classical cuspidal points are given in \S\ref{subsec:regular-compact-twisting}.

For $z\in\Z_p^\times$, write $z=\omega(z)\langle z\rangle$, with $\omega$ the Teichm\"uller character. On its trivial component, $\mathscr L_p^{\mathrm{cpt}}(f,j)$ is the evaluation at $z\mapsto\langle z\rangle^j$. We use $L_p^{\mathrm{KL}}(\beta,s)$ for the standard Kubota--Leopoldt function of an even Dirichlet character $\beta$. Theorem~\ref{thm:compact-weight-one-factorization} and Corollary~\ref{prop:cyclotomic-period-comparison} give the following formula.

\begin{thmA}
Suppose that $\chi_1$, $\chi_2$ and $\chi_1\chi_2^{-1}$ are primitive of the same conductor $M>1$, that $p\geq5$ and $p\nmid M\varphi(M)$, and that $\chi_1(p)=\chi_2(p)$. Then
\[
 \mathscr L_p^{\mathrm{cpt}}(f,j)
   =C_{f,\gamma}\langle M\rangle^{-j}
       L_p^{\mathrm{KL}}(\chi_1,1+j)
       L_p^{\mathrm{KL}}(\chi_2^{-1}\omega,-j),
 \qquad C_{f,\gamma}\in L^\times.
\]
The constant is independent of $j$, and the resulting bounded measure is nonzero. Its normalization is determined by the chosen compact period.
\end{thmA}

To prove the theorem, we write $f=f_0\otimes\eta$, where $\eta=\chi_1$ and $f_0=E_1(1,\chi_2\chi_1^{-1})$. The twisting correspondence gives a regular map from the compact period module of $\mathbf f$ to that of the auxiliary family $\mathbf f_0$. We compute the image measure using the localized Sharifi map, which sends its modular symbols to cup products of cyclotomic units. Explicit reciprocity expresses their local invariants as traces: the logarithm of the first entry gives one Dirichlet factor, and the conjugates of the Coleman logarithmic derivative of the second give the other. The nonvanishing of the resulting product shows that the twisting map is an isomorphism on the weight-one fibre, hence also over the local weight ring, and allows us to recover the factorization in the compact periods of $\mathbf f$.

The reflected argument and the conductor factor can also be understood from the point of view of complex Artin formalism. For an even finite-order character $\rho$ of conductor $p^r$, one has
\[
 L(f\otimes\rho^{-1},s)=L(\chi_1\rho^{-1},s)L(\chi_2\rho^{-1},s).
\]
Applying the functional equation to the second factor at $s=1$ replaces it by a value of $L(\chi_2^{-1}\rho,0)$. The Gauss sum of $\chi_2\rho^{-1}$ contributes $\chi_2(p)^r\rho(M)^{-1}\tau_M(\chi_2)\tau_{p^r}(\rho^{-1})$; the first factor cancels the ordinary eigenvalue $\alpha^r$, and the wild Gauss sum cancels the one in the interpolation formula, so that only $\rho(M)^{-1}$ survives. This accounts for $\langle M\rangle^{-j}$ and the argument $-j$ in the theorem. The first factor, at $s=1$ with even character, is represented $p$-adically by Leopoldt's logarithmic formula. This heuristic explains the shape of the formula, although the proof has to proceed through local reciprocity, since complex interpolation is not available at weight one.

For example, take $M=19$, $p=83$ and a character $\vartheta$ modulo $19$ with $\vartheta(2)=\zeta_{18}$. The pair $(\chi_1,\chi_2)=(\vartheta^2,\vartheta^5)$ satisfies the hypotheses: both characters take the value $\zeta_{18}^{12}$ at $83$, and their ratio is primitive. Since $\zeta_{18}^{12}\ne1$, this example also satisfies the additional condition $\alpha\ne1$ needed in Corollary~B below. We return to these data in Example~\ref{ex:common-conductor-heegner}.

The appearance of logarithms of units also connects this result with the weight-one comparisons of Darmon--Lauder--Rotger \cite{DLRfirst} and Rivero--Rotger \cite{RRstark}, where logarithmic regulators describe specializations or deformations of Euler-system classes. In our setting, the first Kubota--Leopoldt factor arises from a sum of logarithms of cyclotomic units in the evaluation of the compact modular-symbol measure.

\subsection{Kato classes and circular units}

The factorization of Theorem~A suggests a corresponding relation between Euler systems. Recall that Kato's classes are attached to cusp forms, while circular units are attached to Dirichlet characters. When a cuspidal family passes through an Eisenstein point, its Galois representation becomes reducible, and it is natural to ask whether the specialization of Kato's class can be expressed in terms of circular units. Section~\ref{sec:circular-kato} approaches this question from the point of view of Euler systems, following the critical-slope comparison of \cite{PR25}.

Write $\kappa(f)_{\chi_1}$ for the $\chi_1$-component of the specialized Kato class and $c(\chi_1)$ for the circular-unit class. In the cohomological normalization used to compare them, both have coefficients in the cyclotomic deformation of $L(\chi_1^{-1})(1)$: the dual and Tate twist are those of Kummer theory. Under our hypotheses, the other character component vanishes on even cyclotomic components for parity reasons. The comparison is suggested by the explicit reciprocity laws: the regulator of the circular unit gives the first Kubota--Leopoldt factor in Theorem~A, and cancelling this factor should leave the second one as the coefficient of $c(\chi_1)$.

The choice of periods matters here, since Theorem~A uses compact periods, whereas Kato's reciprocity law uses those of the Euler-system construction, and their comparison may introduce a zero at the Eisenstein point. The normalization of \cite{PR25} also contributes an auxiliary $L$-value and an adjoint factor. Our strategy is to keep these weight-dependent factors together: if the class and its regulator specialize regularly, the nonvanishing of the compact measure in Theorem~A forces their combined contribution to be regular.

Hypothesis~\ref{hyp:kato-specialization} specifies the required regular specialization: the family class extends to the chosen split realization, and its ordinary regulator becomes the circular-unit Coleman map, with the stated coordinate factor. Let $\mathscr D_f(j)$ denote the combined weight factor at $f$, together with the tame Euler and coordinate factors. Its expression is given in \eqref{eq:conditional-kato-circular-comparison}.

\begin{thmB}
Assume the hypotheses of Theorem~A, that $\alpha\ne1$, and Hypothesis~\ref{hyp:kato-specialization}. On the trivial Teichm\"uller component,
\[
 \kappa(f)_{\chi_1}
   =\mathscr D_f(j)\langle M\rangle^{-j}
       L_p^{\mathrm{KL}}(\chi_2^{-1}\omega,-j)c(\chi_1)
\]
over the fraction field of the cyclotomic distribution algebra. If the coefficient is analytic, the equality holds in global cohomology. The comparison factor defining $\mathscr D_f$ is regular in weight, but its value may be zero, and $\kappa(f)_{\chi_1}\ne0$ if and only if this value is nonzero.
\end{thmB}

Hence the second Dirichlet $p$-adic $L$-function gives the coefficient of the circular-unit class, up to the factors appearing in the two reciprocity laws. The passage from regulators to global classes uses Leopoldt's theorem for abelian fields, as explained in Proposition~\ref{ass:global-rank-one}. Section~\ref{sec:circular-kato} establishes the comparison before evaluating its analytic coefficient; \S\ref{subsec:factorization-application-kato} then deduces Corollary~B from Theorem~A.

The critical-slope comparisons of Loeffler--Rivero \cite{LR1} and Polo--Rivero \cite{PR25} provide the model for this use of cohomological realizations and crystalline periods. We also study the period map itself, and we show that for the Drinfeld--Manin realization of the auxiliary family $\mathbf f_0$ it has a simple zero, so that it carries a nonzero compact value to a first derivative in the weight direction. Question~\ref{conj:kato-period-comparison} asks for the corresponding comparison with the realization selected by Kato's periods.

\subsection{Heegner classes and elliptic units}

After considering circular units over $\Q$, we turn to elliptic units over an imaginary quadratic field and their relation with Heegner classes. This anticyclotomic comparison has a simpler starting point: the geometry places the localization of a Heegner class in the positive ordinary direction, so its image in the natural rank-one quotient is zero before any factorization of $p$-adic $L$-functions is used, and we can work directly with the positive regulator as in \cite[\S C3]{LR1}.

Here $f=E_1(\chi_1,\chi_2)$ is an irregular Eisenstein series as in \S\ref{subsec:intro-factorization}; this comparison does not use Theorem~A, and the split realization it requires is part of the hypotheses of Section~\ref{sec:heegner-elliptic}. Let $K$ be an imaginary quadratic field in which $p$ and every prime of the tame level split. The construction of \cite{JLZ}, with the twists and specialization hypotheses in Section~\ref{sec:heegner-elliptic}, gives a class $\kappa(f,K)$ varying over anticyclotomic characters. The two Eisenstein constituents give rank-one coefficient modules $\mathbb V_i^{\mathrm{ac}}$, with elliptic-unit classes $\mathfrak u_i$. Their local regulators are the Katz functions $\mathscr L_i$, for $i=1,2$.

Under the interpolation compatibility of Hypothesis~\ref{hyp:heegner-interpolation}, the BDP--Katz comparison gives
\[
 \mathscr L_{\mathrm{BDP}}^2=\mathscr A(f,K)\mathscr L_1\mathscr L_2,
\]
where $\mathscr A(f,K)$ records the Euler and period factors. We therefore form
\[
 \mathfrak u(f,K)
 =v_{f,\chi_1}\otimes\mathscr L_2\mathfrak u_1
       +v_{f,\chi_2}\otimes\mathscr L_1\mathfrak u_2.
\]
Each elliptic-unit class is thus weighted by the Katz function of the other constituent. Under the local compatibility of \S\ref{subsec:weighted-elliptic-unit}, this combination satisfies the Heegner local condition and can be compared with the geometric class.

\begin{thmC}
Assume Hypotheses~\ref{hyp:heegner-specialization} and~\ref{hyp:heegner-interpolation}, and Assumption~\ref{ass:elliptic-local-injectivity}. If the Heegner Selmer module is torsion-free of rank one and its positive regulator is nonzero, then
\[
 \mathscr A(f,K)\,\mathfrak u(f,K)
   =\Omega_{f,\gamma}\,
      \mathscr L_{\mathrm{BDP}}\,\kappa(f,K)
\]
after extension to the fraction field of the anticyclotomic distribution algebra. If $\mathscr A(f,K)$ is analytic and both classes belong to the stated Selmer module, the equality holds in that module.
\end{thmC}

Here $\Omega_{f,\gamma}$ is the period in the Heegner reciprocity law. The extra BDP factor comes from comparing the regulator of the Heegner class, which gives $\mathscr L_{\mathrm{BDP}}$, with that of the weighted elliptic-unit class, which gives its square. Both character summands are retained in this comparison.

\subsection{Double Eisenstein degeneration}

Finally, we ask how these comparisons behave under two successive Eisenstein specializations in the diagonal-cycle construction, keeping the third modular family cuspidal. The first specialization leads to two Beilinson--Flach components and the second to four Kato components, each multiplied at the corresponding stage by the $p$-adic $L$-function of its complementary constituent.

Let $\mathbf f$, $\mathbf g$ and $\mathbf h$ be three ordinary families. For this comparison, and throughout Section~\ref{sec:double-degeneration}, we use $f$ and $h$ for the two Eisenstein points:
\[
 \mathbf f_{x_f}=f=E_1(1,\varepsilon_K),\qquad
 \mathbf h_{x_h}=h=E_1(\psi_1,\psi_2),
\]
where $K$ is imaginary quadratic, $\varepsilon_K$ is its quadratic character, and $p$ splits in $K$. Both points are irregular: for $f$ this follows from $\varepsilon_K(p)=1$, since $p$ splits in $K$. We take $\mathbf f$ to be the CM family through $f$, while $h$ is an irregular weight-one point of $\mathbf h$, and $\mathbf g$ remains cuspidal. Although the two points play parallel roles in the successive decompositions of the coefficient representation, the comparisons used at each step are different: the first one uses the CM construction of \cite{AOR}, and the second uses the Beilinson--Flach/Kato comparison of \cite{BDV}, with the specialization and local compatibility requirements of Proposition~\ref{prop:BF-Kato-components}.

The weights come from two analytic factorizations: first, the square of the $g$-dominant triple-product $p$-adic $L$-function $\mathscr L_p^g$ factors into two Hida--Rankin functions, and then each of them factors into two Mazur--Kitagawa functions. Multiplying at both stages each class by the function of its complementary constituent, we obtain the sum $\operatorname{BK}^{(2)}(f,\mathbf g,h)$ of Definition~\ref{def:doubly-weighted-Kato}.

\begin{thmD}
Under the specialization and Selmer hypotheses of Theorem~\ref{thm:double-degeneration}, and with the compatible periods specified there, one has
\[
 \lambda_{N_g}(\mathbf g)\,\mathscr A_h\,
       \operatorname{BK}^{(2)}(f,\mathbf g,h)
 =\Omega_f\Omega_h\,
       \mathscr L_p^g(f,\mathbf g,h)\,
       \kappa(f,\mathbf g,h).
\]
The equality holds over the fraction field of the coefficient ring of the remaining family $\mathbf g$.
\end{thmD}

Here $\Omega_f$ and $\Omega_h$ are the periods of the two specializations, defined in \S\ref{subsec:first-degeneration} and \S\ref{subsec:second-degeneration}, respectively. The factor $\lambda_{N_g}(\mathbf g)$ is the Atkin--Lehner pseudo-eigenvalue recalled in \S\ref{subsec:first-degeneration}, and $\mathscr A_h$ contains the Euler and normalization factors of the second comparison.

In order to carry out both specializations, we require a regular extension to $x_h$ of the classes and reciprocity laws of the first stage. At the second step, \cite{BDV} supplies the local comparison in its CM setting, while for a general $h$ the stated compatibility remains an input; in particular, the formula concerns the indicated order of specialization. The related constructions of Heegner classes of Castell\`a--Do \cite{castella-do} provide further motivation for this picture.

\subsection{Organization of the article}

Section~\ref{sec:preliminaries} fixes the conventions regarding Galois representations and cyclotomic variables. Sections~\ref{sec:congruences-lattices}--\ref{sec:eisenstein-degeneration} are mostly expository: they compare Eisenstein congruences with degeneration in families, emphasizing the role of the boundary in the cohomological realizations. We then turn to the main comparison. Section~\ref{sec:circular-kato} relates Kato classes with circular units, initially with an unsimplified analytic coefficient, and Section~\ref{sec:weight-one-factorization} proves the compact factorization of Theorem~A and uses it to obtain Corollary~B in the common-conductor setting, without the simple-zero assumption of Section~\ref{sec:circular-kato}. Section~\ref{sec:heegner-elliptic} discusses the anticyclotomic counterpart, where the Heegner local condition directly gives the vanishing of the quotient. Finally, Section~\ref{sec:double-degeneration} applies the same weighting procedure to two successive Eisenstein specializations.

Some questions remain open. From the arithmetic point of view, the main ones are the regular transport of the Kato class and its regulator to the split realization, and the nonvanishing of the combined comparison factor in Corollary~B; the comparison between compact and Kato periods would help to determine the latter, although in the normalization used here the orders of the auxiliary and adjoint factors also play a role. It would also be interesting to extend Theorem~A to characters of different conductors, which would require the corresponding local Fourier and Euler-factor computations.

\subsection{Acknowledgements}

We thank Ra\'ul Alonso for useful conversations on the topic of this article. During the preparation of this article, the authors were supported by PID2023-148398NA-I00, funded by MCIU/AEI/10.13039/501100011033/FEDER, UE.

\subsection{Use of AI}

This project was initiated in 2025, before we began using AI tools in our mathematical research. In particular, the idea of applying the results of Shih--Wang to obtain the factorization formula predates our use of these tools. We subsequently used the language models GPT-6 Astra (OpenAI) and Claude Opus~5.5 (Anthropic) primarily to clarify the justification of tame twisting between the two families in Section~\ref{sec:weight-one-factorization} and the comparison of periods arising from open and compactified cohomology. We also used these models to improve the exposition and rewrite parts of the paper, to develop specific examples (such as Example~\ref{ex:common-conductor-heegner}, which we have verified independently), and to identify errors and inaccuracies in earlier drafts. All mathematical statements and proofs have been checked by the authors, who take full responsibility for the content of this article.

\section{Preliminaries}
\label{sec:preliminaries}

We fix the conventions for Galois representations, modular families and periods. The Euler-system constructions use the cohomology of the open modular curve $Y_1$, whereas the periods in Section~\ref{sec:weight-one-factorization} come from its compactification $X_1$. Their generic cuspidal eigenspaces agree, but their integral or coherent quotients may differ at an Eisenstein point.

Let $p\geq5$ be a prime.  We fix embeddings
\[
 \iota_p:\overline{\Q}\hookrightarrow\overline{\Q}_p,
 \qquad
 \iota_\infty:\overline{\Q}\hookrightarrow\C.
\]
The second embedding determines a complex conjugation $c\in G_{\Q}$. Unless explicitly indicated otherwise, $L(\chi)$ stands for the Galois character with arithmetic Frobenius value $\chi(q)$, and the cyclotomic character $\varepsilon$ satisfies $\varepsilon(\Frob_q)=q$. Since \cite{LR1,PR25} use instead the geometric-Frobenius convention, we record the translation between both conventions in \S\ref{subsec:critical-slope-eisenstein}. As usual, if $L/\Q_p$ is a finite extension, we write $\cO_L$ for its ring of integers; when no confusion can arise, we simply write $\cO$.

\subsection{Galois representations of modular forms}
\label{subsec:deligne-representations}

Let
\[
 \xi=\sum_{n\geq1}a_n(\xi)q^n\in S_k(N_\xi,\chi_\xi)
\]
be a normalized newform of weight $k\geq2$.  After enlarging $L$ so that it contains its Hecke field, the geometric realization used in this paper is the two-dimensional $L$-representation $V_\xi$ obtained from the maximal $\xi$-isotypical quotient of
\begin{equation}
 H^1_{\et}\bigl(Y_1(N_\xi)_{\overline\Q},
             \mathscr L_{k-2}(1)\bigr)\otimes_{\Z_p}L.
 \label{eq:geometric-realization-modular-form}
\end{equation}
Here $\mathscr L_{k-2}$ is the standard degree-$(k-2)$ coefficient sheaf on the modular curve, with the convention of \cite[\S2.3]{BSV}, and the dual Hecke operators act through the eigensystem of $\xi$. With this convention, $V_\xi$ is unramified outside $pN_\xi$ and
\[
 \operatorname{tr}\bigl(\Frob_q;V_\xi\bigr)=a_q(\xi)
 \qquad(q\nmid pN_\xi),
\]
where $\Frob_q$ is arithmetic Frobenius.  All twists and duals are taken with respect to this geometric normalization.

If $\xi$ is ordinary at $p$, its restriction to $G_{\Q_p}$ has an ordinary filtration
\begin{equation}
 0\longrightarrow V_\xi^+
  \longrightarrow V_\xi
  \longrightarrow V_\xi^-
  \longrightarrow0,
 \label{eq:ordinary-filtration-classical}
\end{equation}
in which both graded pieces have dimension one and $V_\xi^-$ is unramified. Arithmetic Frobenius acts on $V_\xi^-$ through the unit root of the Hecke polynomial at $p$, and we use the same superscripts for the corresponding filtration in families. Observe that when a dual or a Tate twist is present the characters on the graded pieces change accordingly; to avoid any confusion, we will always record the resulting character at the point where the projection is used.

\subsection{Weight space, families and specialization}
\label{subsec:families-specialization}

Put $\Lambda_{\mathrm{wt}}=\cO[[1+p\Z_p]]$.  An arithmetic point of weight $k\geq1$ is a continuous homomorphism $\nu:\Lambda_{\mathrm{wt}}\to\overline\Q_p$ whose restriction to $1+p\Z_p$ is of the form
\[
 z\longmapsto z^{k-2}\epsilon(z)
\]
for a finite-order character $\epsilon$.  We use boldface letters such as $\mathbf f$, $\mathbf g$ and $\mathbf h$ for Hida or Coleman families, and we write $\mathbb V_{\mathbf f}$ for the associated rank-two Galois module over its coefficient algebra.  On an ordinary branch one has an exact sequence
\begin{equation}
 0\longrightarrow\mathbb V_{\mathbf f}^{+}
  \longrightarrow\mathbb V_{\mathbf f}
  \longrightarrow\mathbb V_{\mathbf f}^{-}
  \longrightarrow0
 \label{eq:ordinary-filtration-family}
\end{equation}
of $G_{\Q_p}$-modules, whose formation is compatible with specialization at the classical points considered below.

At a $p$-regular cuspidal point of weight one we use the \'etaleness criterion of \cite{bellaiche-dimitrov}, while at an irregular Eisenstein point we use the cuspidal branch of \cite{BeDiPo}. Note that in the latter case the full local eigencurve may be singular, even though this branch is \'etale over weight space.

For the CM family through $E_1(1,\varepsilon_K)$ used in \cite[Proposition~2.2]{BDV}, enlarge $L$ so that the weight-one point $x$ is $L$-rational. The chosen cohomological realization satisfies
\begin{equation}
 \mathbb V_{\mathbf f}\otimes_x L
       \simeq L\oplus L(\varepsilon_K),
 \label{eq:weight-one-specialization}
\end{equation}
where the first summand has trivial Galois action. For a non-CM Eisenstein point, however, a split fibre is a property of the chosen lattice, and it does not follow from the semisimple Artin representation. In the common-conductor setting of Section~\ref{subsec:factorization-setting}, \cite[Theorem~1.5]{SW21} supplies a split fibre for the localized compact realization of the auxiliary family. Tensoring that full rank-two realization by $\chi_1$ gives a split realization of the generic representation of our family. The regularity of the Kato class in this realization is a different question, which we discuss in Section~\ref{subsec:kato-weight-one}, together with the local direction selected by the ordinary filtration in the resulting fibre.

Whenever the weight map is \'etale at $x$, we shrink to a disc $U$ on which the weight is a local coordinate. We write $A_U$ for its coefficient algebra, and $X\in A_U$ for a uniformizer vanishing at $x$. Any leading term in the weight direction is understood with respect to $X$; note that ratios of two first derivatives do not depend on this choice of uniformizer.

\subsection{Irregular weight-one Eisenstein series}
\label{subsec:wt1}

Let $\chi_1$ and $\chi_2$ be primitive Dirichlet characters of conductors $N_1$ and $N_2$, prime to $p$, and assume
\[
 \chi_1\chi_2(-1)=-1.
\]
The associated weight-one Eisenstein series is
\begin{equation}
 E_1(\chi_1,\chi_2)
  =c_0(\chi_1,\chi_2)
   +\sum_{n\geq1}
       \left(\sum_{d\mid n}\chi_1(n/d)\chi_2(d)\right)q^n,
 \label{eq:weight-one-eisenstein-series}
\end{equation}
where
\[
 c_0(\chi_1,\chi_2)=
 \begin{cases}
  \frac12L(0,\chi_2),&\text{if }\chi_1=\mathbf1,\\
  \frac12L(0,\chi_1),&\text{if }\chi_2=\mathbf1,\\
  0,&\text{otherwise}.
 \end{cases}
\]
It has weight $1$, tame level $N_1N_2$ and nebentype $\chi_1\chi_2$. The two roots of its Hecke polynomial at $p$ are $\chi_1(p)$ and $\chi_2(p)$.  We say that it is \emph{$p$-irregular} if
\begin{equation}
 \chi_1(p)=\chi_2(p).
 \label{eq:p-irregular-condition}
\end{equation}
In this case there is a unique ordinary $p$-stabilization, and the results of \cite{BeDiPo} give a generically cuspidal ordinary branch through this point, on which the weight map is \'etale. We denote by $\mathbf f$ a family on this branch.

The semisimplification of the representation at the Eisenstein point is
\begin{equation}
 V_f^{\mathrm{ss}}\simeq L(\chi_1)\oplus L(\chi_2).
 \label{eq:eisenstein-semisimplification-general}
\end{equation}
However, the precise characters appearing in a projected Euler-system class depend on its coefficient representation. For instance, in the Kummer normalization of Section~\ref{sec:circular-kato}, we use $\mathbb W_{\mathbf f}=\mathbb V_{\mathbf f}^{*}(1)$, whose semisimplification at weight one has constituents $L(\chi_i^{-1})(1)$. The Tate twist comes here from the realization used to compare classes with units, and it is not present in the Artin representation \eqref{eq:eisenstein-semisimplification-general}. The splitting of the chosen fibre, and its compatibility with the Kato class and its ordinary regulator, are discussed separately in Section~\ref{subsec:kato-weight-one}.

A particularly important case is
\[
 (\chi_1,\chi_2)=(\mathbf1,\varepsilon_K),
\]
where $K/\Q$ is an imaginary quadratic field in which $p$ splits.  Then $\varepsilon_K(p)=1$, the unique $p$-stabilization is
\[
 E_1(\mathbf1,\varepsilon_K)(q)
  -E_1(\mathbf1,\varepsilon_K)(q^p),
\]
and the chosen cuspidal branch is a CM Hida family associated with $K$. This is the branch that we use in the double-degeneration arguments of \S\ref{sec:double-degeneration}.

\subsection{Periods and the adjoint $p$-adic $L$-function}
\label{subsec:adjoint-p-adic-L}

Let $\mathbf f$ be a generically cuspidal family over $A_U$.  A choice of bases $\gamma_{\mathbf f}^{\pm}$ in the two complex-conjugation eigenspaces determines compatible periods $\Omega_{f_x}^{\pm}$ at its classical cuspidal specializations.  With these choices, Hida's construction \cite{hida-adjoint}, in the eigenvariety formulation of \cite{eigenbook}, gives an adjoint $p$-adic $L$-function
\[
 L_p(\operatorname{Ad}\mathbf f)\in A_U.
\]
At classical points of weight at least two, it interpolates $L^{\mathrm{alg}}(\operatorname{Ad}f_x,1)$ multiplied by the standard $p$-Euler factors, and at a $p$-distinguished classical point it generates, up to a unit, the congruence ideal of the branch.

At an irregular Eisenstein point, however, the residual representation is not $p$-distinguished, and the previous statement about the congruence ideal cannot simply be specialized. We deal with this issue in two different ways. For the comparison with general conductors we specify the order of the adjoint factor in the chosen periods, and whenever a simple zero is required we state it explicitly as an assumption (see Assumption~\ref{ass:simple-adjoint-zero}). For the common-conductor result, instead, we keep the adjoint factor together with the auxiliary and period factors in the full reciprocity coefficient, so that Corollary~\ref{cor:conditional-kato-circular-comparison} does not require any separate assertion about the order of the adjoint factor.

\subsection{Cyclotomic variables and Dirichlet functions}
\label{subsec:cyclotomic-conventions}

We identify $\Gamma=\operatorname{Gal}(\Q(\mu_{p^\infty})/\Q)$ with $\Z_p^\times$ by the cyclotomic character and write
\[
 z=\omega(z)\langle z\rangle,\qquad
 \Gamma=\Delta\times\Gamma_1,\qquad
 \Delta=\mu_{p-1},\quad\Gamma_1=1+p\Z_p.
\]
Here $\omega$ stands for the Teichm\"uller character. Put $\Lambda_{\mathrm{cyc}}=\cO_L[[\Gamma]]$, $\Lambda_L=\cO_L[[\Gamma_1]][1/p]$, and let $\mathcal H_\Gamma$ be the algebra of $L$-valued locally analytic distributions. A fixed $\omega^a$-component is obtained by means of the idempotent $e_a\in L[\Delta]$, and its bounded measures form a copy of $\Lambda_L$; we say that the component is even when $a$ is even. The trivial Teichm\"uller component corresponds to $a=0$, that is, to characters which are trivial on $\Delta$ but vary on $\Gamma_1$.

For a measure $\mu$, the notation $\mu(\rho)$ means $\int_\Gamma\rho(z)\,d\mu(z)$, with $\rho$ of finite order. On the $\omega^a$-component we write
\[
 \mu(j)=\int_\Gamma\omega(z)^a\langle z\rangle^j\,d\mu(z).
\]
Hence finite-order characters and analytic exponents are just two kinds of evaluations of the same measure. The Dirac measure $[m]$ at an integer $m$ prime to $p$ has values $\rho(m)$ and $\omega(m)^a\langle m\rangle^j$, respectively. We use the involution $\iota:[z]\mapsto[z^{-1}]$ when translating the inverse cyclotomic deformation in cohomology into this Mellin convention. On a fixed component, $\mathcal H_\Gamma$ is the ring of analytic functions on the cyclotomic open disc, which is an integral domain, and we write $\mathcal K_\Gamma$ for its fraction field.

For an even Dirichlet character $\beta$, $L_p^{\mathrm{KL}}(\beta,s)$ denotes the standard Kubota--Leopoldt function, with convention
\[
 L_p^{\mathrm{KL}}(\beta,1-n)
 =\bigl(1-(\beta\omega^{-n})(p)p^{n-1}\bigr)
       L(\beta\omega^{-n},1-n),\qquad n\geq1.
\]
We set the Euler factor to $1$ when the character is ramified at $p$. The value at $s=1$ is understood by $p$-adic continuation; for a primitive nontrivial even character $\nu$ of conductor $m$ divisible by $p$, Leopoldt's logarithmic formula reads
\[
 L_p^{\mathrm{KL}}(\nu,1)
   =-\frac{\tau_m(\nu)}{m}
       \sum_{b\in(\Z/m\Z)^\times}\nu^{-1}(b)\log_p(1-\zeta_m^b),
 \qquad \tau_m(\nu)=\sum_b\nu(b)\zeta_m^b.
\]
Here the roots of unity are chosen compatibly. This fixes the normalization of the Dirichlet $p$-adic $L$-functions used both in the circular-unit regulator and in the modular-symbol calculation of Sections~\ref{sec:circular-kato}--\ref{sec:weight-one-factorization}.

\section{Eisenstein congruences and Galois lattices}
\label{sec:congruences-lattices}

In this section we recall the integral comparison between the open and the compactified modular curve developed in \cite{RiRo}. The two curves give lattices in the same rational Galois representation, whose Eisenstein quotients may carry opposite residual extensions. This phenomenon connects Sharifi's conjecture with the specialization of Euler-system classes, and it is the integral prototype of the period comparisons of the following sections. We point out that the weight-two form $f$, the congruence exponent $t$ and the curves $\overline Y$, $Y$ are only used in this section.

\subsection{The Eisenstein congruence}
\label{subsec:eisenstein-congruence}

Let $p \geq 5$ be a prime, let $N$ satisfy $p\nmid N\varphi(N)$, and let \(\theta:(\mathbf Z/N\mathbf Z)^\times\to \mathcal O^\times\) be an even Dirichlet character.  Here $\mathcal O$ stands for the ring of integers of a finite extension $L/\mathbf Q_p$, with maximal ideal $\mathfrak p=(\pi_L)$ and residue field $k$. We impose the usual hypotheses guaranteeing that the weight-two Eisenstein series
\[
  E_2(\theta,1)
\]
is well defined at level $\Gamma_1(N)$ and that the Eisenstein maximal ideal under consideration is non-exceptional; in particular, assume $\theta(p)-1\in\cO^\times$. We keep the congruence and cohomological setup of \cite[\S2]{RiRo}, including the eigensystem at the bad primes, and we let
\[
  f\in S_2(N,\theta;\mathcal O)
\]
be a normalized cuspidal eigenform satisfying
\begin{equation}
  a_\ell(f)\equiv \theta(\ell)+\ell
  \pmod{\mathfrak p^t}
  \qquad (\ell\nmid Np),
  \label{eq:eisenstein-congruence}
\end{equation}
for some $t\geq 1$. We take $t$ to be the exact congruence exponent, that is, the system of Hecke eigenvalues of $f$ factors through the Eisenstein quotient modulo $\mathfrak p^t$, but not modulo $\mathfrak p^{t+1}$.

Write $I_f^*$ for the kernel of the contragredient Hecke eigensystem attached to $f$.  Set
\[
  \overline Y=X_1(N)_{\overline{\mathbf Q}},
  \qquad
  Y=Y_1(N)_{\overline{\mathbf Q}},
  \qquad
  C=\overline Y\smallsetminus Y,
\]
and let $\imath:Y\hookrightarrow \overline Y$ denote the open immersion.  The two natural integral realizations of the representation attached to $f$ are the following, where the subscripts $X$ and $Y$ refer to the compactified curve $\overline Y=X_1(N)_{\overline\Q}$ and to the open curve $Y$, respectively:
\begin{align}
  T_{f,X}
   &:=H^1_{\mathrm{\acute et}}
      \bigl(\overline Y,\mathcal O(1)\bigr)/I_f^*,
  &
  T_{f,Y}
   &:=H^1_{\mathrm{\acute et}}
      \bigl(Y,\mathcal O(1)\bigr)/I_f^*.
  \label{eq:def-two-lattices}
\end{align}
After tensoring with $L$, both modules realize the same two-dimensional representation:
\[
  V_f:=T_{f,X}\otimes_{\mathcal O}L
      \simeq T_{f,Y}\otimes_{\mathcal O}L.
\]
Observe that these Hecke quotients may have $\cO$-torsion, and that their integral structures may differ even though their rational realizations agree. Whenever we use an actual basis below, we explicitly restrict to the case of torsion-free lattices.

\subsection{The boundary quotient}
\label{subsec:boundary-quotient}

The full Hecke eigensystem, together with the saturation statement of \cite[Proposition~2.1]{RiRo}, allows us to pass the localization sequence to these quotients. In particular, restriction induces an exact sequence
\begin{equation}
  0\longrightarrow T_{f,X}
   \longrightarrow T_{f,Y}
   \longrightarrow
   \operatorname{Div}^0(C;\mathcal O)/I_f^*
   \longrightarrow 0.
  \label{eq:localization-f-lattice}
\end{equation}
The boundary module is Eisenstein: away from $Np$, the dual Hecke operator $T_\ell^*$ acts on it through $1+\ell\langle\ell\rangle^*$, and on the eigenspace cut out by $I_f^*$, where $\langle\ell\rangle^*$ acts by $\theta^{-1}(\ell)$, this is the eigenvalue $\theta^{-1}(\ell)(\theta(\ell)+\ell)$ of the contragredient eigensystem, compatibly with \eqref{eq:eisenstein-congruence}. The Galois action on this eigenspace is given by $\theta$. Since $t$ is the exact congruence exponent, we obtain the following description.

\begin{proposition}
\label{prop:two-lattices-quotient}
There is an exact sequence of $\mathcal O[G_{\mathbf Q}]$-modules
\begin{equation}
  0\longrightarrow T_{f,X}
   \longrightarrow T_{f,Y}
   \xrightarrow{\,\pi_{\mathrm{Eis}}\,}
   \mathcal O/\mathfrak p^t(\theta)
   \longrightarrow 0.
  \label{eq:two-lattices-quotient}
\end{equation}
The final identification is canonical only up to multiplication by an element of $(\mathcal O/\mathfrak p^t)^\times$.
\end{proposition}

\begin{proof}
The localization sequence gives an exact sequence before taking the $f$-isotypical quotient, and the exactness on the left after taking this quotient is precisely the saturation statement
\[
 I_f^*H^1(Y,\mathcal O(1))\cap H^1(\overline Y,\mathcal O(1))
   =I_f^*H^1(\overline Y,\mathcal O(1));
\]
here we use the saturation statement in \cite[Proposition~2.1]{RiRo}, under its running hypotheses, to obtain \eqref{eq:localization-f-lattice}. Moreover, the Hecke action on the cusps factors through the Eisenstein quotient, and the relevant $\theta$-isotypical boundary module is cyclic. Since $t$ is the exact congruence exponent, this module has length $t$ and is therefore isomorphic to $\mathcal O/\mathfrak p^t$.  Finally, the standard calculation of the Galois action on the cusps identifies its character with $\theta$; see \cite[Theorem~1.3.1]{St-book}.
\end{proof}

Two a priori different integers appear here. The congruence exponent $t$ measures the length of the quotient in \eqref{eq:two-lattices-quotient}, while the valuation of a generalized Bernoulli number---for example \(v_p(B_2(\overline\theta))\) in the present normalization---measures the size of an Eisenstein congruence module. Identifying these two quantities requires an additional numerical criterion, typically furnished by an Eisenstein-ideal or a Gorenstein theorem, and it is not a formal consequence of the existence of the congruence.

\subsection{Global and ordinary filtrations}
\label{subsec:two-filtrations}

For an $\mathcal O$-module $T$, write \(\overline T=T/\mathfrak p^tT\).  The Eisenstein congruence endows \(\overline T_{f,X}\) with a global filtration.  With the Tate twist used in \eqref{eq:def-two-lattices}, it takes the form
\begin{equation}
  0\longrightarrow \overline T_{f,X}^{+}
   \longrightarrow \overline T_{f,X}
   \longrightarrow \overline T_{f,X}^{-}
   \longrightarrow 0,
  \label{eq:global-filtration-X}
\end{equation}
where the signs refer to the action of complex conjugation. The apparent reversal of signs with respect to some formulations of Sharifi's conjecture is due to the Tate twist in our convention; compare \cite[Theorem~4.3]{sharifi12} and \cite[\S\,6.3.1, \S\,7.1.11]{FK}.

On the other hand, ordinarity gives, for $*\in\{X,Y\}$, an exact sequence of $\mathcal O[G_{\mathbf Q_p}]$-modules
\begin{equation}
  0\longrightarrow T_{f,*}^{\mathrm{sub}}
   \longrightarrow T_{f,*}
   \longrightarrow T_{f,*}^{\mathrm{quo}}
   \longrightarrow 0,
  \label{eq:ordinary-filtration-XY}
\end{equation}
where the quotient is unramified. In the weight-two normalization used here, $G_{\mathbf Q_p}$ acts on the submodule through the cyclotomic character, up to the unramified character dictated by the chosen $p$-stabilization.  The global sequence \eqref{eq:global-filtration-X} splits upon restriction to $G_{\mathbf Q_p}$, and the two filtrations are transverse in the precise sense that
\begin{equation}
  \overline T_{f,X}^{+}
    \simeq \overline T_{f,X}^{\mathrm{quo}},
  \qquad
  \overline T_{f,X}^{\mathrm{sub}}
    \simeq \overline T_{f,X}^{-}.
  \label{eq:global-local-identifications}
\end{equation}
In particular, the ordinary submodule and the corresponding negative global quotient are free of rank one over $\mathcal O/\mathfrak p^t$. The ordinary quotient, however, may still contain torsion inherited from the Hecke quotient, and we come back to this issue in \S\ref{subsec:torsion-free-ordinary-quotient}.

Since the quotient in \eqref{eq:two-lattices-quotient} is unramified at $p$, the map between the two lattices induces an isomorphism on the ordinary submodules:
\begin{equation}
  T_{f,X}^{\mathrm{sub}}=T_{f,Y}^{\mathrm{sub}}.
  \label{eq:same-ordinary-submodule}
\end{equation}
Hence the passage from $\overline Y$ to $Y$ only enlarges the ordinary quotient, and this is the local manifestation of the global change of extension that we describe next.

\subsection{Changing the lattice reverses the extension}
\label{subsec:lattice-switch}

In order to display the change of extension in terms of matrices, suppose in this paragraph that the two Hecke quotients are free lattices of rank two and that their ordinary filtrations are saturated. This is the torsion-free model of the previous construction, and the general quotient modules are treated in the next subsection. Choose an $\mathcal O$-basis $\{u,v\}$ of $T_{f,X}$ such that, modulo $\mathfrak p^t$, the class of $u$ generates $\overline T_{f,X}^{\mathrm{sub}}$ and the class of $v$ generates $\overline T_{f,X}^{+}$.  With respect to this basis, the Galois action has the shape
\begin{equation}
  \rho_{T_{f,X}}(\sigma)=
  \begin{pmatrix}
    A(\sigma)&\pi_L^t B(\sigma)\\
    C(\sigma)&D(\sigma)
  \end{pmatrix},
  \qquad \sigma\in G_{\mathbf Q}.
  \label{eq:matrix-X}
\end{equation}
The divisibility of the upper-right entry expresses the stability of the line generated by $v$ modulo $\mathfrak p^t$, and the local splitting at $p$ says, in addition, that $C|_{G_{\mathbf Q_p}}$ is divisible by $\pi_L^t$.

Proposition~\ref{prop:two-lattices-quotient} and \eqref{eq:same-ordinary-submodule} show that \(\{u,\pi_L^{-t}v\}\) is a basis of $T_{f,Y}$.  In this new basis the same rational representation is described by
\begin{equation}
  \rho_{T_{f,Y}}(\sigma)=
  \begin{pmatrix}
    A(\sigma)&B(\sigma)\\
    \pi_L^t C(\sigma)&D(\sigma)
  \end{pmatrix}.
  \label{eq:matrix-Y}
\end{equation}
Consequently, the line generated by $u$ is now globally stable modulo $\mathfrak p^t$. Since the common semisimplification is
\[
  \overline T_{f,X}^{\mathrm{ss}}
  \simeq \overline T_{f,Y}^{\mathrm{ss}}
  \simeq
  \mathcal O/\mathfrak p^t(1)
  \oplus \mathcal O/\mathfrak p^t(\theta),
\]
we obtain an exact sequence
\begin{equation}
  0\longrightarrow \mathcal O/\mathfrak p^t(1)
   \longrightarrow \overline T_{f,Y}
   \xrightarrow{\,\overline\pi_{\mathrm{Eis}}\,}
   \mathcal O/\mathfrak p^t(\theta)
   \longrightarrow 0.
  \label{eq:residual-extension-Y}
\end{equation}
The comparison between \eqref{eq:matrix-X} and \eqref{eq:matrix-Y} is the central linear-algebra observation of this section: enlarging one basis vector by $\pi_L^{-t}$ moves the factor $\pi_L^t$ from one off-diagonal entry to the other. In this way, it changes which residual character is a subrepresentation and which one is a quotient, while the ambient rational representation stays the same.

\subsection{The torsion-free ordinary quotient}
\label{subsec:torsion-free-ordinary-quotient}

Recall that the ordinary quotients in \eqref{eq:ordinary-filtration-XY} need not be free over $\mathcal O$. Let $\Sigma_X$ and $\Sigma_Y$ denote their respective torsion submodules, and put
\[
  T_{f,Y,\circ}^{\mathrm{quo}}
  :=T_{f,Y}^{\mathrm{quo}}/\Sigma_Y.
\]
This is a free $\mathcal O$-module of rank one, and if $\psi_f$ denotes the unramified character giving the action of $G_{\mathbf Q_p}$ on the ordinary quotient, then
\begin{equation}
  T_{f,Y,\circ}^{\mathrm{quo}}
  \simeq \mathcal O(\psi_f),
  \qquad
  T_{f,Y,\circ}^{\mathrm{quo}}/\mathfrak p^t
  \simeq \mathcal O/\mathfrak p^t(\theta).
  \label{eq:free-ordinary-quotient}
\end{equation}
The second isomorphism, like the last arrow in \eqref{eq:two-lattices-quotient}, is well defined only up to a unit.

Under the weak Gorenstein hypothesis used in \cite{RiRo}---equivalently, in the present formulation, under the hypothesis that $\mathcal O/\mathfrak p^t(\theta)$ does not occur as a quotient of the torsion discrepancy $\Sigma_Y/\Sigma_X$---the Eisenstein projection kills $\Sigma_Y$, and therefore it factors as
\begin{equation}
  T_{f,Y}
   \longrightarrow T_{f,Y,\circ}^{\mathrm{quo}}
   \longrightarrow
   T_{f,Y,\circ}^{\mathrm{quo}}/\mathfrak p^t
   \xrightarrow{\ \sim\ }
   \mathcal O/\mathfrak p^t(\theta).
  \label{eq:factorization-eisenstein-projection}
\end{equation}
This factorization is what makes the quotient compatible with the ordinary local condition, and hence suitable for Euler-system arguments.

\subsection{Relation with Sharifi's conjecture}
\label{subsec:relation-sharifi}

The compact and the open curves thus give two integral models of the same rational representation. The Eisenstein congruence appears in the extension carried by $\overline T_{f,X}$, and enlarging the lattice moves the extension data to the other off-diagonal entry, as in \eqref{eq:matrix-X}--\eqref{eq:matrix-Y}; this change of lattice is controlled by the boundary quotient. The map $\pi_{\mathrm{Eis}}$ then extracts the Eisenstein character through the ordinary quotient. In this framework, comparing the projected Euler-system class with the cyclotomic class predicted by Sharifi's conjecture is the congruence-theoretic counterpart of the degeneration in families that we consider in the next section.

\section{Eisenstein degeneration in families}
\label{sec:eisenstein-degeneration}

We now pass from a fixed congruence to a one-parameter deformation. The basic strategy remains the same: we place an Eisenstein specialization inside a generically cuspidal family, construct the relevant Galois representation or Euler-system class away from the Eisenstein point, and study its leading term at that point. The geometry differs in the two cases we consider. A critical-slope Eisenstein point of weight at least two lies on a generically cuspidal Coleman family; an irregular weight-one point is reached along an ordinary cuspidal branch, which in the quadratic setting of \cite{BDV} is a CM Hida family.

\subsection{Critical-slope Eisenstein points}
\label{subsec:critical-slope-eisenstein}

Let $r\geq 0$, and let $\psi$ and $\tau$ be primitive Dirichlet characters of conductors prime to $p$ satisfying \(\psi\tau(-1)=(-1)^r\).  Consider the Eisenstein series
\[
  E_{r+2}(\psi,\tau)
\]
and assume that its level is prime to $p$.  Its Hecke polynomial at $p$ has roots
\[
  \alpha=\psi(p),
  \qquad
  \beta=p^{r+1}\tau(p).
\]
We write
\[
  f_\beta=E_{r+2}^{\mathrm{crit}}(\psi,\tau)
\]
for the $p$-stabilization with $U_p$-eigenvalue $\beta$, which is the critical-slope refinement since $v_p(\beta)=r+1$.

Following \cite{LR1}, and the study of the geometry at critical points of \cite{BC06}, we impose two hypotheses. First, $f_\beta$ is \emph{$p$-decent}: either $r>0$, or $r=0$ and, for every prime $\ell\mid Np$, the conductor of $\psi/\tau$ is divisible by $\ell$ or $(\psi/\tau)(\ell)\neq 1$.  Second, $f_\beta$ is \emph{non-critical}, which means that its generalized eigenspace in the space of overconvergent modular forms is one-dimensional. In this Eisenstein setting, the latter condition is equivalent to
\begin{equation}
  L_p(\psi^{-1}\tau,r+1)\neq 0.
  \label{eq:noncritical-p-adic-L-value}
\end{equation}

Here $L_p(\psi^{-1}\tau,s)$ denotes the Kubota--Leopoldt $p$-adic $L$-function of $\psi^{-1}\tau$ in the normalization of \cite{LR1}. Since $\psi^{-1}\tau$ has parity $(-1)^r$, for odd $r$ this normalization includes the appropriate power of the Teichm\"uller character, so this function should not be confused with the function $L_p^{\mathrm{KL}}$ of \S\ref{subsec:cyclotomic-conventions}, which we only use for even characters. These assumptions provide the local geometry required for the degeneration arguments.

\begin{theorem}[Bella\"iche; Bella\"iche--Chenevier]
\label{thm:etale-critical-eisenstein}
Assume that $f_\beta$ is $p$-decent and non-critical.  Then the cuspidal eigencurve is smooth at the point defined by $f_\beta$, and the weight map is \'etale there.  Consequently, after shrinking around this point, there is a closed weight disc $U$ and a unique Coleman family
\[
  \mathcal F=\sum_{n\geq 1}a_n(\mathcal F)q^n
  \]
over $U$ specializing to $f_\beta$ and cuspidal at every classical point of $U$ distinct from the Eisenstein point.
\end{theorem}

This statement is a particular case of the results on the eigencurve of \cite{bellaiche11a,LR1}, and the formulation above is the one used in \cite{LR1}. Choose a weight coordinate $X$ on $U$ which vanishes at $f_\beta$, and write $\mathcal O(U)$ for the ring of rigid-analytic functions on $U$.  Then the punctured disc $U^\times=U\smallsetminus\{X=0\}$ is entirely cuspidal, and it is over $U^\times$ that one first constructs the family of Galois representations and the associated Euler-system classes.

The special fibre is reducible, and this reducibility is measured to first order by the weight parameter. More precisely, the reducibility ideal of the two-dimensional pseudocharacter over the local ring of the eigencurve is the maximal ideal generated by $X$.  As a consequence, the two natural cohomological realizations of the family have different specializations at $X=0$. In order to compare with \cite{LR1,PR25}, write $V^c(f_\beta)$ and $V(f_\beta)$ for the modules used there. Their generic cuspidal normalization is dual to our arithmetic-Frobenius normalization, and those papers also attach Dirichlet characters to geometric Frobenius, so that their characters $\psi$ and $\tau$ become $\psi^{-1}$ and $\tau^{-1}$ in our notation. Hence the exact sequences, translated into the character convention of \S\ref{sec:preliminaries}, are
\begin{align}
  0&\longrightarrow \tau^{-1}\varepsilon^{-1-r}
    \longrightarrow V^c(f_\beta)
    \longrightarrow \psi^{-1}
    \longrightarrow 0,
  \label{eq:critical-compact-support-extension}\\
  0&\longrightarrow \psi^{-1}
    \longrightarrow V(f_\beta)
    \longrightarrow \tau^{-1}\varepsilon^{-1-r}
    \longrightarrow 0.
  \label{eq:critical-usual-extension}
\end{align}
Here $\varepsilon$ stands for the $p$-adic cyclotomic character. The natural map \(V^c(f_\beta)\to V(f_\beta)\) has one-dimensional image, so compactly supported and usual cohomology select opposite extensions of the same two characters. This is the analogue, in characteristic zero and in families, of the switch between $T_{f,X}$ and $T_{f,Y}$ in \S\ref{subsec:lattice-switch}.

\subsection{Leading terms of cuspidal Euler systems}
\label{subsec:leading-terms-critical}

Let $\mathcal V$ denote one of the coherent Galois modules attached to $\mathcal F$ over $U$, and let
\[
  \mathbf z_{\mathcal F}\in
  H^1\bigl(\mathbf Q,\mathcal V|_{U^\times}\bigr)
\]
be a class obtained by varying a cuspidal Euler system. Although $\mathbf z_{\mathcal F}$ is analytic on the punctured disc, it need not lie in the integral lattice over the whole disc $U$, and after choosing a lattice one may only have
\begin{equation}
  \mathbf z_{\mathcal F}
  \in X^{-m}H^1(\mathbf Q,\mathcal V)
  \qquad (m\geq 0).
  \label{eq:pole-order-euler-family}
\end{equation}
The smallest such $m$ is the order of the possible pole. Multiplication by $X^m$ then produces an integral class, whose specialization at $X=0$ is the leading term of the degeneration.

It is convenient to distinguish two phenomena. First, a pole may come from the mismatch between the natural cohomological lattices at the reducible point. Second, even when the family class is integral, its specialization may vanish for representation-theoretic reasons, and one must then divide by an additional power of $X$ before specializing. Explicit reciprocity laws relate these orders with the order of the associated $p$-adic $L$-function, although recovering the order of the class itself also requires some control of the kernel of the regulator and of the period lattice. The output is typically a cyclotomic unit or a circular-unit class, multiplied by the appropriate $p$-adic $L$-value or logarithmic factor.

\subsection{Irregular Eisenstein points of weight one}
\label{subsec:irregular-weight-one}

The weight-one situation has a different geometry. Let \(E_1(\chi_1,\chi_2)\) be a weight-one Eisenstein series, and suppose that its two $U_p$-roots coincide:
\[
  \chi_1(p)=\chi_2(p).
\]
The resulting refinement is $p$-irregular. In the quadratic case relevant to \cite{BDV}, let $K$ be an imaginary quadratic field in which $p$ splits and let $\varepsilon_K$ be its quadratic character.  Then
\[
  f_K=E_1(1,\varepsilon_K)(q)
      -E_1(1,\varepsilon_K)(q^p)
\]
is the unique $p$-stabilization, since its Hecke polynomial at $p$ is \((T-1)^2\).

Recall that the full eigencurve at this point may be singular and non-Gorenstein (see \cite{pozzi-thesis}), so for the Euler-system comparison we use its cuspidal branch. Since $\varepsilon_K$ is primitive and odd, with $\varepsilon_K(p)=1$, \cite[Theorem~A(i)]{BeDiPo} shows that this branch is \'etale over weight space. Therefore a sufficiently small weight disc carries an ordinary Hida family $\mathbf f_K$ through $f_K$, which is unique up to conjugation and has CM by $K$.

For the CM cohomological realization used here, the special representation itself splits, by \cite[Proposition~2.2]{BDV}:
\begin{equation}
  V_{f_K}\simeq \mathbf 1\oplus\varepsilon_K.
  \label{eq:weight-one-split-representation}
\end{equation}
This is an essential difference with respect to the opposite non-split extensions \eqref{eq:critical-compact-support-extension} and \eqref{eq:critical-usual-extension}: the repeated $U_p$-root does not lead to the geometry of the critical-slope case, and instead the CM Hida family provides a distinguished ordinary direction through a point at which the ambient eigencurve may have several components.

The method of \cite{BDV} places a Beilinson--Flach class in a family involving $\mathbf f_K$, and then specializes it at the irregular weight-one point. Using the decomposition \eqref{eq:weight-one-split-representation}, the specialized class has two projections, corresponding to the trivial and to the quadratic character, and the explicit reciprocity laws identify these projections with Kato-type cyclotomic classes, up to the expected $p$-adic $L$-values and logarithmic factors. Hence one again compares a generically cuspidal Euler system with rank-one Eisenstein data, but now the comparison is organized by the CM decomposition, and not by a pair of opposite non-split extensions.

The CM pair $(1,\varepsilon_K)$, which we use in Section~\ref{sec:double-degeneration}, lies outside the common-conductor setting of Theorem~A. Moreover, its even constituent is trivial and its ordinary eigenvalue is $1$, so the non-exceptional comparison with circular units of Section~\ref{sec:circular-kato} does not apply: in this case the first Dirichlet factor is a zeta function, and the local comparison requires the normalizations of the exceptional-zero setting.

\subsection{Comparison of the two degenerations}
\label{subsec:comparison-degenerations}

Both constructions place an Eisenstein point in a generically cuspidal family, but their cohomological specializations are different. At critical slope, compactly supported and usual cohomology give opposite extensions, and the weight parameter measures the passage between them. At irregular weight one, instead, the repeated Hecke root leads to a distinguished ordinary direction in a split realization, as in the CM case of \cite{BDV}.

In the comparisons of Sections~\ref{sec:circular-kato}--\ref{sec:weight-one-factorization}, the periods are fixed by means of a compact cohomological quotient of rank one. The resulting analytic factorization can then be compared with the regulator of a family class, provided that both the class and its regulator specialize in the chosen realization, and this is precisely where the difference between compact and open cohomology enters into the Euler-system argument.

\section{Circular units and Kato elements}\label{sec:circular-kato}

In this section we compare the weight-one specialization of Kato's Euler system with the Euler system of circular units. Following the critical-slope argument of \cite{PR25}, the idea is to use the two explicit reciprocity laws to determine the coefficient of the surviving circular-unit class. In contrast with the critical-slope case, however, the analytic factorization is not available a priori, and it has to be established separately once the period normalization is fixed. For this reason, we first obtain the cohomological comparison in terms of a quotient of $p$-adic $L$-functions, and then Section~\ref{sec:weight-one-factorization} evaluates its analytic coefficient in the common-conductor setting.

Let
\[
    f=E_1(\chi_1,\chi_2)
\]
be an irregular Eisenstein series of weight one and tame level $N=N_1N_2$, where $\chi_1$ and $\chi_2$ are Dirichlet characters of conductors prime to $p$ and
\[
    \chi_1(p)=\chi_2(p).
\]
Since $f$ has weight one, its nebentype $\chi_1\chi_2$ is odd. We order the characters in such a way that $\chi_1$ is even and $\chi_2$ is odd, and assume throughout that $\chi_1$ is non-trivial. Let $\hf$ be the generically cuspidal ordinary family passing through $f$, and let $x_f$ denote the corresponding point of its weight space. We write $A_U$ for the coefficient ring of $\hf$ on a sufficiently small neighbourhood $U$ of $x_f$, and fix a local parameter $X\in A_U$ vanishing at $x_f$.

In the arithmetic-Frobenius convention of Section~\ref{sec:preliminaries}, the semisimple Artin representation at weight one is
\[
 V_f^{\mathrm{ss}}=L(\chi_1)\oplus L(\chi_2).
\]
Here $L(\chi_i)$ is a one-dimensional $L$-vector space on which $\sigma\in G_{\Q}$ acts by $\chi_i(\sigma)$, and there is no Tate twist in this formula. In contrast, the Kummer map for a number field $F$ takes a unit to $H^1(F,L(1))$, and hence projecting a Kummer class to a Dirichlet character naturally leads to coefficients $L(\chi_i^{-1})(1)$, where the inverse character compensates for the character projector on the units.

Accordingly, we use the \emph{Kummer normalization}
\[
 \mathbb W_{\hf}:=\mathbb V_{\hf}^{*}(1),\qquad
 W_f^{\mathrm{ss}}=L(\chi_1^{-1})(1)\oplus L(\chi_2^{-1})(1).
\]
The coefficient line $L(\chi_i^{-1})(1)$ has Galois action $\chi_i^{-1}\varepsilon$, and in order to project the specialized class to it we need the splitting of the chosen family lattice specified in Hypothesis~\ref{hyp:kato-specialization}.

\begin{remark}[Comparison with \cite{PR25}]
The notation in \cite{PR25} uses a different duality convention: its module $V(\mathbf f)^*$ agrees generically with our $\mathbb V_{\hf}$, whereas the comparison with units uses $\mathbb W_{\hf}=\mathbb V_{\hf}^*(1)$. At a split weight-one fibre, two-dimensional duality identifies the latter with $V_f\otimes L(\epsilon_f^{-1})(1)$, where $\epsilon_f=\chi_1\chi_2$. Thus, to use the class of \cite[Definition~4.6]{PR25} here, we transport it by this nebentype and Tate twist, together with the corresponding change of cyclotomic variable. The compatibility of this transport with the reciprocity law is part of Hypothesis~\ref{hyp:kato-specialization}.
\end{remark}

Write $\mathbb W_{\hf}^{-}$ for the chosen rank-one ordinary quotient in this normalization; observe that it is different from the unramified quotient of $\mathbb V_{\hf}$ of Section~\ref{subsec:deligne-representations}. We require that
\begin{equation}\label{eq:ordinary-quotient-weight-one}
    \mathbb W_{\hf}^{-}/X\mathbb W_{\hf}^{-}
        \simeq L(\chi_1^{-1})(1).
\end{equation}
At weight one the two characters are globally distinct, since their parities differ, but their restrictions to $G_{\Q_p}$ agree, because $\chi_1(p)=\chi_2(p)$ and both are unramified at $p$. Consequently, the character of the quotient does not determine the two coordinates of the quotient map. On the even cyclotomic component used below, the odd character component of the global class vanishes. We therefore need only the first coordinate of the quotient and the specialization of a regularly normalized rank-one regulator.

\subsection{Circular units and cyclotomic cohomology}
\label{subsec:circular-units}

We keep the cyclotomic algebras and components of Section~\ref{subsec:cyclotomic-conventions}, and we assume that $L$ contains the values of the characters. Let $\mathbf j:G_\Q\to\Lambda_{\mathrm{cyc}}^\times$ be the tautological character obtained from the cyclotomic character and the inclusion $\Gamma\hookrightarrow\Lambda_{\mathrm{cyc}}^\times$. We work on a fixed even $\omega^a$-component; the displayed analytic formulas correspond to $a=0$, and the other components are obtained as in Corollary~\ref{prop:cyclotomic-period-comparison}.

In order to use the analytic exponent $j$ corresponding to the complex argument $1+j$, we pull back the usual inverse cyclotomic deformation by the involution $\iota:[z]\mapsto[z^{-1}]$, which we apply simultaneously to the coefficients, the classes and the regulators. The resulting deformation has twist $(\mathbf j)$ and specializes at an algebraic exponent $j$ to $L(\chi_1^{-1})(1+j)$. Equivalently, the Tate twist may be included in the cyclotomic notation:
\[
 L(\chi_1^{-1})(1)\widehat\otimes_L\mathcal H_\Gamma(\mathbf j)
   =L(\chi_1^{-1})\widehat\otimes_L\mathcal H_\Gamma(1+\mathbf j).
\]

Consider the one-dimensional representation
\[
    T_{\chi_1}=\mathcal O(\chi_1^{-1})(1)
\]
and its cyclotomic deformation
\[
    \mathbb T_{\chi_1}
      =T_{\chi_1}\widehat\otimes_{\mathcal O}\Lambda_{\mathrm{cyc}}(\mathbf j).
\]
After inverting $p$ and extending scalars to $\mathcal H_\Gamma$, we put
\[
    \mathbb V_{\chi_1}
      =L(\chi_1^{-1})(1)
         \widehat\otimes_L\mathcal H_\Gamma(\mathbf j).
\]

The Euler system of circular units gives rise to a canonical class
\begin{equation}\label{eq:circular-unit-class}
    c(\chi_1)\in H^1(\mathbb Q,\mathbb V_{\chi_1}).
\end{equation}
Before extending scalars, it is obtained from the inverse limit of Kummer images below, followed by the involution that we have just introduced:
\[
    c(\chi_1)\in
    \varprojlim_n
       H^1\bigl(\mathbb Q(\mu_{p^n}),T_{\chi_1}\bigr).
\]
We refer the reader to \cite[\S1.1]{BCDDPR} for the details of the construction and the normalizations.

\begin{ass}\label{ass:no-local-exceptional-zero}
We assume that
\[
 H^0\bigl(\mathbb Q_p(\mu_{p^\infty}),T_{\chi_1}\bigr)=0.
\]
For the unramified character considered here, this is equivalent to $\chi_1(p)\neq1$, and we use this explicit condition throughout the local argument.
\end{ass}

Let $S$ contain $p$ and all tame conductor primes, and put
\[
 H^1_{\mathrm{Iw}}(\Z[1/S],T)
   :=\varprojlim_n H^1\bigl(\cO_{\Q(\mu_{p^n})}[1/S],T\bigr),
\]
with corestriction as transition maps, where cohomology stands for Galois cohomology unramified outside the primes above $S$. Recall that the idempotent $e_a$ of Section~\ref{subsec:cyclotomic-conventions} selects the fixed even cyclotomic component. A parity argument gives the following vanishing result. Note that, since $\chi_1(p)=\chi_2(p)$, the condition $\chi_2(p)\ne1$ in it is equivalent to Assumption~\ref{ass:no-local-exceptional-zero}, and that the class $\kappa(f)_{\chi_2}$ appearing in its last assertion is introduced in \S\ref{subsec:kato-weight-one}.

\begin{lemma}[The odd character component]\label{lem:odd-component-zero}
Put $T_{\chi_2}=\cO(\chi_2^{-1})(1)$ and assume $\chi_2(p)\ne1$. On every fixed even cyclotomic component,
\[
 e_a H^1_{\mathrm{Iw}}(\Z[1/S],T_{\chi_2})[1/p]=0.
\]
Consequently, a regular global specialization with split coefficients $W_f$ has $\kappa(f)_{\chi_2}=0$ on that component.
\end{lemma}
\begin{proof}
We use Kummer theory in the non-exceptional range. Let $C$ be the conductor of $\chi_2$ and put $F_n=\Q(\mu_{Cp^n})$. Since the $p$-adic Tate module of the finite $S$-class group of $F_n$ is zero, the Kummer sequence identifies the finite-level $\Z_p(1)$-cohomology with the $p$-adic completion of its $S$-units. Tensoring with $\chi_2^{-1}$ selects the $\chi_2$-eigenspace of these units, and we take norm-compatible sequences before inverting $p$. Since $\chi_2$ is odd and the cyclotomic component is even, this is an odd character component of the units. The units of the CM field $F_n$ have no such component after removing the roots of unity, and the $p$-power roots of unity do not contribute either, since their tame action is trivial while $\chi_2$ is nontrivial.

It remains to consider the valuations at the primes in $S$. At the primes above $p$, the tame decomposition group acts on the character line through $\chi_2(p)\ne1$, so this part of the divisor module has zero $\chi_2$-projection. For a prime $\ell\ne p$, its decomposition group in the cyclotomic $\Z_p$-extension is open, since $\langle\ell\rangle$ has infinite order. Hence, in sufficiently high layers, the residue degree is $p$ at every step and the norm on the valuation module is multiplication by $p$, so the inverse limit over integral $p$-adic valuations is zero. As a consequence, every norm-compatible odd sequence of $S$-units has zero valuations, and it belongs to the unit component, which we have already seen to vanish. This proves the assertion, since the involution and the extension of scalars used in our convention preserve this vanishing.
\end{proof}

The same non-exceptional condition gives the following result.

\begin{lemma}\label{lem:local-iwasawa-free}
Under Assumption~\ref{ass:no-local-exceptional-zero}, the local Iwasawa cohomology group
\[
    H^1_{\mathrm{Iw}}(\mathbb Q_p,T_{\chi_1})
      :=\varprojlim_n
        H^1\bigl(\mathbb Q_p(\mu_{p^n}),T_{\chi_1}\bigr)
\]
is free of rank one over the chosen component of $\Lambda_{\mathrm{cyc}}[1/p]$. Consequently,
\[
    H^1(\mathbb Q_p,\mathbb V_{\chi_1})
      \simeq
    H^1_{\mathrm{Iw}}(\mathbb Q_p,T_{\chi_1})
       \otimes_{\Lambda_{\mathrm{cyc}},\iota}\mathcal H_\Gamma
\]
is a free $\mathcal H_\Gamma$-module of rank one.
\end{lemma}

\begin{proof}
The local Iwasawa Euler characteristic formula shows that the rank is one. Under the stated hypothesis the rational cyclotomic invariants vanish, and for the unramified finite-order character considered here the dual invariants also vanish when $\chi_1(p)\ne1$. Hence the standard structure theorem for local Iwasawa cohomology gives freeness over the rational Iwasawa algebra, that is, after inverting $p$. Extending scalars and applying $\iota$ gives the displayed isomorphism.
\end{proof}

There is a Coleman map
\begin{equation}\label{eq:coleman-circular}
    \operatorname{Col}_{\chi_1}:
       H^1(\mathbb Q_p,\mathbb V_{\chi_1})
          \longrightarrow \mathcal H_\Gamma
\end{equation}
whose specializations are given by Bloch--Kato logarithms and dual exponential maps, together with their interpolation factors. We use the measure-valued normalization of the Coleman map of \cite[Propositions~3.3--3.4]{RiRo}, translated by the involution above. Thus the circular-unit class and the Dieudonn\'e coordinate are fixed so that
\begin{equation}\label{eq:reciprocity-circular}
 \operatorname{Col}_{\chi_1}(\operatorname{res}_p c(\chi_1))
   =\mathcal F_{\chi_1}(j)
   :=L_p^{\mathrm{KL}}(\chi_1,1+j).
\end{equation}
The Euler factors in this identity are those of the standard Kubota--Leopoldt function, so there is no additional normalization factor. At a finite-order character $\rho$ of $\Gamma_1$, the same identity reads
\begin{equation}\label{eq:circular-finite-character-dictionary}
 \mathcal F_{\chi_1}(\rho)
     =L_p^{\mathrm{KL}}(\chi_1\rho^{-1},1).
\end{equation}
Indeed, writing $G_\beta(u^s-1)=L_p^{\mathrm{KL}}(\beta,1-s)$, twisting by $\rho^{-1}$ replaces the variable $1+T$ by $\rho(u)^{-1}(1+T)$, as in \cite[equations~(2.8)--(2.9)]{SW22}, and this is precisely the involution in our cohomological convention. Observe also that, since $\chi_1$ is nontrivial, $\mathcal F_{\chi_1}$ is a nonzero bounded measure.

\begin{corollary}\label{cor:coleman-injective}
The Coleman map \eqref{eq:coleman-circular} is injective.
\end{corollary}

\begin{proof}
By Lemma~\ref{lem:local-iwasawa-free}, its source is free of rank one over the integral domain $\mathcal H_\Gamma$, and the map is non-zero by \eqref{eq:reciprocity-circular}. Hence it is injective.
\end{proof}

To compare global classes through their regulators, we also need localization at $p$ to be injective. We work with the module obtained from the same Iwasawa cohomology:
\[
 \mathcal S_{\chi_1}
   :=e_aH^1_{\mathrm{Iw}}(\Z[1/S],T_{\chi_1})[1/p]
       \otimes_{\Lambda_L,\iota}\mathcal H_\Gamma
   \subset H^1(\Q,\mathbb V_{\chi_1}).
\]
Here and below, $\mathcal H_\Gamma$ refers to the fixed component. This module contains the circular-unit class and, under Hypothesis~\ref{hyp:kato-specialization}(i), the projected Kato specialization; the following proposition shows that a global class in it is determined by its Coleman image.

\begin{proposition}\label{ass:global-rank-one}
Assume $\chi_1$ is nontrivial and even, and $\chi_1(p)\ne1$. Then $\mathcal S_{\chi_1}$ is torsion-free of rank one, and
\[
 \operatorname{Col}_{\chi_1}\circ\operatorname{res}_p:
        \mathcal S_{\chi_1}\longrightarrow\mathcal H_\Gamma
\]
is injective.
\end{proposition}
\begin{proof}
Let $C$ be the conductor of $\chi_1$ and put $F_n=\Q(\mu_{Cp^n})$. Kummer theory identifies the global Iwasawa module with the corresponding character component of norm-compatible $S$-units. As in Lemma~\ref{lem:odd-component-zero}, the inverse limit of the valuation modules at primes away from $p$ is zero, while the $\chi_1$-part at $p$ vanishes because $\chi_1(p)\ne1$, so only units remain; moreover, since the tame character is nontrivial, the $p$-power roots of unity do not contribute either.

Leopoldt's theorem for the abelian fields $F_n$, due to Brumer \cite{Brumer67}, makes localization of these units into the product of the completions at $p$ injective after tensoring with $\Q_p$. In fact, this localization is already injective on $\cO_{F_n}^\times\otimes\Z_p$: its kernel is finite, hence contained in the $p$-power roots of unity, which inject into the local units. Therefore the injectivity is preserved when passing to the inverse limit. Taking norm-compatible systems and the character component, and using Shapiro's lemma to identify the semilocal target, gives an injection into $e_aH^1_{\mathrm{Iw}}(\Q_p,T_{\chi_1})[1/p]$. This target is free of rank one by Lemma~\ref{lem:local-iwasawa-free}, and localization remains injective under the flat extension to locally analytic distributions. It follows that $\mathcal S_{\chi_1}$ is torsion-free of rank at most one; the circular-unit class has nonzero regulator by \eqref{eq:reciprocity-circular}, so the rank is exactly one. Finally, composing localization with the injective Coleman map proves the assertion.
\end{proof}

\subsection{The analytic input to Kato reciprocity}
\label{subsec:factorization-kato}

The analytic identity used below involves three factors: a modular $p$-adic $L$-function varying cyclotomically, an auxiliary modular factor depending only on the weight, and an adjoint factor. Following \cite[\S4.3]{PR25}, we describe their arguments explicitly. We denote the auxiliary characters by $(\vartheta_1,\vartheta_2)$ to distinguish them from the constituents $(\chi_1,\chi_2)$ of $f$.

For the fixed cohomological periods $\gamma_{\hf}^{\pm}$, let
\[
 \mathscr L_p(\hf,\vartheta_i)
       \in A_U\widehat\otimes_L\mathcal H_\Gamma
\]
be the twisted Mazur--Kitagawa functions, and write $\mathscr L_p(\hf)$ for the untwisted one. Their argument $j$ follows the exponent convention of Section~\ref{subsec:circular-units}: at a classical weight $k\geq2$ they interpolate $L(f_k\otimes\vartheta_i,1+j)$ for $0\leq j\leq k-2$. In contrast, the Rankin variable $s$ below stands for the complex $L$-argument, so that a modular factor at complex argument $s$ is written $\mathscr L_p(\hf,\vartheta_i)(k,s-1)$.

The adjoint function $L_p(\operatorname{Ad}\hf)\in A_U$ compares the Petersson and the cohomological period normalizations, as recalled in Section~\ref{subsec:adjoint-p-adic-L}, and its possible zero at $f$ is precisely the reason why we keep it in the formula. Let $\mathbf E=\mathbf E(\vartheta_1,\vartheta_2)$ be the auxiliary ordinary Eisenstein family, of weight $\ell$. Artin formalism gives the two complex arguments $s$ and $s-\ell+1$ in its Rankin product with $f_k$; compare \cite[Theorem~3.4]{BD1} and \cite[Proposition~4.4]{PR25}.

Write $L_p^{\mathrm{HR}}(\hf,\mathbf E)(k,\ell,s)$ for the three-variable Hida--Rankin $p$-adic $L$-function, with $\hf$ as the dominant family and $s$ the complex Rankin argument. We use the Eisenstein-family normalization of \cite[\S4.3]{PR25}; for the general Hida--Rankin interpolation convention and its dependence on the dominant variable, see \cite[\S3.6]{Das}. The period and tame Euler factors entering the comparison with our Mazur--Kitagawa functions are retained below.

Restrict the Rankin variables to the surface
\begin{equation}\label{eq:critical-surface}
 s=s(k,\ell):=\frac{k-1}{2}+\ell.
\end{equation}
Then the two complex arguments are $(k-1)/2+\ell$ and $(k+1)/2$, which at $k=1$ become $\ell$ and $1$; their modular exponents are $s-1$ and $(k-1)/2$. Observe that along this surface the second modular factor depends only on the weight of $\hf$. At classical points we take $k$ odd, so that the two displayed arguments are integers. On a fixed weight disc, the half-weight is the analytic character $z\mapsto\exp((k-1)\log_p(z)/2)$ of $\Gamma_1$; the finite character components remain fixed. Classical points of this parity in the Rankin interpolation range form the interpolation set. 

\begin{proposition}\label{prop:hida-rankin-factorization}
With these conventions, there is a meromorphic period factor $\mathfrak d\in\operatorname{Frac}(A_U)$ depending only on $k$ such that
\begin{align}
 &L_p^{\mathrm{HR}}(\hf,\mathbf E)
       \bigl(k,\ell,s(k,\ell)\bigr)
       L_p(\operatorname{Ad}\hf)(k)\notag\\
 &\qquad=
 \mathfrak d(k)\,
 \mathscr L_p(\hf,\vartheta_1)
       \left(k,\frac{k-3}{2}+\ell\right)
 \mathscr L_p(\hf,\vartheta_2)
       \left(k,\frac{k-1}{2}\right),
 \label{eq:hida-rankin-factorization}
\end{align}
up to the explicit Euler factors at the primes dividing the tame levels.
\end{proposition}
\begin{proof}
At classical points in the Rankin interpolation range, the complex factorization is the product of the twists at $s$ and $s-\ell+1$, and the comparison of interpolation formulas in \cite[Proposition~4.4]{PR25}, following \cite[Theorem~3.4]{BD1}, accounts for the adjoint period and the tame factors. Restricting to \eqref{eq:critical-surface} and replacing each complex modular argument by its exponent, that is, subtracting $1$, we obtain the displayed identity on the classical interpolation set, and hence, by analytic continuation, as an identity of meromorphic functions in the family.
\end{proof}

In order to obtain the untwisted cyclotomic factor used below, we take $\vartheta_1=1$ and fix $\vartheta_2$ with the required parity. We then pull back the surface by $s=1+j$, equivalently $\ell=j+(3-k)/2$, and put
\[
 \mathscr L_p^{\mathrm{aux}}(\hf)(k)
   :=\mathscr L_p(\hf,\vartheta_2)
                         \left(k,\frac{k-1}{2}\right).
\]
The first factor is then $\mathscr L_p(\hf)(k,j)$, while the auxiliary factor only depends on $k$. Abbreviating the pulled-back Rankin function by the same notation, the formula becomes
\begin{equation}\label{eq:factorization-schematic}
 L_p^{\mathrm{HR}}(\hf,\mathbf E)L_p(\operatorname{Ad}\hf)
    =\mathscr E_N\,
        \mathscr L_p^{\mathrm{aux}}(\hf)\mathscr L_p(\hf),
\end{equation}
where $\mathscr E_N$ collects the tame factors and $\mathfrak d$. This is a generic identity of meromorphic functions, and the question of its regularity at $f$ is addressed below.

\begin{ass}[Simple-zero normalization]\label{ass:simple-adjoint-zero}
In the fixed cohomological period modules and bases of this subsection, $\mathscr E_N$ is regular at $x_f$, and
\begin{equation}\label{eq:two-vanishings}
 \operatorname{ord}_{X=0}L_p(\operatorname{Ad}\hf)=1,
 \qquad
 \operatorname{ord}_{X=0}\mathscr L_p^{\mathrm{aux}}(\hf)\geq1.
\end{equation}
\end{ass}

This assumption concerns the chosen analytic normalization, and it is used in Propositions~\ref{thm:kato-leading-term-weight-one} and~\ref{thm:kato-circular-comparison}, whose conductor hypotheses are more general than those of Section~\ref{sec:weight-one-factorization}. In the common-conductor setting, Corollary~\ref{cor:conditional-kato-circular-comparison} dispenses with it by treating the full reciprocity factor at once. Observe that changes of basis within a fixed period module preserve the orders of vanishing, while changing the module by a power of $X$ may modify them.

Under this assumption the quotient has a removable singularity, and we set
\begin{equation}\label{eq:def-Cf}
 \mathfrak C_f
  :=\left.
       \frac{\mathscr L_p^{\mathrm{aux}}(\hf)}
            {L_p(\operatorname{Ad}\hf)}
     \right|_{x_f}
   =\frac{
       \left.\dfrac{d}{dX}\mathscr L_p^{\mathrm{aux}}(\hf)\right|_{X=0}}
      {\left.\dfrac{d}{dX}L_p(\operatorname{Ad}\hf)\right|_{X=0}}.
\end{equation}
This ratio does not depend on the choice of the local parameter $X$, and it vanishes when the auxiliary function has order greater than one, a case which is allowed in the comparisons below.

\subsection{Kato classes and the weight-one specialization}
\label{subsec:kato-weight-one}

We now turn to the specialization of Kato classes at $f$. The splitting of the weight-one representation allows us to project a specialized class to its two character components and compare the $\chi_1$-component with circular units. We use the class with auxiliary characters of \cite[Definition~4.6]{PR25}, normalized as in Proposition~4.9 of loc.\,cit. Its reciprocity law involves the auxiliary Mazur--Kitagawa value and the adjoint period factor introduced above.

To carry out this comparison, fix a family of coefficient modules $\mathbb W_{\hf}$ with fibre
\[
 W_f\simeq L(\chi_1^{-1})(1)\oplus L(\chi_2^{-1})(1),
\]
and bases of its two character lines. In the common-conductor setting this realization is available: if $\mathbb V_{\mathbf f_0}^{\mathrm{cpt}}$ denotes the full rank-two localized compact module of the auxiliary family, \cite[Theorem~1.5]{SW21} gives the required splitting for
\[
 \mathbb W_{\mathbf f}^{\mathrm{tw}}
 :=\bigl(\mathbb V_{\mathbf f_0}^{\mathrm{cpt}}\otimes L(\chi_1)\bigr)^*(1).
\]
We take $\mathbb W_{\hf}=\mathbb W_{\mathbf f}^{\mathrm{tw}}$ in this case, and assume the existence of such a split realization for general conductors. Specialization and character projection then give maps on cohomology. To apply them to the Kato class, we need it to be regular in this realization at $X=0$.

With the auxiliary characters and period bases fixed above, write the transported family class as
\begin{equation}\label{eq:big-kato-class}
 \kappa(\hf)\in H^1\left(
   \Q,\mathbb W_{\hf}\widehat\otimes_L\mathcal H_\Gamma(\mathbf j)
 \right)[1/X].
\end{equation}
Here we have made the change of cyclotomic variable of Section~\ref{subsec:circular-units}. The notation $[1/X]$ allows the transported class to have a finite-order pole at $f$; the coefficient module itself is defined at $f$.

Choose a regular basis of the Dieudonn\'e line of the rank-one quotient $\mathbb W_{\hf}^{-}$. We write
\begin{equation}\label{eq:big-pr-map}
 \mathcal L_{\hf}^{-}:
 H^1\left(\Q_p,\mathbb W_{\hf}^{-}
            \widehat\otimes_L\mathcal H_\Gamma(\mathbf j)\right)
       \longrightarrow A_U\widehat\otimes_L\mathcal H_\Gamma
\end{equation}
for the regulator in this basis. We allow a fixed cyclotomic denominator, independent of $X$, and include projection to the ordinary quotient in the notation $\operatorname{loc}_p$. With this convention, the congruence and cohomological period factors occur on the analytic side of the reciprocity law. The following hypothesis specifies the regularity and compatibility needed at $f$.

\begin{hyp}[Specialization of the Kato class]\label{hyp:kato-specialization}
For the coefficient realization and normalizations fixed above, assume that:
\begin{enumerate}[label=\textup{(\roman*)}]
\item \emph{Regularity.} The generic class \eqref{eq:big-kato-class} is the image of a class in cyclotomic Iwasawa cohomology of the fixed coefficient realization, unramified outside a fixed finite set $S$, without inverting $X$.
\item \emph{Specialization of the regulator.} The map \eqref{eq:big-pr-map}, in a regular Dieudonn\'e basis, is defined on the chosen local cohomology module and specializes on every regular class to the rank-one Coleman map \eqref{eq:coleman-circular}, where the two character lines of the split fibre are identified locally, since $\chi_1$ and $\chi_2$ agree on $G_{\Q_p}$, and $\operatorname{Col}_{\chi_1}$ is applied to both projections through this identification. The ordinary quotient on the split local fibre is
\[
 \ell_0=a_f\,\operatorname{pr}_1+b_f\,\operatorname{pr}_2,
 \qquad a_f,b_f\in L,\quad a_f\ne0.
\]
\item \emph{Reciprocity and periods.} On the cuspidal locus, transport of the representation, cyclotomic variable and periods gives \eqref{eq:kato-reciprocity-family}. The factor $\mathscr E_N$ includes all weight-dependent period factors introduced by this transport, which may be meromorphic. The modular function uses the chosen cohomological period line; its comparison with compact periods depends only on weight.
\end{enumerate}
\end{hyp}

These conditions connect the chosen split realization with Kato's construction and its reciprocity law: (i) allows us to specialize the class, (ii) computes its Coleman image, and (iii) expresses that image in terms of modular $p$-adic $L$-functions.

Since the ordinary zeta construction takes values in open cohomology \cite[\S3.2]{FK}, condition~(i) asks that transport to the split realization preserve regularity. For Beilinson--Flach classes, the corresponding passage is supplied by the lift to parabolic coefficients in \cite[Lemma~2.4]{BDV}, which places the classes in the realization split by Proposition~2.2 of loc.\,cit. In our setting, the signed twisting map of Section~\ref{subsec:regular-compact-twisting} compares the period lines used in (iii), while regularity of the transported cohomology class is the content of (i).

In our even cyclotomic component, the regulator also provides a way to verify (i). Assume $\chi_1(p)\ne1$, and let $\mathcal M$ be the family Iwasawa cohomology module, unramified outside $S$, with the chosen coefficients and the flat scalar extensions used above. For a generic class in $\mathcal M[1/X]$, condition~(ii) gives the criterion
\begin{equation}\label{eq:kato-regularity-test}
 \kappa(\hf)\in\operatorname{im}\bigl(\mathcal M\longrightarrow\mathcal M[1/X]\bigr)
 \quad\Longleftrightarrow\quad
 \mathcal L_{\hf}^{-}(\operatorname{loc}_p\kappa(\hf))
 \text{ has no pole at }X=0.
\end{equation}
Indeed, multiplication by $X$ on the universal cyclotomic coefficient module is injective. The cohomology sequence therefore gives
\[
 \mathcal M/X\mathcal M\lhook\joinrel\longrightarrow
 \mathcal S_{\chi_1}\oplus\mathcal S_{\chi_2},
\]
where $\mathcal S_{\chi_2}$ is defined as in \S\ref{subsec:circular-units} with $\chi_2$ in place of $\chi_1$. The flat extension from bounded Iwasawa modules to distributions, and any fixed cyclotomic localization, preserve this injection. By Lemma~\ref{lem:odd-component-zero}, the second summand vanishes; on the first, the specialized regulator is injective by (ii) and Proposition~\ref{ass:global-rank-one}. Now suppose that $\kappa(\hf)$ has a pole of order $n>0$. A representative of $X^n\kappa(\hf)$ in $\mathcal M$ has nonzero specialization, and hence nonzero specialized regulator. If the regulator of $\kappa(\hf)$ were regular, that same value would be zero, since it is the specialization of $X^n\mathcal L_{\hf}^{-}(\operatorname{loc}_p\kappa(\hf))$. This proves the implication $\Leftarrow$ in \eqref{eq:kato-regularity-test}; the converse follows from (ii).

Once the class is regular, (ii) identifies its specialized regulator with the Coleman image of its $\chi_1$-component. The odd component vanishes by Lemma~\ref{lem:odd-component-zero}, so
\begin{equation}\label{eq:regular-kato-specialization}
 \left.\mathcal L_{\hf}^{-}(\operatorname{loc}_p\kappa(\hf))\right|_{x_f}
    =a_f\operatorname{Col}_{\chi_1}
            (\operatorname{res}_p\kappa(f)_{\chi_1}).
\end{equation}
The scalar $a_f$ records the ordinary quotient in the chosen character basis.

Hypothesis~\ref{hyp:kato-specialization}(i) allows us to work with the specialization of the original Kato class. If instead the class has a pole of order $n>0$ in $\mathcal M[1/X]$, the natural object to consider is the specialization of $X^n\kappa(\hf)$, followed by character projection. Under the same assumptions on the cyclotomic component and regulator as in \eqref{eq:kato-regularity-test}, its Coleman image is $a_f^{-1}$ times the coefficient of $X^{-n}$ in the scalar regulator. This would lead to a comparison for the leading term, with the coefficient realization and parameter $X$ fixed. The criterion itself does not determine whether a pole occurs. We treat the regular case here; Proposition~\ref{thm:kato-leading-term-weight-one} verifies regularity under Assumption~\ref{ass:simple-adjoint-zero} and the stated compatibility conditions.

To compute the regulator, we use the reciprocity law with the period normalization in (iii).

\begin{proposition}\label{prop:kato-reciprocity-family}
With the normalization in part~\textup{(iii)} of the specialization hypothesis, the family reciprocity identity is
\begin{equation}\label{eq:kato-reciprocity-family}
 \mathcal L_{\hf}^{-}\bigl(\operatorname{loc}_p\kappa(\hf)\bigr)
       L_p(\operatorname{Ad}\hf)
  =\mathscr E_N\,
       \mathscr L_p^{\mathrm{aux}}(\hf)\mathscr L_p(\hf).
\end{equation}
\end{proposition}
\begin{proof}
On the cuspidal locus, Kato's explicit reciprocity law and the period comparison give this identity in the form of \cite[Proposition~4.9]{PR25} (see also \cite{Kato,ochiai03,BB22}), and part~(iii) identifies its transported terms with our chosen functions and bases. The equality on the cuspidal locus then gives the meromorphic identity over $U$. When (i) and (ii) hold, the regulator on the left is regular at $f$.
\end{proof}

For a regular Kato class, specialization followed by projection gives
\begin{equation}\label{eq:projected-kato-class}
 \kappa(f)_{\chi_1}\in H^1(\Q,\mathbb V_{\chi_1}).
\end{equation}
It has the same coefficients as the circular-unit class $c(\chi_1)$ of Section~\ref{subsec:circular-units}. Under Assumption~\ref{ass:simple-adjoint-zero}, both the adjoint and auxiliary factors in \eqref{eq:kato-reciprocity-family} vanish at $f$. Cancelling their common factor $X$ gives a regular expression for the regulator. This proves regularity of the class by \eqref{eq:kato-regularity-test} and computes its specialization.

\begin{proposition}[General conductors]\label{thm:kato-leading-term-weight-one}
Assume parts~\textup{(ii)} and~\textup{(iii)} of Hypothesis~\ref{hyp:kato-specialization} and Assumptions~\ref{ass:no-local-exceptional-zero} and~\ref{ass:simple-adjoint-zero}. If the generic Kato class belongs to $\mathcal M[1/X]$ for the split realization fixed above, then it extends regularly to $X=0$, and
\begin{equation}\label{eq:kato-leading-term-weight-one}
 \operatorname{Col}_{\chi_1}
       \bigl(\operatorname{res}_p\kappa(f)_{\chi_1}\bigr)
   =a_f^{-1}\mathfrak C_f\,
      \mathscr E_N(f)\mathscr L_p(f),
\end{equation}
where $\mathscr L_p(f):=\mathscr L_p(\hf)|_{x_f}\in\mathcal H_\Gamma$ and $a_f$ is the quotient-coordinate scalar in part~\textup{(ii)} of the hypothesis.
\end{proposition}
\begin{proof}
Write $L_p(\operatorname{Ad}\hf)=Xa(X)$ and $\mathscr L_p^{\mathrm{aux}}(\hf)=Xb(X)$, with $a(0)\ne0$. The meromorphic reciprocity identity gives
\[
 \mathcal L_{\hf}^{-}(\operatorname{loc}_p\kappa(\hf))
   =\mathscr E_N\,\frac{b(X)}{a(X)}\,\mathscr L_p(\hf).
\]
Its right-hand side is regular: $\mathscr E_N$ is regular by Assumption~\ref{ass:simple-adjoint-zero}, and the modular function is analytic in the fixed cohomological periods of \S\ref{subsec:factorization-kato}. Criterion~\eqref{eq:kato-regularity-test} therefore gives regularity of the class. Evaluating at $X=0$ gives
\[
 \left.\mathcal L_{\hf}^{-}(\operatorname{loc}_p\kappa(\hf))\right|_{x_f}
   =\frac{b(0)}{a(0)}\mathscr E_N(f)\mathscr L_p(f).
\]
The quotient $b(0)/a(0)$ is $\mathfrak C_f$ by \eqref{eq:def-Cf}, and part~(ii) identifies the left-hand side with $a_f$ times the circular-unit regulator of the projected class.
\end{proof}

\begin{remark}\label{rem:not-derived-kato-class}
The class in Proposition~\ref{thm:kato-leading-term-weight-one} is obtained by specializing $\kappa(\hf)$ and projecting to the $\chi_1$-component. The ratio of derivatives
\[
 \mathfrak C_f
   =\frac{(\mathscr L_p^{\mathrm{aux}}(\hf))'(0)}
          {(L_p(\operatorname{Ad}\hf))'(0)}
\]
comes from cancelling the common zero of the two analytic factors. If the auxiliary factor has a zero of higher order, then $\mathfrak C_f=0$, as allowed in the proposition and the comparison below.
\end{remark}

\subsection{Comparison with circular units}
\label{subsec:kato-circular-comparison}

We can now compare \eqref{eq:reciprocity-circular} and \eqref{eq:kato-leading-term-weight-one}, keeping the actual circular-unit regulator $\mathcal F_{\chi_1}$.

\begin{proposition}[Comparison for general conductors]\label{thm:kato-circular-comparison}
Under the hypotheses of Proposition~\ref{thm:kato-leading-term-weight-one}, one has
\begin{equation}\label{eq:kato-circular-comparison}
    \kappa(f)_{\chi_1}
       =
       \frac{
          \mathfrak C_f\,\mathscr E_N(f)\,
          \mathscr L_p(f)
       }{
          a_f\,
          \mathcal F_{\chi_1}
       }
       \;c(\chi_1)
\end{equation}
in $\mathcal S_{\chi_1}\otimes_{\mathcal H_\Gamma} \mathcal K_\Gamma$.

If the quotient on the right-hand side belongs to $\mathcal H_\Gamma$, then the equality already holds in $\mathcal S_{\chi_1}$.
\end{proposition}

\begin{proof}
By the explicit reciprocity law for circular units,
\[
    \operatorname{Col}_{\chi_1}
       \bigl(\operatorname{res}_p c(\chi_1)\bigr)
       =\mathcal F_{\chi_1}.
\]
On the other hand, Proposition~\ref{thm:kato-leading-term-weight-one} gives
\[
    \operatorname{Col}_{\chi_1}
       \bigl(\operatorname{res}_p\kappa(f)_{\chi_1}\bigr)
       =
       a_f^{-1}\mathfrak C_f
       \mathscr E_N(f)\mathscr L_p(f).
\]
Hence both sides of \eqref{eq:kato-circular-comparison} have the same image under $\operatorname{Col}_{\chi_1}\circ\operatorname{res}_p$, and since this morphism is injective on $\mathcal S_{\chi_1}$ by Proposition~\ref{ass:global-rank-one}, the result follows.
\end{proof}

Formula \eqref{eq:kato-circular-comparison} is the comparison that we obtain under the stated hypotheses, without using any factorization into Dirichlet $p$-adic $L$-functions. Its coefficient contains the quotient
\begin{equation}\label{eq:unsimplified-quotient}
 \frac{\mathscr L_p(f)}{\mathcal F_{\chi_1}}
 =\frac{\mathscr L_p(f)}
        {L_p^{\mathrm{KL}}(\chi_1,1+j)}.
\end{equation}
Artin formalism suggests that the first Dirichlet factor should cancel, leaving a function attached to the other Eisenstein character. In Section~\ref{sec:weight-one-factorization} we make this statement precise in the common-conductor setting, separating the analytic computation from the period comparison needed to apply it here; the specialization input for the Kato class is kept, but the assumption on the leading term of the adjoint factor is removed.

Finally, observe that if $\chi_1(p)=1$ the local hypothesis used above fails. In that case, an exceptional-zero comparison would require a separate treatment of the local torsion and of the relevant $p$-unit or extended cohomology class, and we do not pursue this extension here.

\section{Artin factorization and compact periods}
\label{sec:weight-one-factorization}

We prove Theorem~A by evaluating the compact modular-symbol measure at the irregular weight-one point. Since there are no complex critical values at this weight, the factorization must come from the cohomological construction itself.

Write $f$ as a tame twist of an auxiliary form $f_0$ with one trivial constituent. We first construct the corresponding map of compact period modules over the weight ring. The localized Sharifi map then turns the twisted modular symbols of $f_0$ into cup products of cyclotomic units. Local reciprocity evaluates these cup products as a product of two sums, giving the two Kubota--Leopoldt factors. Their nonvanishing also shows that the twisting map is an isomorphism at weight one. After proving the factorization, we compare period modules and return to Kato classes; only this last application uses Hypothesis~\ref{hyp:kato-specialization}.

\subsection{The common-conductor setting}
\label{subsec:factorization-setting}

Throughout this section we assume that $p\geq5$, and we let $\chi_1$ and $\chi_2$ be primitive Dirichlet characters of the same conductor $M>1$, with $\chi_1$ even and $\chi_2$ odd. We further assume that
\begin{equation}\label{eq:factorization-conductor-hypotheses}
 p\nmid M\varphi(M),\qquad
 \chi_1(p)=\chi_2(p)=\alpha,\qquad
 \chi:=\chi_1\chi_2^{-1}\ \text{is primitive of conductor }M.
\end{equation}
We are interested in the form $f=E_1(\chi_1,\chi_2)^{\mathrm{ord}}$, whose minimal tame level is $M^2$. In order to relate it with the setting of Shih--Wang, put $\eta=\chi_1$ and consider the auxiliary form
\[
 f_0=E_1(1,\chi^{-1})^{\mathrm{ord}},\qquad\text{so that}\qquad
 f=f_0\otimes\eta.
\]
The constituents of $f_0$ are $1$ and $\chi^{-1}$, which places $f_0$ in the setting of the Sharifi comparison of \cite{SW21}; moreover, $\chi(p)=1$, so we are in their exceptional weight-one case. (In the dual-Hecke convention of loc.\,cit., the even character $\Theta=\chi\omega$ gives the eigenvalues $1+\chi^{-1}(\ell)$, which explains the inverse in the definition of $f_0$.) Twisting by the even character $\eta$ recovers the constituents $\chi_1$ and $\chi_2$ of $f$.

The hypotheses of \cite[Theorem~1.2]{SW21} are satisfied: $\Theta$ is primitive and even of conductor $Mp$, $p\geq5$, $p\nmid M\varphi(M)$, and $\chi(p)=1$. Similarly, the geometric input of \cite[Theorem~A(i)]{BeDiPo} requires precisely that $\chi^{-1}$ is primitive and odd and that $\chi(p)=1$. Since we work after localization in characteristic zero, no integral Gorenstein assumption is needed. Observe also that we do not assume $\alpha\ne1$ here: the twisting and reciprocity computations below only use $\chi(p)=1$, while the condition $\alpha\ne1$ is needed for the comparison with circular units, and we impose it explicitly in \S\ref{subsec:factorization-application-kato}. When $\alpha=1$, the second factor in Theorem~A has a trivial zero at $j=0$, which is compatible with the formula.

The ordinary eigenvalues of $f_0$ and $f$ are $1$ and $\alpha$, respectively. Let $\mathbf f_0$ be the cuspidal ordinary branch through $f_0$, and set $\mathbf f=\mathbf f_0\otimes\eta$. By \cite[Theorem~A(i)]{BeDiPo}, the branch $\mathbf f_0$ is \'etale over weight space, so after enlarging $L$ if necessary we may use the same completed local weight ring $A=L[[X]]$ for both branches, with $X=0$ corresponding to weight one. Here $X$ varies the weight; the cyclotomic variable remains independent. All the period modules are considered after inverting $p$, as in the characteristic-zero localization of \cite[Remark~1.3]{SW22}.

Finally, we use the cyclotomic conventions of Section~\ref{subsec:cyclotomic-conventions}. All the measures in this section are considered separately on each even $\omega^a$-component, and the analytic variable $j$ stands for evaluation at $\omega^a\langle\cdot\rangle^j$; the formulas displayed in terms of $j$ correspond to $a=0$.

\subsection{Cohomological periods and tame twisting}
\label{subsec:regular-compact-twisting}

Since the reciprocity computation will take place on the auxiliary family $\mathbf f_0$, we need a map from the periods of $\mathbf f$ to those of $\mathbf f_0$ which is defined over $A=L[[X]]$. Its regularity follows from the geometry of the twisting correspondence (Proposition~\ref{prop:compact-twisting-regular}), while its nonvanishing at weight one will be a consequence of the reciprocity computation of \S\ref{subsec:twisted-reciprocity-calculation}. We begin by introducing the relevant cohomology modules and the modular-symbol measure.

\subsubsection*{Compact period modules}
For $\mathbf g\in\{\mathbf f_0,\mathbf f\}$, put $M_{\mathbf f_0}=M$ and $M_{\mathbf f}=M^2$, and consider the ordinary inverse limit
\[
 H(M_{\mathbf g})
 =\varprojlim_s
 H^1_{\et}\bigl(X_1(M_{\mathbf g}p^s)_{\overline\Q},\cO_L\bigr)^{\ord},
\]
with trace transition maps and dual Hecke operators. We take the finite diamond-character component of the branch and the signed module $H(M_{\mathbf g})^-(1)$, as in \cite{FK,SW21}. Here the sign refers to the $-1$-eigenspace of complex conjugation before the Tate twist, which becomes the $+1$-eigenspace after twisting; Poincar\'e duality identifies it with the plus part of modular-symbol homology. Since complex conjugation sends $[u:v]$ to $[-u:v]$, this is the sign which carries the even cyclotomic moments. We stress that, throughout, compact cohomology refers to the $H^1$ of the compactified curve $X_1$, and not to the $H_c^1$ of the open curve $Y_1$.

Let $\Lambda_{\mathrm{wt}}$ stand for the ordinary weight algebra on this finite component, and let $\mathscr H_{\mathbf g}$ be the previous signed module after extension of scalars from $\Lambda_{\mathrm{wt}}$ to $A$. Let $\mathbb T_{M_{\mathbf g}}$ be the cuspidal Hecke algebra generated by the operators $T_\ell$ for $\ell\nmid Mp$, the tame diamond operators and $U_p$, and let $\lambda_{\mathbf g}:\mathbb T_{M_{\mathbf g}}\to A$ be the eigensystem of $\mathbf g$. We do not include the operators $U_\ell$ for $\ell\mid M$, since tame twisting changes their local newvector eigenvalues. Define
\begin{equation}\label{eq:compact-period-quotient}
 \mathscr P_{\mathbf g,\mathrm{cpt}}
 =\left(
 \mathscr H_{\mathbf g}\Big/
 \bigl(T-\lambda_{\mathbf g}(T):T\in\mathbb T_{M_{\mathbf g}}\bigr)
             \mathscr H_{\mathbf g}
 \right)_{\mathrm{tf}},
 \qquad
 q_{\mathbf g}:\mathscr H_{\mathbf g}\twoheadrightarrow\mathscr P_{\mathbf g,\mathrm{cpt}},
\end{equation}
where $(-)_{\mathrm{tf}}$ denotes the quotient by the $A$-torsion. In other words, we impose the Hecke eigenvalues of the family and then discard the elements killed by a power of $X$. This finite $A$-module is free of rank one; this follows from ordinary control and multiplicity one at the minimal levels, as explained in the proof of Proposition~\ref{prop:compact-twisting-regular} below. We refer to $\mathscr P_{\mathbf g,\mathrm{cpt}}$ as the \emph{compact period module}. Its special fibre $\overline{\mathscr P}_{\mathbf g,\mathrm{cpt}}:=\mathscr P_{\mathbf g,\mathrm{cpt}}/X\mathscr P_{\mathbf g,\mathrm{cpt}}$ is a one-dimensional $L$-vector space, which records one complex-conjugation sign of the cohomology.

\subsubsection*{The modular-symbol measure}
We now define the ordinary modular-symbol measure $\boldsymbol{\mathscr L}_{\mathbf g,\mathrm{cpt}}$, with values in $\mathscr P_{\mathbf g,\mathrm{cpt}}$. At level $M_{\mathbf g}p^r$, let $[u:v]_r^+$ stand for the plus Manin symbol and let $e^{\ord}$ be the ordinary projector. For $a\in(\Z/p^r\Z)^\times$, the value on $a+p^r\Z_p$ is the image of $U_p^{-r}e^{\ord}[M_{\mathbf g}a:1]_r^+$ in the family quotient. Equivalently, the measure is the system
\begin{equation}\label{eq:compact-modular-measure-definition}
 \boldsymbol{\mathscr L}_{\mathbf g,\mathrm{cpt}}
 =\left(
   \sum_{a\in(\Z/p^r\Z)^\times}
      q_{\mathbf g}\bigl(U_p^{-r}e^{\ord}
          [M_{\mathbf g}a:1]_r^+\bigr)\otimes[a]
       \right)_{r\geq1}.
\end{equation}
Following \cite[\S2.2]{SW22}, the symbols are identified along the ordinary tower before passing to the family quotient, and $[M_{\mathbf g}a:1]_r$ represents the path $\{\infty,a/p^r\}$; the trace relation, normalized by $U_p^{-r}$, gives
\[
 \boldsymbol{\mathscr L}_{\mathbf g,\mathrm{cpt}}(a+p^r\Z_p)
 =\sum_{b=0}^{p-1}\boldsymbol{\mathscr L}_{\mathbf g,\mathrm{cpt}}
       (a+bp^r+p^{r+1}\Z_p).
\]
Thus the values are compatible with subdivision into residue classes. The fact that these symbols lie in compact homology is \cite[Lemma~3.1]{sharifi12}: the ordinary projection of $[u:v]_r$ belongs to compact homology as soon as both $u$ and $v$ are nonzero modulo $p^r$, as follows from the ordinary cusp computation of \cite[Proposition~4.3.4]{ohta99}. Both entries of the symbols in \eqref{eq:compact-modular-measure-definition} satisfy this condition, and since finite tame diamond projectors preserve compact homology, the same holds on our branch. In contrast, the symbol $\{0,\infty\}$ has a zero entry; it gives rise to the boundary class of Section~\ref{subsec:period-boundary}.

In order to obtain a scalar-valued measure, we choose an $A$-basis
\[
 \gamma_{\mathbf g}\in\mathscr P_{\mathbf g,\mathrm{cpt}},
 \qquad \mathscr P_{\mathbf g,\mathrm{cpt}}=A\gamma_{\mathbf g},
\]
and we evaluate against the dual basis $\gamma_{\mathbf g}^{\vee}\in\operatorname{Hom}_A(\mathscr P_{\mathbf g,\mathrm{cpt}},A)$:
\[
 \mathscr L_p^{\mathrm{cpt}}(\mathbf g)
   =(\gamma_{\mathbf g}^{\vee}\otimes1)
          (\boldsymbol{\mathscr L}_{\mathbf g,\mathrm{cpt}}),
 \qquad\text{so that}\qquad
 \boldsymbol{\mathscr L}_{\mathbf g,\mathrm{cpt}}
   =\gamma_{\mathbf g}\otimes\mathscr L_p^{\mathrm{cpt}}(\mathbf g).
\]
If $\gamma_{\mathbf g}$ is replaced by $v(X)\gamma_{\mathbf g}$, with $v\in A^\times$, the scalar measure gets multiplied by $v(X)^{-1}$.

We write $\mathscr L_p^{\mathrm{cpt}}(f)$ for the specialization of $\mathscr L_p^{\mathrm{cpt}}(\mathbf f)$ at $X=0$. Before choosing a basis, the specialized measure belongs to $\overline{\mathscr P}_{\mathbf f,\mathrm{cpt}}\otimes_L\Lambda_L$, and to $\Lambda_L$ afterwards. Let us check that it is bounded. The integral ordinary cohomology module is finitely generated over the integral weight algebra, and its modular measure is integral. After specialization at weight one, its image is contained in a finitely generated $\cO_L$-module of the finite-dimensional specialized cohomology. Projecting to $\overline{\mathscr P}_{\mathbf f,\mathrm{cpt}}$ and applying a fixed $L$-linear coordinate introduces a fixed power of $p$ in the denominator, independent of the conductor of the cyclotomic character, and hence the scalar measure is bounded.

\subsubsection*{Interpolation and the generic comparison of periods}
We follow the direct Mellin convention for modular symbols, as in \cite[\S2.1--2.2]{SW22}. At a classical point $y$ of weight $k\geq2$, a primitive character $\rho$ of conductor $p^r$ and an integer $0\leq j\leq k-2$ with $\rho(-1)(-1)^j=1$, the scalar moment at $z\mapsto\rho(z)z^j$ is
\[
 \frac{j!\,p^{rj}\tau(\rho)}{(-2\pi i)^j a_p(g_y)^r}
       \frac{L(g_y,\rho^{-1},1+j)}
            {\Omega_{g_y,\gamma_{\mathbf g}(y)}},
\]
where the period convention includes the corresponding factor at $j=0$. This agrees with the Mazur--Kitagawa interpolation of Section~\ref{subsec:factorization-kato} once the same conventions for finite characters and tame Euler factors are fixed. On the generic cuspidal locus, the two cohomological realizations identify the same rational signed eigenspace, so their period bases differ by an element $q(X)\in L((X))^\times$, independent of the cyclotomic moment. Comparing the two interpolation formulas we obtain
\begin{equation}\label{eq:generic-compact-period-comparison}
 \mathscr L_p(\mathbf f)=q(X)\mathscr L_p^{\mathrm{cpt}}(\mathbf f).
\end{equation}
Indeed, this identity holds at the classical interpolation pairs, which are Zariski-dense in the product of weight space and cyclotomic space, and hence it holds on the family. This is the generic period comparison used in Hypothesis~\ref{hyp:kato-specialization}(iii). Observe, however, that at weight one there are no complex critical integers: the specialization $\mathscr L_p^{\mathrm{cpt}}(f)$ is defined through the cohomological construction above, and not through an interpolation property.

\subsubsection*{The twisting map}
Write $\tau_m(\beta)=\sum_{a\bmod m}\beta(a)\zeta_m^a$ for the Gauss sum of a character $\beta$ modulo $m$, where the primitive roots of unity are chosen compatibly. We can now construct the twisting map between the two compact period modules.

\begin{proposition}\label{prop:compact-twisting-regular}
Under \eqref{eq:factorization-conductor-hypotheses}, tame twisting induces an $A$-linear map
\[
 t_\eta:\mathscr P_{\mathbf f,\mathrm{cpt}}
       \longrightarrow\mathscr P_{\mathbf f_0,\mathrm{cpt}},
 \qquad
 t_\eta(\gamma_{\mathbf f})=b_\eta(X)\gamma_{\mathbf f_0},
 \quad b_\eta(X)\in A.
\]
Both modules are free of rank one, and the map is compatible with reduction modulo $X^n$ for every $n\geq1$.
\end{proposition}
\begin{proof}
We first check that both modules are free of rank one. At a prime $\ell\mid M$, the auxiliary branch has fixed inertial characters $1$ and $\chi_\ell^{-1}$ and zero monodromy, by \cite[Proposition~4.4]{BeDiPo}. Hence its generic local representation is a principal series with one unramified character and one character of conductor exponent $n=v_\ell(M)$. After twisting by $\eta$, the inertial characters become $\chi_{1,\ell}$ and $\chi_{2,\ell}$, both of conductor exponent $n$, so that the conductor exponent is $2n$. Therefore the generic tame levels are exactly $M$ and $M^2$, and multiplicity one at these minimal levels gives rank one for the signed generic Hecke quotients. Their torsion-free quotients over the discrete valuation ring $A$ are consequently free of rank one. Note also that imposing the omitted $U_\ell$-eigenvalues gives the same torsion-free quotient, since in both cases its kernel is the saturation of the kernel of the same generic quotient.

We next describe the correspondence explicitly. For $b\in(\Z/M\Z)^\times$, put $u_b=\left(\begin{smallmatrix}1&b/M\\0&1\end{smallmatrix}\right)$. Then
\[
 u_b\Gamma_1(M^2p^s)u_b^{-1}\subseteq\Gamma_1(Mp^s).
\]
Indeed, conjugating $\left(\begin{smallmatrix}a&d\\c&e\end{smallmatrix}\right)$ gives diagonal entries $a+bc/M$ and $e-bc/M$ and upper-right entry $d+b(e-a)/M-b^2c/M^2$, and the defining congruences at level $M^2p^s$ give the assertion. Hence
\[
 C_{b,s}:X_1(M^2p^s)\longrightarrow X_1(Mp^s),
 \qquad z\longmapsto z+b/M
\]
is a finite map, which extends across the cusps. Under Poincar\'e duality, its cohomological trace corresponds to pushforward on modular-symbol homology. The normalized sum
\[
 \widetilde t_\eta
   =\tau_M(\eta^{-1})^{-1}
       \sum_{b\in(\Z/M\Z)^\times}\eta^{-1}(b)(C_{b,s})_*
\]
is the transpose of the usual twisting operator on modular forms; in particular, on paths it is given by the corresponding sum of translated paths.

The maps $C_{b,s}$ commute with the forgetful maps in the $p$-power tower, so functoriality of the trace gives compatible maps on the inverse limit. The twisting identities read
\[
 \widetilde t_\eta T_\ell
       =\eta(\ell)T_\ell\widetilde t_\eta\quad(\ell\nmid Mp),
 \qquad
 \widetilde t_\eta U_p
       =\eta(p)U_p\widetilde t_\eta,
\]
where the second one follows by permuting the tame translation indices under multiplication by $p$. Since $\eta(p)$ is a unit, $\widetilde t_\eta$ preserves the ordinary summand, and since $\eta$ is even, it preserves the sign. All the coefficients of the correspondence, including the inverse Gauss sum, lie in $L$, so it is defined over $A$.

Finally, we descend the correspondence to the period modules. By the Hecke identities and multiplicity one, the composition $q_{\mathbf f_0}\widetilde t_\eta$ factors through $q_{\mathbf f}$ on the generic fibre. Hence, if an element lies in $\ker q_{\mathbf f}$, its image vanishes after inverting $X$; since the target is torsion-free, this image is already zero over $A$. This gives the map $t_\eta$ and proves its compatibility with reduction.
\end{proof}

Therefore, the comparison of compact periods under twisting is regular at weight one. It remains to prove that its special fibre is nonzero, and this will follow from the computation of the next subsection.

\subsection{The modular-symbol and reciprocity calculation}
\label{subsec:twisted-reciprocity-calculation}

In this subsection we evaluate the twisted modular symbols of the auxiliary family by means of the Sharifi map. More precisely, we compute the composition
\[
 \mathscr P_{\mathbf f,\mathrm{cpt}}
   \xrightarrow{\ t_\eta\ }\mathscr P_{\mathbf f_0,\mathrm{cpt}}
   \longrightarrow\overline{\mathscr P}_{\mathbf f_0,\mathrm{cpt}}
   \xrightarrow[\sim]{\ \varpi\ }\mathscr B_\chi
   \xrightarrow{\ \ell_\chi\ }L
\]
on the modular measure, where the target $\mathscr B_\chi$ and its coordinate $\ell_\chi$ are introduced below. The first arrow is the map of Proposition~\ref{prop:compact-twisting-regular}, the second one is specialization at weight one, and $\varpi$ is the localized Sharifi isomorphism. The twisting identity below expresses the image measure as a tame-twisted modular measure on the auxiliary family. After specialization, Sharifi sends its symbols to cup products of cyclotomic units, and $\ell_\chi$ evaluates them by local Hilbert symbols. In Lemma~\ref{lem:twisted-local-reciprocity}, explicit reciprocity separates this evaluation into the logarithm of a twisted cyclotomic unit and a sum of logarithmic derivatives of Coleman series. These give the two Dirichlet factors.

\subsubsection*{Twisted modular measures}
Let $\rho$ be an even primitive character of conductor $p^r$, with $r\geq2$. On the auxiliary branch, define the tame-twisted moment
\[
 \boldsymbol{\mathscr M}_{\mathbf f_0,\eta}(\rho)
 =q_{\mathbf f_0}\left(
    \sum_{b\in(\Z/Mp^r\Z)^\times}
       \eta^{-1}(b)\rho(b)U_p^{-r}e^{\ord}[b:1]_r^+
       \right).
\]
As in \eqref{eq:compact-modular-measure-definition}, each ordinary symbol is compact because both entries are nonzero modulo $p^r$. The factor $\eta^{-1}$ selects the tame twist, while $U_p^{-r}$ gives compatibility along the ordinary tower. The values lie in $\mathscr P_{\mathbf f_0,\mathrm{cpt}}\otimes_L L(\rho)$, and we write $\boldsymbol{\mathscr M}_{\mathbf f_0,\eta}|_{X=0}(\rho)$ for the result of specializing this module-valued measure at $X=0$ and then evaluating at $\rho$. Here $[M]$ is the Dirac measure on $\Gamma$, so that $[M](\rho)=\rho(M)$; this is the source of the conductor factor in Theorem~A.

The twisting correspondence gives the identity of module-valued measures
\begin{equation}\label{eq:twisted-modular-measures}
 \boldsymbol{\mathscr M}_{\mathbf f_0,\eta}
    =\tau_M(\eta^{-1})[M]\,
          t_\eta(\boldsymbol{\mathscr L}_{\mathbf f,\mathrm{cpt}}).
\end{equation}
In order to prove it, and in particular to check the direction of $[M]$, apply the normalized correspondence to the moment in \eqref{eq:compact-modular-measure-definition}. In modular-symbol notation, we get
\[
 t_\eta\bigl(\boldsymbol{\mathscr L}_{\mathbf f,\mathrm{cpt}}(\rho)\bigr)
 =\frac{a_p(\mathbf f)^{-r}}{\tau_M(\eta^{-1})}
   \sum_{\substack{a\bmod p^r\,,\ b\bmod M\\(a,p)=(b,M)=1}}
        \rho(a)\eta^{-1}(b)
        \left\{\infty,\frac a{p^r}+\frac bM\right\},
\]
where the ordinary projection and the auxiliary quotient are understood on the right-hand side. Set $c=Ma+p^rb$. The Chinese remainder theorem identifies these pairs $(a,b)$ with $(\Z/Mp^r\Z)^\times$, and
\[
 \rho(a)\eta^{-1}(b)
   =\rho(M)^{-1}\eta(p)^r\rho(c)\eta^{-1}(c),
 \qquad a_p(\mathbf f)=\eta(p)a_p(\mathbf f_0).
\]
Since $a/p^r+b/M=c/(Mp^r)$, it follows that
\[
 t_\eta\bigl(\boldsymbol{\mathscr L}_{\mathbf f,\mathrm{cpt}}(\rho)\bigr)
   =\tau_M(\eta^{-1})^{-1}\rho(M)^{-1}
      \boldsymbol{\mathscr M}_{\mathbf f_0,\eta}(\rho),
\]
which is \eqref{eq:twisted-modular-measures}. In particular, solving for the modular measure of $\mathbf f$ produces $[M]^{-1}$, which is the origin of the factor $\langle M\rangle^{-j}$ in Theorem~A.

\subsubsection*{The target of the Sharifi map}
Put $\Theta=\chi\omega$, and let $\mathfrak p$ be the exceptional height-one weight prime of the ordinary weight algebra $\Lambda_\Theta$. Consider
\[
 S_\Theta=
   \left(\varprojlim_n
      H^2_{\et}(\Z[\zeta_{Mp^n},1/p],\Z_p(2))\right)_\Theta,
 \qquad
 \mathscr B_\chi=S_{\Theta,\mathfrak p}
          \otimes_{k_{\Theta,\mathfrak p}}L,
\]
where $k_{\Theta,\mathfrak p}$ is the residue field of $\Lambda_{\Theta,\mathfrak p}$. Recall that local invariants identify the local groups $H^2(\,\cdot\,,\Z_p(1))$ at the primes $v\mid p$ with $\Z_p$, and that the global sum of these invariants is zero. Hence, if $\Sigma$ denotes the set of primes above $p$ in $\Q(\mu_{Mp^\infty})$ and $\oplus^0$ stands for the kernel of the summation map, \cite[Lemma~3.1 and Proposition~3.2]{SW21} give
\[
 S_{\Theta,\mathfrak p}
  \simeq\left((\mathop{\bigoplus\!{}^0}_{v\in\Sigma}\Z_p)_\chi(1)\right)_{\mathfrak p}
  \simeq k_{\Theta,\mathfrak p}.
\]
Here the class-group term vanishes by \cite[Lemma~3.1]{SW21}, since its Tate twist has characteristic ideal prime to the exceptional weight prime. Since $\chi$ is nontrivial and $\chi(p)=1$, the $\chi$-component of the permutation module on the primes above $p$ has dimension one and lies in the kernel of summation. This explains both the dimension of $\mathscr B_\chi$ and its interpretation in terms of local invariants.

In order to fix a coordinate, choose the prime $v_0\mid p$ determined by $\iota_p$, and take the compatible roots $\zeta_{p^n}$ as a basis of the Tate twist. Put $G=(\Z/M\Z)^\times$, let $D_p=\langle p\rangle\subset G$, and write $f_M=|D_p|$. If $\delta_{v_0}$ denotes the vector supported at $v_0$, define
\[
 e_\chi=\frac{\varphi(M)}{f_M}
       \left(\frac1{\varphi(M)}\sum_{g\in G}\chi^{-1}(g)g\right)
           \delta_{v_0}\otimes(\zeta_{p^n})_n.
\]
Since $\chi|_{D_p}=1$ and $\chi\ne1$, this is a nonzero vector in the kernel of summation, whose $v_0$-coordinate is $1$ in the chosen Tate basis. We let $\ell_\chi:\mathscr B_\chi\to L$ be this coordinate, so that $\ell_\chi(e_\chi)=1$. Concretely, $\ell_\chi$ is the local invariant at $v_0$, and we use the usual convention in which the exponent of the Hilbert symbol is identified with the corresponding element of $\Z/p^n\Z$. Thus the coordinate of a cup product is computed by its local Hilbert symbol at $v_0$.

\subsubsection*{The localized Sharifi map}
The source of \cite[Theorem~1.2]{SW21} is
\[
 H^-_{\Theta,\mathfrak p}(1)\big/
       I^*_{\Theta,\mathfrak p}H^-_{\Theta,\mathfrak p}(1),
 \qquad
 H=\varprojlim_n H^1_{\et}
       (X_1(Mp^n)_{\overline\Q},\Z_p)^{\ord},
\]
and in particular $H$ is the cohomology of the compactified curve. Their theorem identifies this quotient with $S_{\Theta,\mathfrak p}$, by sending a compact ordinary Manin symbol $[u:v]$ to the $\Theta$-projection of the cup product $\{1-\zeta^u,1-\zeta^v\}$ of the Kummer classes of the corresponding cyclotomic units. Moreover, the localized cuspidal Hecke algebra is the weight ring by \cite[Theorem~A(i)]{BeDiPo}, and \cite[Proposition~3.5 and the paragraph preceding it]{SW21} shows that the signed compact module is free of rank one over this algebra and identifies its Eisenstein ideal with the exceptional weight prime. After completion, this ideal becomes $(X)$, and the source of the theorem is exactly $\overline{\mathscr P}_{\mathbf f_0,\mathrm{cpt}}$. We thus obtain an isomorphism
\begin{equation}\label{eq:localized-Sharifi-period}
 \varpi:\overline{\mathscr P}_{\mathbf f_0,\mathrm{cpt}}
                    \xrightarrow{\ \sim\ }\mathscr B_\chi.
\end{equation}
We emphasize that this identification comes from \cite[Theorem~1.2]{SW21}, which is formulated after localization at $\mathfrak p$, and that it is defined on compact cohomology itself; the Drinfeld--Manin enlargement of compact cohomology, which we consider in Proposition~\ref{prop:auxiliary-boundary-length}, plays no role here.

\subsubsection*{The reciprocity computation}
We now evaluate the image under the localized Sharifi map. The resulting formula is a tame-twisted version of \cite[Proposition~2.6]{SW22}, which evaluates, in its notation,
\[
 \sum_{a\in(\Z/p^r\Z)^\times}
  \psi^{-1}(a)(1-\zeta_{p^r}^{a},
                   1-\zeta_M^{p^{-r}}\zeta_{p^r})_{r,\Theta}
   =-\frac{\tau_M(\chi^{-1})}{\varphi(M)}
        L_p^{\mathrm{KL}}(\psi,1)
        L_p^{\mathrm{KL}}(\chi\omega\psi^{-1},0)
\]
modulo $p^r$, for an even nontrivial character $\psi$ of conductor $p^r$ with $r\geq2$. In our setting, the first entry of the cup product also carries the tame character $\eta^{-1}$, and the main point of the computation below is to keep track of this character through the trace. As in loc.\,cit., we restrict to characters with $r\geq2$; they are enough to determine the bounded measures.

Two kinds of sums will appear. For a modulus $m$ and characters $\nu$, $\vartheta$ modulo $m$, put
\begin{equation}\label{eq:twisted-local-raw}
 S_m(\nu)=\sum_{a\in(\Z/m\Z)^\times}
       \nu^{-1}(a)\log_p(1-\zeta_m^a),
 \qquad
 B_m(\vartheta)=\sum_{b\in(\Z/m\Z)^\times}
       \vartheta^{-1}(b)\frac{\zeta_m^b}{1-\zeta_m^b}.
\end{equation}
The first one is the logarithmic sum of Leopoldt's formula for $L_p^{\mathrm{KL}}(\nu,1)$, and it arises from the logarithm of the first entry of the cup product. The second one is a finite Fourier expression for $L(\vartheta,0)$, and it arises from the logarithmic derivative of the second entry.

\begin{lemma}\label{lem:twisted-local-reciprocity}
With the coordinate $\ell_\chi$ fixed above, put $d_\chi=f_M/(M\varphi(M))$. Then, for every even primitive character $\rho$ of conductor $p^r$ with $r\geq2$,
\begin{equation}\label{eq:twisted-local-product}
 \ell_\chi\varpi\bigl(\boldsymbol{\mathscr M}_{\mathbf f_0,\eta}|_{X=0}(\rho)\bigr)
 =d_\chi M\frac{\tau_M(\chi_2)}{\tau_M(\chi_1)}
       L_p^{\mathrm{KL}}(\chi_1\rho^{-1},1)
       L_p^{\mathrm{KL}}(\chi_2^{-1}\omega\rho,0).
\end{equation}
\end{lemma}
\begin{proof}
Fix $\rho$ of conductor $p^r$, and put $m=Mp^r$, $\nu=\chi_1\rho^{-1}$ and $\vartheta=\chi_2^{-1}\rho$. Then $\nu$ is even, $\vartheta$ is odd, $\nu\vartheta=\chi$, and both characters are primitive of conductor $m$. Since $\eta=\chi_1$, we have $\eta^{-1}\rho=\nu^{-1}$. We divide the proof into five steps.

\emph{Step 1: from modular symbols to a single cup product.} Fix $n\geq r$. To compute Hilbert symbols modulo $p^n$, we work at level $Mp^n$, and we take as primitive $Mp^n$-th root of unity
\[
 q_n=\zeta_M^{p^{-n}}\zeta_{p^n}^{M^{-1}},
\]
where $p^{-n}$ is taken modulo $M$ and $M^{-1}$ modulo $p^n$. Then $q_n^p=q_{n-1}$, so $(q_n)_n$ is a compatible system of primitive $Mp^n$-th roots of unity, which we use in the description of $\varpi$; in particular, the root $\zeta_m$ is $q_r$. Extend $\nu$ to a character modulo $Mp^n$. By the compatibility of the measure along the tower, the $\rho$-moment can be computed at level $n$; note that at $X=0$ the operator $U_p$ acts on $\overline{\mathscr P}_{\mathbf f_0,\mathrm{cpt}}$ through $a_p(f_0)=1$, so that the normalization by $U_p^{-n}$ does not contribute. Then, by the description of $\varpi$ recalled before \eqref{eq:localized-Sharifi-period}, at this level the element $\varpi\bigl(\boldsymbol{\mathscr M}_{\mathbf f_0,\eta}|_{X=0}(\rho)\bigr)$ is represented by the $\Theta$-projection of
\[
 \sum_{b\in(\Z/Mp^n\Z)^\times}
       \eta^{-1}(b)\rho(b)
       \{1-q_n^b,1-q_n\}
 =\sum_{b\in(\Z/Mp^n\Z)^\times}
       \nu^{-1}(b)
       \{1-q_n^b,1-q_n\}
 =\{u_{\nu,n},1-q_n\}.
\]
Here the braces denote the cup product of the two Kummer classes, with values in $H^2$, and the last equality is the bilinearity of the cup product, with
\[
 u_{\nu,n}=\sum_{b\in(\Z/Mp^n\Z)^\times}\nu^{-1}(b)\,(1-q_n^b),
\]
written additively in the units of $\Q(\mu_{Mp^n})$ tensored with $L(\rho)$. Observe that each $1-q_n^b$ is a unit, since $Mp^n$ is not a prime power (recall that $M>1$); this is needed to take $p$-adic logarithms and to apply the reciprocity law below. This is the tame-character version of the compact-symbol construction of \cite[\S2.4]{SW22}. The moment is therefore represented by the projected cup product $e_\Theta\{u_{\nu,n},1-q_n\}$. Reindexing the first entry gives
\[
 \sigma_c(u_{\nu,n})
 =\sum_b\nu^{-1}(c^{-1}b)(1-q_n^b)
 =\nu(c)u_{\nu,n},
 \qquad
 \log_p(u_{\nu,n})=S_{Mp^n}(\nu),
\]
where $\sigma_c(q_n)=q_n^c$ and the logarithm is computed through $\iota_p$. The Galois action is on the unit factor, with the character coefficients fixed. We will use this eigenvector identity to separate the two entries in the reciprocity formula.

\emph{Step 2: from the cup product to Hilbert symbols.} For $g\in G=(\Z/M\Z)^\times$ and $\delta\in\Delta$, let $\sigma_{g,\delta}$ be the automorphism raising the tame root of unity to the power $g$ and the $p$-power root of unity to the Teichm\"uller lift $\omega(\delta)$, and put $q_n^{(g,\delta)}=\sigma_{g,\delta}(q_n)$. The $\Theta$-projector is
\[
 e_\Theta=\frac{1}{\varphi(M)(p-1)}\sum_{g\in G,\,\delta\in\Delta}
        \chi^{-1}(g)\omega^{-1}(\delta)\,\sigma_{g,\delta}.
\]
Recall that $\ell_\chi$ is the invariant at $v_0$, and that the local invariant of the cup product of the Kummer classes of two units is their Hilbert symbol. Since the cup product is Galois equivariant and $\sigma_{g,\delta}(u_{\nu,n})=\nu(\sigma_{g,\delta})u_{\nu,n}$ by Step~1, we find, modulo $p^n$,
\[
 \ell_\chi\bigl(e_\Theta\{u_{\nu,n},1-q_n\}\bigr)
 =\frac{1}{\varphi(M)(p-1)}\sum_{g,\delta}
        \chi^{-1}(g)\omega^{-1}(\delta)\nu(\sigma_{g,\delta})\,
        \bigl(u_{\nu,n},\,1-q_n^{(g,\delta)}\bigr)_{v_0}.
\]
This is the precise sense in which the left-hand side of \eqref{eq:twisted-local-product} is a twisted sum of Hilbert symbols, always with the same first entry $u_{\nu,n}$.

\emph{Step 3: explicit reciprocity and separation of the two sums.} We now apply the explicit reciprocity law used in \cite[Proposition~2.6]{SW22}. Put $K_n=\Q_p(\mu_{Mp^n})$ and $Z_n=\zeta_{p^n}-1$. If $D_{g,\delta}$ is the logarithmic derivative of the Coleman series for $1-q_n^{(g,\delta)}$, evaluated at $Z_n$ with the Frobenius convention below, the law computes its Hilbert symbol against a unit $x$ by
\[
 (x,1-q_n^{(g,\delta)})_{v_0}
 =p^{-n}\operatorname{Tr}_{K_n/\Q_p}
       \bigl(\log_p(x)D_{g,\delta}\bigr)\pmod{p^n}.
\]
The law applies here since both entries are units of $K_n$. We use the exponent convention for Hilbert symbols fixed above and extend linearly in the character coefficients. Note that the constant $d_\chi$ obtained below depends on the normalization of the coordinate $\ell_\chi$ and on this sign convention, so a direct comparison with the constant of \cite[Proposition~2.6]{SW22} requires matching these normalizations. Thus we must compute the logarithm of the first entry, the Coleman derivative of the second, and their trace.

For $x=u_{\nu,n}$, Step~1 gives $\tau(\log_p u_{\nu,n})=\nu(\tau)\log_p(u_{\nu,n})$. Expanding the trace therefore yields
\[
 \begin{aligned}
 \operatorname{Tr}_{K_n/\Q_p}\bigl(\log_p(u_{\nu,n})D_{g,\delta}\bigr)
 &=\sum_\tau\tau(\log_p u_{\nu,n})\tau(D_{g,\delta})\\
 &=S_{Mp^n}(\nu)\sum_\tau\nu(\tau)\tau(D_{g,\delta}).
 \end{aligned}
\]
This extracts the logarithmic sum as the first factor. The remaining sum carries the character weights left by the first entry and by $e_\Theta$.

To compute its terms, consider the Coleman series
\[
 \mathcal C(Z)=1-\zeta_M(1+Z)^{1/M}
       \in\cO_{\Q_p(\mu_M)}[[Z]],
\]
where the binomial expansion is defined since $M\in\Z_p^\times$. Let $\varphi$ be unramified Frobenius on the coefficients, so that $\varphi(\zeta_M)=\zeta_M^p$, and set
\[
 \mathcal C_n(Z)=\varphi^{-n}\mathcal C(Z)
       =1-\zeta_M^{p^{-n}}(1+Z)^{1/M}.
\]
Then $\mathcal C_n(Z_n)=1-q_n$. With $\partial=(1+Z)d/dZ$, differentiation followed by evaluation gives
\[
 \partial\log\mathcal C_n(Z)
 =-\frac1M\frac{\zeta_M^{p^{-n}}(1+Z)^{1/M}}
                  {1-\zeta_M^{p^{-n}}(1+Z)^{1/M}},
 \qquad
 (\partial\log\mathcal C_n)(Z_n)
       =-\frac1M\frac{q_n}{1-q_n}.
\]
For the conjugate $1-q_n^{(g,\delta)}$, replace the coefficient by $\zeta_M^{gp^{-n}}$ and the exponent by $\omega(\delta)/M$. The same calculation gives
\[
 D_{g,\delta}=-\frac{\omega(\delta)}M
       \frac{q_n^{(g,\delta)}}{1-q_n^{(g,\delta)}}.
\]
Each evaluation has produced one rational term. To obtain the Fourier sum, we combine its conjugates in the trace with the weights of the projector. The local trace runs over the unramified Frobenius powers $p^i$, $0\leq i<f_M$, and the automorphisms $t\in(\Z/p^n\Z)^\times$ of the cyclotomic part.

The Teichm\"uller factor $\omega(\delta)$ in $D_{g,\delta}$ cancels with the coefficient $\omega^{-1}(\delta)$ of the projector. Put $b=gp^i\in G$ and $c=\delta t\in(\Z/p^n\Z)^\times$. Since $\chi(p)=1$, we have $\chi^{-1}(g)=\chi^{-1}(b)$, and the product of the two character values $\nu(\sigma_{g,\delta})\nu(\tau)$ equals $\nu(b,c)=\eta(b)\rho^{-1}(c)$. Combined with $\chi^{-1}(b)$, this leaves
\[
 \chi^{-1}(b)\eta(b)\rho^{-1}(c)
       =\chi_2(b)\rho^{-1}(c)=\vartheta^{-1}(b,c),
\]
and the conjugate of $q_n$ indexed by $(b,c)$ is $q_n^{(b,c)}$. For each fixed pair $(b,c)$, the equations $b=gp^i$ and $c=\delta t$ have exactly $f_M(p-1)$ solutions with $0\leq i<f_M$ and $\delta\in\Delta$, since $g$ and $t$ are determined once $i$ and $\delta$ are fixed. Together with the normalization $1/(\varphi(M)(p-1))$ of the projector and the factor $-1/M$ of the logarithmic derivative, the second factor becomes
\[
 -\frac{f_M}{M\varphi(M)}\sum_{(b,c)}\vartheta^{-1}(b,c)\frac{q_n^{(b,c)}}{1-q_n^{(b,c)}}
 =-d_\chi B_{Mp^n}(\vartheta),
\]
since $q_n$ is a primitive $Mp^n$-th root of unity. In particular, the scalar $d_\chi=f_M/(M\varphi(M))$ is independent of $n$, $\eta$ and $\rho$. Including the factor $p^{-n}$ of the reciprocity law, we conclude that, modulo $p^n$, the left-hand side of \eqref{eq:twisted-local-product} is given by
\[
 -d_\chi p^{-n}S_{Mp^n}(\nu)B_{Mp^n}(\vartheta).
\]

\emph{Step 4: passage to characteristic zero.} The previous expression was obtained from Hilbert symbols modulo $p^n$, while the left-hand side of \eqref{eq:twisted-local-product} is a characteristic-zero quantity. We pass to the limit with $\rho$ fixed and $n\geq r$ tending to infinity. Regarding $\nu$ and $\vartheta$ as characters modulo $Mp^n$, as in Step~1, norm compatibility and the corresponding identity for logarithmic derivatives give
\[
 S_{Mp^n}(\nu)=S_m(\nu),\qquad
 B_{Mp^n}(\vartheta)=p^{n-r}B_m(\vartheta).
\]
Hence $p^{-n}S_{Mp^n}(\nu)B_{Mp^n}(\vartheta)=p^{-r}S_m(\nu)B_m(\vartheta)$ is independent of $n$. Moreover, for this fixed $\rho$ there is an integer $c\geq0$ such that the character projectors and the coordinate $\ell_\chi$ send the integral lattice of local invariants into $p^{-c}\cO_{L(\rho)}$. The reciprocity congruences then show that the difference between the left-hand side of \eqref{eq:twisted-local-product} and $-d_\chi p^{-r}S_m(\nu)B_m(\vartheta)$ belongs to $p^{n-c}\cO_{L(\rho)}$ for every $n\geq r$. Since the intersection of these ideals is zero, we get the equality
\[
 \ell_\chi\varpi\bigl(\boldsymbol{\mathscr M}_{\mathbf f_0,\eta}|_{X=0}(\rho)\bigr)
   =-d_\chi p^{-r}S_m(\nu)B_m(\vartheta)
\]
in $L(\rho)$.

\emph{Step 5: the two Kubota--Leopoldt values.} It remains to identify the two sums with $p$-adic $L$-values. We have
\[
 S_m(\nu)=-\frac{m}{\tau_m(\nu)}L_p^{\mathrm{KL}}(\nu,1),
 \qquad
 B_m(\vartheta)=\tau_m(\vartheta^{-1})
                        L_p^{\mathrm{KL}}(\omega\vartheta,0).
\]
The first identity is Leopoldt's formula, in the convention of Section~\ref{subsec:cyclotomic-conventions}. For the second one, use that
\[
 \frac{z}{1-z}=\sum_{a=0}^{m-1}\left(\frac12-\frac a m\right)z^a-1
 \quad(z^m=1,\ z\ne1)
\]
together with $\sum_b\vartheta^{-1}(b)\zeta_m^{ab}=\vartheta(a)\tau_m(\vartheta^{-1})$; the constant term drops out because $\vartheta$ is nontrivial, and we recover the expression $L(\vartheta,0)=-\frac1m\sum_a\vartheta(a)a$, which also fixes the sign. In both identities the $p$-Euler factors are equal to $1$, since the characters are ramified at $p$. Substituting into the formula of Step~4, and using $m=Mp^r$, we obtain
\[
 -d_\chi p^{-r}S_m(\nu)B_m(\vartheta)
 =d_\chi M\frac{\tau_m(\vartheta^{-1})}{\tau_m(\nu)}
       L_p^{\mathrm{KL}}(\nu,1)L_p^{\mathrm{KL}}(\omega\vartheta,0).
\]
Finally, the multiplicativity of primitive Gauss sums gives
\[
 \frac{\tau_{Mp^r}(\vartheta^{-1})}{\tau_{Mp^r}(\nu)}
 =\left(\frac{\chi_2(p)}{\chi_1(p)}\right)^r
       \frac{\tau_M(\chi_2)}{\tau_M(\chi_1)}
 =\frac{\tau_M(\chi_2)}{\tau_M(\chi_1)},
\]
since the wild Gauss sums and the factors at $M$ cancel. As $\nu=\chi_1\rho^{-1}$ and $\omega\vartheta=\chi_2^{-1}\omega\rho$, this proves \eqref{eq:twisted-local-product}.
\end{proof}

\subsection{Factorization and the choice of periods}
\label{subsec:factorization-compact-periods}

The reciprocity calculation determines the image of the compact measure under twisting and shows that at least one of its finite-order moments is nonzero. Since this image factors through $t_\eta$, the twisting map is nonzero at weight one and therefore identifies the two compact period modules, allowing us to recover the measure of $\mathbf f$ itself.

\begin{theorem}\label{thm:compact-weight-one-factorization}
Under \eqref{eq:factorization-conductor-hypotheses}, the map $t_\eta$ is an isomorphism of compact period modules. Moreover, if the compact period at weight one is chosen to be
\[
 \bar\gamma_{\mathbf f}^{\mathrm{cpt}}
      =\bar t_\eta^{-1}\varpi^{-1}(e_\chi),
\]
then, for every even finite-order character $\rho$ of conductor $p^r$, with $r\geq2$,
\begin{equation}\label{eq:compact-finite-factorization}
 \mathscr L_p^{\mathrm{cpt}}(f)(\rho)
 =d_\chi\tau_M(\chi_2)\rho(M)^{-1}
       L_p^{\mathrm{KL}}(\chi_1\rho^{-1},1)
       L_p^{\mathrm{KL}}(\chi_2^{-1}\omega\rho,0).
\end{equation}
These values determine an equality of nonzero bounded measures on every even component.
\end{theorem}
\begin{proof}
We begin with arbitrary bases, and set $s_{\mathbf f_0}=\ell_\chi\varpi(\bar\gamma_{\mathbf f_0})$, which is nonzero by \eqref{eq:localized-Sharifi-period}. Applying $\ell_\chi\varpi$ to the specialization of \eqref{eq:twisted-modular-measures} and using Lemma~\ref{lem:twisted-local-reciprocity}, we get
\[
 \tau_M(\chi_1^{-1})\rho(M)s_{\mathbf f_0}b_\eta(0)
       \mathscr L_p^{\mathrm{cpt}}(f)(\rho)
 =d_\chi M\frac{\tau_M(\chi_2)}{\tau_M(\chi_1)}
       L_p^{\mathrm{KL}}(\chi_1\rho^{-1},1)
       L_p^{\mathrm{KL}}(\chi_2^{-1}\omega\rho,0).
\]
The two Dirichlet factors are nonzero bounded measures: their tame characters are nontrivial, so neither of them has the pole of the $p$-adic zeta function. On a fixed component, a nonzero element of $\cO_L[[\Gamma_1]][1/p]$ has only finitely many zeros at finite-order characters, by the Weierstrass preparation theorem. Consequently, the right-hand side is not identically zero, which forces $b_\eta(0)\ne0$, and Proposition~\ref{prop:compact-twisting-regular} then shows that $t_\eta$ is an isomorphism over $A$. Since $\chi_1$ is even, we have $\tau_M(\chi_1)\tau_M(\chi_1^{-1})=\chi_1(-1)M=M$. With the period specified in the statement we have $s_{\mathbf f_0}b_\eta(0)=1$, and the formula follows. The same uniqueness argument for bounded measures promotes the identity at finite-order characters to an identity of measures.
\end{proof}

The theorem normalizes the period using the Sharifi map. For other bases, the constant changes as follows.

\begin{remark}[Other compact periods]\label{rem:compact-period-bases}
For arbitrary bases, put
\[
 s_{\mathbf f_0}=\ell_\chi\varpi(\bar\gamma_{\mathbf f_0}),\qquad
 C_{f,\gamma}=
       \frac{d_\chi\tau_M(\chi_2)}{s_{\mathbf f_0}b_\eta(0)}.
\]
Then \eqref{eq:compact-finite-factorization} holds with $d_\chi\tau_M(\chi_2)$ replaced by $C_{f,\gamma}$. The two factors in the denominator compare the auxiliary basis with $e_\chi$, and the two compact bases under twisting, respectively. A change $\gamma_{\mathbf f}\mapsto v(X)\gamma_{\mathbf f}$ multiplies both the function and $C_{f,\gamma}$ by $v(0)^{-1}$, while the auxiliary basis cancels out. Finally, any two lifts to $A$ of the period specified at weight one differ by a unit in $1+XA$, and hence they give the same weight-one function.
\end{remark}

We finally express the factorization in terms of the analytic cyclotomic variable, which gives the formula of Theorem~A.

\begin{corollary}[The cyclotomic variable]\label{prop:cyclotomic-period-comparison}
On the trivial Teichm\"uller component, the factorization of Theorem~\ref{thm:compact-weight-one-factorization} reads
\begin{equation}\label{eq:weight-one-bd-factorization}
 \mathscr L_p^{\mathrm{cpt}}(f,j)
 =C_{f,\gamma}\langle M\rangle^{-j}
       L_p^{\mathrm{KL}}(\chi_1,1+j)
       L_p^{\mathrm{KL}}(\chi_2^{-1}\omega,-j).
\end{equation}
On the even $\omega^a$-component, the two characters have to be replaced by $\chi_1\omega^{-a}$ and $\chi_2^{-1}\omega^{1+a}$, and the conductor factor has to be multiplied by $\omega^{-a}(M)$.
\end{corollary}
\begin{proof}
Choose a topological generator $u$ of $\Gamma_1$, write $G_\beta(u^s-1)=L_p^{\mathrm{KL}}(\beta,1-s)$, and put $T=u^j-1$. In this variable, the identity at finite-order characters involves the two power series
\[
 G_{\chi_1}\bigl((1+T)^{-1}-1\bigr),\qquad
 G_{\chi_2^{-1}\omega}\bigl(u(1+T)-1\bigr),
\]
which correspond to the arguments $1+j$ and $-j$ in our convention, while the Dirac measure $[M]^{-1}$ gives $\langle M\rangle^{-j}$. Since the bounded measures are determined by their values at finite-order characters of conductor $p^r$ with $r\geq2$, this proves the formula. Evaluation at $\omega^a\langle\cdot\rangle^j$ gives the statement for the other components.
\end{proof}

\subsection{Changing the period module}
\label{subsec:period-boundary}

The compact specialization of Theorem~\ref{thm:compact-weight-one-factorization} is nonzero, but it may acquire a zero when it is transported to a different period module. For the cohomology of the open modular curve, the discrepancy comes from the cusps: restriction is an isomorphism on the generic cuspidal quotients, but its cokernel at the Eisenstein point may be nonzero. This is the family counterpart of the lattice comparison of \S\ref{subsec:boundary-quotient}, and the length of this cokernel measures the order of vanishing of the image of the compact measure.

As in \cite[\S2.2 and Appendix~A]{PR25} (see in particular Propositions~A.7 and~A.11 and Corollary~A.12 there), the comparison has to be made before specialization. Observe, however, that the module $V^c(\mathbf f)$ of loc.\,cit.\ is defined using compact support on the open curve, while our period module uses the $H^1$ of the compactification and a signed quotient. The boundary length must therefore be computed for the present realization, and we do so below for the auxiliary Drinfeld--Manin module.

\begin{proposition}\label{prop:factorization-open-periods}
For $\mathbf g\in\{\mathbf f,\mathbf f_0\}$, let $\mathscr P_{\mathbf g,\mathrm{op}}$ be the torsion-free signed family quotient of the cohomology of the open curve, formed with the Hecke algebra used in \eqref{eq:compact-period-quotient}. Both modules are free of rank one over $A$, and restriction induces an injection
\[
 j_{\mathbf g}:\mathscr P_{\mathbf g,\mathrm{cpt}}\longrightarrow\mathscr P_{\mathbf g,\mathrm{op}}
\]
with finite-length cokernel. Moreover, twisting descends to both realizations, and the square
\[
 \begin{tikzcd}[column sep=large]
 \mathscr P_{\mathbf f,\mathrm{cpt}}\arrow[r,"t_{\eta,\mathrm{cpt}}"]\arrow[d,"j_{\mathbf f}"']
   &\mathscr P_{\mathbf f_0,\mathrm{cpt}}\arrow[d,"j_{\mathbf f_0}"]\\
 \mathscr P_{\mathbf f,\mathrm{op}}\arrow[r,"t_{\eta,\mathrm{op}}"']
   &\mathscr P_{\mathbf f_0,\mathrm{op}}
 \end{tikzcd}
\]
commutes. Choose bases and write
\[
 j_{\mathbf g}(\gamma_{\mathbf g,\mathrm{cpt}})=c_{\mathbf g}(X)\gamma_{\mathbf g,\mathrm{op}},
 \qquad
 t_{\eta,\bullet}(\gamma_{\mathbf f,\bullet})
       =b_{\eta,\bullet}(X)\gamma_{\mathbf f_0,\bullet}.
\]
Then $c_{\mathbf g}\in A\setminus\{0\}$, and
\begin{equation}\label{eq:compact-open-period-relation}
 c_{\mathbf f_0}b_{\eta,\mathrm{cpt}}=b_{\eta,\mathrm{op}}c_{\mathbf f},
 \qquad
 \mathscr L^{\mathrm{im}}_{\mathbf g,\mathrm{op}}=c_{\mathbf g}(X)\mathscr L_p^{\mathrm{cpt}}(\mathbf g),
\end{equation}
where $\boldsymbol{\mathscr L}^{\mathrm{im}}_{\mathbf g,\mathrm{op}}:=(j_{\mathbf g}\otimes1)\boldsymbol{\mathscr L}_{\mathbf g,\mathrm{cpt}}$ is the image of the compact measure and the scalar functions are taken in the indicated bases. Finally,
\[
 \operatorname{length}_A\operatorname{coker}(j_{\mathbf g})
       =\operatorname{ord}_X(c_{\mathbf g}).
\]
\end{proposition}
\begin{proof}
The localization sequence for $Y_1\hookrightarrow X_1$ has Eisenstein boundary terms, so restriction is an isomorphism on the generic cuspidal quotients. Consequently, the signed open quotient has generic rank one, and it is free after removing the torsion. For restriction followed by the open quotient map, every relation in the compact family quotient vanishes generically, and hence it vanishes over $A$ because the target is torsion-free, as in the proof of Proposition~\ref{prop:compact-twisting-regular}. This gives $j_{\mathbf g}$ over $A$. Since a nonzero map between free rank-one modules over a discrete valuation ring is injective, with cokernel $A/(c_{\mathbf g})$ in our case, the length formula follows.

The same argument descends twisting to open cohomology. The geometric trace correspondences commute with restriction, so the square commutes both before and after passing to the family quotients. Applying the square to $\gamma_{\mathbf f,\mathrm{cpt}}$ gives the first identity in \eqref{eq:compact-open-period-relation}, and applying $j_{\mathbf g}$ to the module-valued measure and taking its coordinate in $\gamma_{\mathbf g,\mathrm{op}}$ gives the second one.
\end{proof}

The following remark describes the leading term of this image measure. Identifying it with an independently constructed measure on open cohomology requires a further comparison, which we formulate below for Kato's periods.

\begin{remark}[The leading term in an image measure]
Write $c_{\mathbf f}(X)=X^mv(X)$, with $v\in A^\times$, in Proposition~\ref{prop:factorization-open-periods}. Then the image measure satisfies
\begin{equation}\label{eq:factorization-Kato-period-input}
 \operatorname{ord}_X\mathscr L^{\mathrm{im}}_{\mathbf f,\mathrm{op}}=m,
 \qquad
 \left.X^{-m}\mathscr L^{\mathrm{im}}_{\mathbf f,\mathrm{op}}\right|_{X=0}
       =v(0)\mathscr L_p^{\mathrm{cpt}}(f).
\end{equation}
Indeed, the image is $X^mv(X)\mathscr L_p^{\mathrm{cpt}}(\mathbf f)$, and its compact factor has nonzero special fibre by Theorem~\ref{thm:compact-weight-one-factorization}. Observe that the order refers to the measure-valued function, and an individual cyclotomic moment may have an additional zero. Hence $m=0$ gives a formula for the value, while $m=1$ gives a formula for the first derivative in the weight direction. In order to apply this observation to Kato's $p$-adic $L$-function, one needs to identify its module-valued measure with such an image.
\end{remark}

The order $m$ does not change if we choose a different basis of either module, and in order to determine it one needs the boundary quotient of the particular realization. For the auxiliary family, Shih--Wang compute such a quotient by means of the Drinfeld--Manin construction, which we now recall. In \cite[\S3.2]{SW21}, $\widetilde H_\Theta$ is defined as the inverse image, inside open ordinary cohomology, of the single Eisenstein boundary line $\Lambda_\Theta e_{\Theta,1}(-1)$. Its Drinfeld--Manin quotient is
\[
 \widetilde H_{\Theta,\mathrm{DM}}
   =\widetilde H_\Theta\otimes_{\mathfrak H_\Theta^*}\mathfrak h_\Theta^*,
\]
where $\mathfrak H_\Theta^*$ is the full ordinary Hecke algebra and $\mathfrak h_\Theta^*$ is its cuspidal quotient. The Drinfeld--Manin splitting realizes this torsion-free module inside rational compact cohomology, where it contains the integral compact module. Its single boundary quotient is $\Lambda_\Theta/(\xi_\Theta)$, generated by the image of $\{0,\infty\}_{\Theta,\mathrm{DM}}$.

\begin{proposition}[The auxiliary Drinfeld--Manin period map]\label{prop:auxiliary-boundary-length}
For $f_0=E_1(1,\chi^{-1})^{\mathrm{ord}}$, consider the one-boundary Drinfeld--Manin module $\widetilde H_{\Theta,\mathrm{DM}}$ of \cite[\S3.2]{SW21}, where $\Theta=\chi\omega$. Take the negative sign and the Tate twist $(1)$, localize at the exceptional weight prime and complete to $A$. If $P_{f_0,\mathrm{DM}}$ denotes the resulting line, then the natural map fits in an exact sequence
\[
 0\longrightarrow\mathscr P_{\mathbf f_0,\mathrm{cpt}}\longrightarrow P_{f_0,\mathrm{DM}}
       \longrightarrow A/(X)\longrightarrow0.
\]
Consequently, the image of the nonzero compact modular measure in this Drinfeld--Manin line has a simple zero in the weight variable.
\end{proposition}
\begin{proof}
By \cite[Proposition~3.3]{SW21}, the quotient of the Drinfeld--Manin module by compact cohomology is $\Lambda_\Theta/(\xi_\Theta)$, generated by the image of $\{0,\infty\}_{\Theta,\mathrm{DM}}$, and this generator lies in the negative sign; see also \cite[Corollary~3.4]{SW21}. The same section identifies $(\xi_\Theta)$ with the exceptional weight prime after localization, so that completion gives $A/(X)$. The localized cuspidal Hecke algebra is $A$, and its compact signed module is free of rank one by \cite[Proposition~3.5]{SW21}. Since the Drinfeld--Manin module is torsion-free of the same rank, we obtain the displayed sequence. Therefore, a nonzero compact specialization acquires exactly one factor of $X$ in its image.
\end{proof}

This period map behaves differently from the Sharifi map \eqref{eq:localized-Sharifi-period}. The inclusion into the Drinfeld--Manin line becomes zero modulo $X$, whereas the Sharifi map is an isomorphism on the special fibre of compact cohomology; in particular, the latter does not factor through the former. This is consistent with \cite[proof of Theorem~1.2]{SW21}: the classes $(1-U_q^*)\{0,\infty\}_{\Theta,\mathrm{DM}}$ used there already belong to compact cohomology, while the unmodified boundary symbol does not.

We can now formulate the remaining comparison. Fix the modular-symbol realization in the transported Kato reciprocity law of Hypothesis~\ref{hyp:kato-specialization}, and let $\mathscr P_{\mathbf f}^{\mathrm K}$ denote its signed torsion-free family quotient, which is a rank-one $A$-module, and $\boldsymbol{\mathscr L}_{\mathbf f}^{\mathrm K}$ its modular measure before choosing a basis. Since Kato's construction starts from Siegel units on the open modular curve, open cohomology and its Drinfeld--Manin quotient are the natural candidates for this realization. The generic cuspidal identification determines a comparison over $L((X))$, and the question is whether, and how, it extends across the Eisenstein point.

\begin{question}[Comparison of Kato and compact periods]\label{conj:kato-period-comparison}
Does the generic cuspidal comparison extend to a map induced by the chosen cohomological realizations
\[
 j_{\mathrm K}:\mathscr P_{\mathbf f,\mathrm{cpt}}
       \longrightarrow\mathscr P_{\mathbf f}^{\mathrm K},
\]
for which the module-valued measures satisfy
\[
 \boldsymbol{\mathscr L}_{\mathbf f}^{\mathrm K}
   =(j_{\mathrm K}\otimes1)\boldsymbol{\mathscr L}_{\mathbf f,\mathrm{cpt}}\,?
\]
\end{question}

Both modules and their measures are fixed by their constructions. If this comparison holds, write
\[
 j_{\mathrm K}(\gamma_{\mathbf f})=X^{m_{\mathrm K}}v_{\mathrm K}(X)\gamma_{\mathbf f}^{\mathrm K}.
\]
Then $m_{\mathrm K}$ is the length of the cokernel, and \eqref{eq:factorization-Kato-period-input} gives the leading-term formula with $m=m_{\mathrm K}$ and $v=v_{\mathrm K}$. The auxiliary Drinfeld--Manin computation suggests a first derivative if this is the realization selected by Kato's periods, but one needs to establish the identification, including that of the module-valued measures, before assigning this order to Kato's $p$-adic $L$-function.

\subsection{Application to the circular-unit coefficient}
\label{subsec:factorization-application-kato}

We finally apply the compact factorization to the Kato class of Section~\ref{subsec:kato-weight-one}. Recall that its reciprocity law involves the auxiliary value, the adjoint factor and the comparison with compact periods. Treating these three factors together, we obtain the coefficient of the circular-unit class without the simple-zero assumption of Section~\ref{subsec:factorization-kato}.

Fix the trivial Teichm\"uller component and a nonzero cyclotomic denominator $d$, independent of $X$, which suffices for the regulator. Set
\[
 D=(e_0\mathcal H_\Gamma)[1/d],\qquad \mathcal R=D[[X]].
\]
The ring $D$ carries the cyclotomic variable, while $X$ varies the weight. Then $\mathcal R$ is $X$-torsion-free and $\mathcal R/X\mathcal R=D$ is a domain, and working in this completion allows us to compare orders in the weight variable without choosing a cyclotomic moment.

Part~(iii) of Hypothesis~\ref{hyp:kato-specialization}, together with \eqref{eq:generic-compact-period-comparison}, compares the generic modular period line with the compact one. Gathering all the weight-dependent factors, including this comparison, into a single scalar, we may write
\begin{equation}\label{eq:full-kato-multiplier}
 F(X):=\mathcal L_{\hf}^{-}(\operatorname{loc}_p\kappa(\hf))
   =\mathscr E_N^\circ(X,j)\,r(X)\,\mathscr L_p^{\mathrm{cpt}}(\mathbf f)(X),
 \qquad r\in L((X)).
\end{equation}
Here the scalar $r(X)$ is independent of $j$ by Hypothesis~\ref{hyp:kato-specialization}(iii), and $\mathscr E_N^\circ\in\mathcal R$ contains the retained tame Euler factors, normalized so that $\mathscr E_N^\circ(0,j)$ is not identically zero (individual cyclotomic zeros are allowed); its common power of $X$ has been included in $r$. More concretely, if $\mathscr L_p(\mathbf f)=q(X)\mathscr L_p^{\mathrm{cpt}}(\mathbf f)$ and $\mathscr E_N=\mathscr E_N^\circ e(X)$ in the normalization fixed before, then
\[
 r(X)=e(X)q(X)
       \frac{\mathscr L_p^{\mathrm{aux}}(\mathbf f)}
            {L_p(\operatorname{Ad}\mathbf f)}.
\]
In the analytic period normalization of \S\ref{subsec:factorization-kato}, $q$ has no pole: analyticity of $\mathscr L_p(\mathbf f)$ and the nonzero compact specialization in \eqref{eq:generic-compact-period-comparison} imply this directly. The order of $r$, however, also involves $e$ and the auxiliary-to-adjoint quotient, so we consider the product as a whole.

For a generic class in $\mathcal M[1/X]$ satisfying Hypothesis~\ref{hyp:kato-specialization}(ii), criterion~\eqref{eq:kato-regularity-test} identifies regularity of the class with $\operatorname{ord}_Xr\geq0$, because the other two factors have nonzero product at $X=0$. In Proposition~\ref{thm:kato-leading-term-weight-one}, this regularity followed from the simple-zero normalization. Here we instead assume regularity of the class as in part~(i), use it to specialize the full multiplier $r$, and apply the compact factorization to compute the resulting circular-unit coefficient.

\begin{corollary}[Kato and circular units]\label{lem:full-multiplier-regular}\label{cor:conditional-kato-circular-comparison}
Assume \eqref{eq:factorization-conductor-hypotheses}, $\alpha\ne1$, and Hypothesis~\ref{hyp:kato-specialization}. Then the function $r$ in \eqref{eq:full-kato-multiplier} is regular at $X=0$, and on the trivial Teichm\"uller component,
\begin{equation}\label{eq:conditional-kato-circular-comparison}
 \kappa(f)_{\chi_1}
 =\frac{r(0)\mathscr E_N^\circ(0,j)C_{f,\gamma}}
        {a_f}
       \langle M\rangle^{-j}
       L_p^{\mathrm{KL}}(\chi_2^{-1}\omega,-j)c(\chi_1)
\end{equation}
over $\mathcal K_\Gamma$. If the coefficient is analytic, the equality holds in $\mathcal S_{\chi_1}$. The specialized class is nonzero if and only if $r(0)\ne0$.
\end{corollary}
\begin{proof}
Parts~(i) and~(ii) of Hypothesis~\ref{hyp:kato-specialization} make the scalar regulator $F$ regular. To see what this implies for $r$, put $G=\mathscr E_N^\circ\mathscr L_p^{\mathrm{cpt}}(\mathbf f)$, so that $F=rG$ and $G(0)\ne0$ by Theorem~\ref{thm:compact-weight-one-factorization} and the normalization of $\mathscr E_N^\circ$. If $r=X^nv$ with $v\in A^\times$ and $n<0$, then reducing $X^{-n}F=vG$ modulo $X$ gives $0=v(0)G(0)$, a contradiction. Thus $r\in A$, with the case $r=0$ immediate. We may now specialize \eqref{eq:full-kato-multiplier}; using \eqref{eq:regular-kato-specialization}, we obtain
\[
 a_f\operatorname{Col}_{\chi_1}(\operatorname{res}_p\kappa(f)_{\chi_1})
    =r(0)\mathscr E_N^\circ(0,j)\mathscr L_p^{\mathrm{cpt}}(f,j).
\]
Substituting the compact factorization \eqref{eq:weight-one-bd-factorization} expresses the right-hand side as a product of the two Dirichlet functions. By \eqref{eq:reciprocity-circular} and its finite-character convention, the first is exactly the regulator of $c(\chi_1)$ in the same cyclotomic variable. It is nonzero, so it can be cancelled over $\mathcal K_\Gamma$, and the two sides of \eqref{eq:conditional-kato-circular-comparison} have the same global regulator. Proposition~\ref{ass:global-rank-one} then identifies the classes. If the coefficient is analytic, both sides lie in $\mathcal S_{\chi_1}$ and torsion-freeness gives the equality there; the nonvanishing of the factors and of the circular-unit regulator gives the final assertion.
\end{proof}

Although $r$ and the periods depend on the chosen bases, the coefficient of the circular-unit class does not. If $\gamma_{\mathbf f}$ is replaced by $u(X)\gamma_{\mathbf f}$, then $\mathscr L_p^{\mathrm{cpt}}(\mathbf f)$ gets multiplied by $u^{-1}$ and $r$ by $u$, so the product $r(0)C_{f,\gamma}$ is unchanged. Similarly, changing the regular Dieudonn\'e basis rescales $F$, $r$ and $a_f$ compatibly, and the coefficient in \eqref{eq:conditional-kato-circular-comparison} is again unchanged. On the other hand, a leading coefficient of order $m>0$ also depends on the choice of parameter: replacing $X$ by $X'=u(X)X$ multiplies $\left.X^{-m}\mathscr L\right|_0$ by $u(0)^{-m}$.

We close with two comments on the normalization. First, the auxiliary and adjoint factors come from the Kato class of \cite[Definition~4.6 and Proposition~4.9]{PR25}. For a zeta element whose reciprocity law involves a single modular $p$-adic $L$-function, the multiplier would instead consist of the period comparison and the retained Euler factors. Comparing both normalizations at $X=0$ requires transporting the classes before dividing by the auxiliary value, including their regularity in the split realization.

Second, suppose that the period map of Question~\ref{conj:kato-period-comparison} has order one and that no other factor contributes a zero or a pole. Then the injectivity of the regulator forces the class to vanish at weight one, and a comparison for its weight derivative would require control of the division by $X$ and of base change. The auxiliary Drinfeld--Manin computation suggests this possibility, although it does not determine the order for Kato's realization.

It remains to determine whether the specialized coefficient is nonzero, which amounts to computing the order of $r$. If the factors are nonzero, the definition of $r$ gives
\[
 \operatorname{ord}_Xr
   =\operatorname{ord}_Xe+\operatorname{ord}_Xq
       +\operatorname{ord}_X\mathscr L_p^{\mathrm{aux}}(\mathbf f)
       -\operatorname{ord}_XL_p(\operatorname{Ad}\mathbf f).
\]
Corollary~\ref{cor:conditional-kato-circular-comparison} shows that this integer is nonnegative, and $r(0)\ne0$ precisely when it is zero. Hence the remaining problem is to determine the balance between the auxiliary, adjoint and period orders in the fixed Kato normalization. The compact factorization supplies the nonzero weight-one value, but it does not determine these orders separately.

\section{Heegner cycles and elliptic units}
\label{sec:heegner-elliptic}

Having compared Kato classes with circular units in the cyclotomic setting, we turn to elliptic units over an imaginary quadratic field. The corresponding geometric classes are Heegner classes, and the comparison takes place on an anticyclotomic translate of Katz's character space. Their localization already lies in the positive ordinary direction, so we can use its regulator directly, as in \cite[\S C3]{LR1}.

At weight one, the two Eisenstein constituents give two elliptic-unit classes. The factorization of the $p$-adic $L$-function suggests multiplying each class by the Katz function of the other constituent. In order to compare this combination with the specialized Heegner class, we first need to check that it satisfies the same local condition at $p$.

\subsection{Katz $p$-adic $L$-functions and elliptic units}
\label{subsec:katz-elliptic-units}

Let $K$ be an imaginary quadratic field in which
\[
      p\cO_K=\mathfrak p\bar{\mathfrak p}
\]
splits, where $\mathfrak p$ is determined by our fixed embedding $\overline{\mathbb Q}\hookrightarrow\overline{\mathbb Q}_p$. We assume the Heegner hypothesis of \cite[\S2.1]{JLZ}, namely that every prime dividing the tame level $N$ splits in $K$, and choose an ideal $\mathfrak N\subseteq\cO_K$ such that $\cO_K/\mathfrak N\simeq\mathbb Z/N\mathbb Z$. For simplicity of exposition, we shall assume that $K$ has class number one and that $K\neq\mathbb Q(\sqrt{-1}),\mathbb Q(\sqrt{-3})$, although neither assumption is essential for the argument.

Let $\sigma$ and $\bar\sigma$ denote the two embeddings $K\hookrightarrow\overline{\mathbb Q}_p$, ordered so that $\sigma$ corresponds to $\mathfrak p$. We follow the conventions of \cite[\S 2.1]{JLZ} for infinity types. Thus, if $a-b\equiv0\pmod{w_K}$, where $w_K=\#\cO_K^\times$, the character $\sigma^a\bar\sigma^b$ determines an algebraic Gr\"ossencharacter of conductor one and infinity type $(a,b)$.

Let
\[
  \Gamma_K=\operatorname{Gal}(K(\mathfrak m p^\infty)/K),
  \qquad
  \mathcal W_K=\operatorname{Hom}_{\mathrm{cts}}
       (\Gamma_K,\mathbb G_m),
\]
where $\mathfrak m$ is a fixed ideal prime to $p$, divisible by $\mathfrak N$ and by the conductors of the finite-order Hecke characters used below, and $K(\mathfrak m p^\infty)$ is the union of the ray class fields of modulus $\mathfrak m p^n$. We work componentwise, fixing the finite characters. The subspace of characters unramified away from $p$ is the space of \cite[\S B4.2]{LR1}, and the enlarged ray modulus allows us to consider the nontrivial tame part of the character $\xi_0$ introduced below. Each component has two $p$-adic directions; the anticyclotomic quotient used below has one. Write $\mathcal H(\Gamma_K)=\cO(\mathcal W_K)$ for its algebra of locally analytic distributions. We denote by
\[
  \mathbf j_K:G_K\longrightarrow\mathcal H(\Gamma_K)^\times
\]
the universal character. Our reciprocity-map conventions are those of \cite[\S2.3.2]{JLZ}: geometric Frobenius corresponds to a uniformizer. For a Hecke character $\xi$, write $\widehat\xi=V_p(\xi)$ for its Galois avatar; in particular $V_p(\mathrm{Nm})=L(-1)$. Following \cite[\S B4.2]{LR1}, the point of character space labelled by $\xi$ is the point at which $\mathbf j_K$ specializes to $\widehat\xi^{-1}$. This inverse is part of the Mellin convention, and hence the specialization of $V_p(\Psi)^*(1-\mathbf j_K)$ at $\xi$ is $V_p(\Psi)^*(1)\otimes\widehat\xi$.

Let $\Psi$ be a finite-order Gr\"ossencharacter of $K$, of conductor prime to $p$, with values in a finite extension $L/\mathbb Q_p$. We use Katz's two-variable $p$-adic $L$-function and the elliptic-unit system with the fixed ray modulus $\mathfrak m$, retaining the local factors at its primes. This gives an analytic function
\[
  L_{\mathfrak p}^{\mathrm{Katz}}(\Psi)
     \in \mathcal H(\Gamma_K)
          \widehat\otimes_{\mathbb Q_p}\widehat{\mathbb Q}_p^{\mathrm{nr}}.
\]
Its interpolation formula, in the range relevant to us, takes the following form. If $\xi$ is a character of conductor one and infinity type $(a,b)$, with $a\geq1$ and $b\leq0$, then
\begin{equation}\label{eq:katz-interpolation}
 \frac{L_{\mathfrak p}^{\mathrm{Katz}}(\Psi)(\xi)}
      {\Omega_p^{a-b}}
  = \mathfrak a(\xi)\mathfrak e_p(\Psi,\xi)
       \mathfrak d_K(\xi)
       \frac{L(\Psi\xi^{-1},0)}{\Omega^{a-b}},
\end{equation}
where
\[
 \mathfrak a(\xi)=(a-1)!\pi^{-b},
 \qquad
 \mathfrak e_p(\Psi,\xi)
   =\bigl(1-p^{-1}\Psi^{-1}\xi(\mathfrak p)\bigr)
      \bigl(1-\xi^{-1}\Psi(\bar{\mathfrak p})\bigr),
\]
and $\mathfrak d_K(\xi)=D_K^{b/2}2^{-b}$. Here $D_K$ is the absolute discriminant of $K$, and $\Omega$ and $\Omega_p$ are compatible complex and $p$-adic CM periods. We refer the reader to \cite{katz76} and \cite[Theorem 3.1]{BDP12} for the precise normalization. When the fixed tame part is nontrivial, the corresponding local factors must also be included, the conductor-one formula above being only the reference normalization. For both the Katz and the Heegner comparisons, we use a common coefficient extension containing the chosen CM periods.

Observe that the finite part of the anticyclotomic translate is determined by the Heegner datum, and it cannot be chosen independently. Write $\epsilon_f=\chi_1\chi_2$ for the odd nebentype of the weight-one form, and choose a Gr\"ossencharacter $\lambda_0$ of finite type $(1,\mathfrak N,\epsilon_f)$ and infinity type $(0,-1)$ in the convention of \cite[Definition~2.1.1]{JLZ}. Its finite-type restriction is the one prescribed there; existence follows from \cite[Remark~2.1.2]{JLZ}, since $\epsilon_f(-1)=-1$ and $\cO_K^\times=\{\pm1\}$. Put $\xi_0=\mathrm{Nm}\,\lambda_0$, of infinity type $(1,0)$, where $\mathrm{Nm}$ is the Hecke norm character of type $(1,1)$. Let $\Gamma^{\mathrm{ac}}$ be the Galois group of the anticyclotomic $\Z_p$-extension of $K$, and let $\mathcal W^{\mathrm{ac}}$ be its character space, with the finite component fixed as above. We use the translate
\[
 \iota_{\mathrm{ac}}:\mathcal W^{\mathrm{ac}}\longrightarrow\mathcal W_K,
 \qquad j\longmapsto\xi_j:=\xi_0(\sigma/\bar\sigma)^j.
\]
Its infinity type is $(1+j,-j)$, and the Katz interpolation range contains the integers $j\geq0$ of the fixed finite component. The norm factor, together with the inverse Mellin convention that we have just specified, supplies the Tate twist needed to match the elliptic-unit coefficients with those of the Heegner classes (see \eqref{eq:heegner-self-duality} below). Since the character $\xi_0$ has a nontrivial finite part, the local factors of that finite part must be included in the interpolation comparison. We write $\mathcal H_{\mathrm{ac}}$ for the chosen component of $\mathcal H(\Gamma^{\mathrm{ac}})$ and $\mathcal K_{\mathrm{ac}}$ for its fraction field. Pullbacks along $\iota_{\mathrm{ac}}$ will often be left implicit.

The Euler system of elliptic units gives a class
\begin{equation}\label{eq:elliptic-unit-big-class}
 \kappa(\Psi,K)
   \in H^1\bigl(K,V_p(\Psi)^*(1-\mathbf j_K)\bigr)^{\sim},
\end{equation}
where $(-)^{\sim}$ denotes the reflexive hull. At finite-order characters of the full two-variable space, this class specializes to the Kummer image of a suitable linear combination of units or $S$-units in an abelian extension of $K$, where $S$ contains the conductor primes. We refer the reader to \cite[\S 15]{Kato} or \cite[\S 1.2]{BCDDPR} for an account of the construction.

The two-variable Perrin--Riou regulator of \cite{loefflerzerbes14} gives a Coleman map
\[
 \operatorname{Col}_{\mathfrak p,\Psi}:
 H^1\bigl(K_{\mathfrak p},V_p(\Psi)^*(1-\mathbf j_K)\bigr)^{\sim}
    \longrightarrow\mathcal H(\Gamma_K),
\]
and the explicit reciprocity law for elliptic units reads
\begin{equation}\label{eq:elliptic-unit-reciprocity}
 \operatorname{Col}_{\mathfrak p,\Psi}
   \bigl(\operatorname{res}_{\mathfrak p}\kappa(\Psi,K)\bigr)
    =L_{\mathfrak p}^{\mathrm{Katz}}(\Psi).
\end{equation}
As usual, this equality incorporates the local Euler factors into the normalization of the Katz function.

\subsection{Heegner classes in families and their local projection}
\label{subsec:heegner-families-local}

Let
\[
       f=E_1(\chi_1,\chi_2)
\]
be the irregular weight-one Eisenstein series fixed above, with nebentype $\epsilon_f=\chi_1\chi_2$. Let $\Psi_i$ be the Hecke base change of $\chi_i$ to $K$. In the geometric-Frobenius convention of \cite[\S B4.1]{LR1}, we have $V_p(\Psi_i)=L(\chi_i^{-1})|_{G_K}$ and $V_p(\Psi_i)^*=L(\chi_i)|_{G_K}$. Note that the Heegner hypothesis rules out the case $\epsilon_f=\varepsilon_K$, since then the ramified primes of $K$ would divide $N$; instead, the finite part of $\lambda_j$ accounts for the nebentype in the self-duality relation below.

Let $\mathbf f$ be the chosen generically cuspidal ordinary family on $U$. We work on its inverse image in $\mathcal E_K(N,\epsilon_f)$, the space of conjugate-self-dual twists of \cite[Definition~2.4.2]{JLZ}; this construction allows for a nontrivial nebentype, and does not require a square root of the weight character. On the weight-one fibre take $\lambda_j=\xi_j\,\mathrm{Nm}^{-1}$, of infinity type $(j,-1-j)$. Write hats for the associated $p$-adic Galois characters. The finite-type condition gives
\begin{equation}\label{eq:heegner-self-duality}
 \widehat\lambda_j\widehat\lambda_j^{c}
      =\epsilon_f^{-1}|_{G_K}\,\varepsilon,
 \qquad
 \widehat\lambda_j=\varepsilon\widehat\xi_j.
\end{equation}
The second identity follows from $V_p(\mathrm{Nm})=L(-1)$. Once the split coefficient realization has been fixed, it gives
\[
 V_f\otimes\widehat\lambda_j
   \simeq\bigoplus_{i=1}^2
        V_p(\Psi_i)^*(1)\otimes\widehat\xi_j.
\]
The right-hand side consists precisely of the specializations of the elliptic-unit coefficients $V_p(\Psi_i)^*(1-\mathbf j_K)$. On the other hand, in \cite[\S2.3]{JLZ} the representation $V_p(f)$ is defined using geometric Frobenius, and it agrees with our $V_f^*$, so their Heegner coefficients $V_p(f)^*\otimes\widehat\lambda_j$ form the left-hand side. In this way the dual, the Tate twist and the inverse Mellin convention are fixed simultaneously.

Let $\mathbb V_{\mathbf f,K}^{\mathrm{ac}}$ denote the resulting two-variable Galois representation over $\cO(U)\widehat\otimes\mathcal H(\Gamma^{\mathrm{ac}})$ after choosing a local section of the twist space through $\lambda_0$, normalized to be conjugate-self-dual. The construction of \cite{JLZ} gives a big Heegner class
\begin{equation}\label{eq:big-heegner-class}
 \kappa(\mathbf f,K)
    \in H^1\bigl(K,\mathbb V_{\mathbf f,K}^{\mathrm{ac}}\bigr),
\end{equation}
whose arithmetic specializations are the Abel--Jacobi images of generalized Heegner cycles, up to the explicit $p$-stabilization factors appearing in \emph{loc.\,cit.}

Writing $\mathbb D_{\mathbf f,K}^{\mathrm{ac}}$ for the local $(\varphi,\Gamma)$-module of $\mathbb V_{\mathbf f,K}^{\mathrm{ac}}$, the ordinary filtration of $\mathbb V_{\mathbf f}$ induces an exact sequence of trianguline $(\varphi,\Gamma)$-modules at $\mathfrak p$,
\[
 0\longrightarrow
 \mathscr F^+\mathbb D_{\mathbf f,K}^{\mathrm{ac}}
 \longrightarrow
 \mathbb D_{\mathbf f,K}^{\mathrm{ac}}
 \longrightarrow
 \mathscr F^-\mathbb D_{\mathbf f,K}^{\mathrm{ac}}
 \longrightarrow0.
\]
The family construction and its reciprocity law are available on the cuspidal locus. At the Eisenstein point we need the following input, which we keep for the rest of the section.

\begin{hyp}[Regular specialization of Heegner classes]\label{hyp:heegner-specialization}
The family, its finite character components and its self-dual twists lie in the range of the construction and reciprocity law of \cite{JLZ} on the cuspidal locus. Moreover, the Heegner class, its positive local direction and its positive regulator extend regularly to $x_f$ in the chosen realization, whose fibre splits as in \eqref{eq:eisenstein-ac-decomposition} below.
\end{hyp}

We fix the nonzero scalar $\Omega_{f,\gamma}$ comparing the specialized regulator with the regulator in the chosen weight-one bases. Hypothesis~\ref{hyp:heegner-specialization} allows the family local condition to pass to weight one.

\begin{proposition}\label{prop:heegner-natural-projection-zero}
Under Hypothesis~\ref{hyp:heegner-specialization}, the localization of $\kappa(\mathbf f,K)$ at $\mathfrak p$ belongs to the image of
\[
 H^1\bigl(K_{\mathfrak p},
       \mathscr F^+\mathbb D_{\mathbf f,K}^{\mathrm{ac}}\bigr)
 \longrightarrow
 H^1\bigl(K_{\mathfrak p},
       \mathbb D_{\mathbf f,K}^{\mathrm{ac}}\bigr).
\]
Equivalently, its image under the natural projection
\[
 H^1\bigl(K_{\mathfrak p},
       \mathbb D_{\mathbf f,K}^{\mathrm{ac}}\bigr)
 \longrightarrow
 H^1\bigl(K_{\mathfrak p},
       \mathscr F^-\mathbb D_{\mathbf f,K}^{\mathrm{ac}}\bigr)
\]
is zero.
\end{proposition}

\begin{proof}
At the crystalline classical points in the interpolation range, the Heegner classes are images of algebraic cycles under the $p$-adic Abel--Jacobi map, and hence they satisfy the Bloch--Kato local condition; in particular, their projections to the negative rank-one quotient vanish. Since these points are Zariski-dense and the cohomology of the negative quotient is torsion-free over the coefficient algebra, the assertion follows by rigid-analytic continuation. This is precisely \cite[Theorem 6.3.4]{JLZ}; see also \cite[Proposition C3.3]{LR1}.
\end{proof}

Hence we evaluate the Heegner class by means of the regulator on the positive piece, as in \cite[\S C3.4]{LR1}. Denote the specialized coefficient module by $\mathbb V_{f,K}^{\mathrm{ac}}$ and the Heegner Selmer condition by $\mathcal F_{\mathrm{Hg}}$: it uses the specialized positive direction at $\mathfrak p$, the conjugate condition at $\bar{\mathfrak p}$, and the chosen unramified or conductor-allowed conditions away from $p$. By Hypothesis~\ref{hyp:heegner-specialization}, the fixed period isomorphism $\gamma$ gives a class at $x_f$
\[
 \kappa(f,K)\in H^1_{\mathcal F_{\mathrm{Hg}}}
       \bigl(K,\mathbb V_{f,K}^{\mathrm{ac}}\bigr)
\]
and a regulator
\[
 \mathcal L_{f,K}^{+}:
 H^1\bigl(K_{\mathfrak p},
      \mathscr F^+\mathbb V_{f,K}^{\mathrm{ac}}\bigr)
   \longrightarrow\mathcal H(\Gamma^{\mathrm{ac}}).
\]
The explicit reciprocity law of \cite{JLZ} takes the form
\begin{equation}\label{eq:heegner-reciprocity}
 \mathcal L_{f,K}^{+}
   \bigl(\operatorname{res}_{\mathfrak p}\kappa(f,K)\bigr)
   =\Omega_{f,\gamma}^{-1}
      L_{\mathfrak p}^{\mathrm{BDP}}(f),
\end{equation}
where $L_{\mathfrak p}^{\mathrm{BDP}}(f)$ denotes the weight-one specialization of the family BDP $p$-adic $L$-function in the reciprocity law, and $\Omega_{f,\gamma}\in L^\times$ measures the discrepancy between the specialization of the family regulator and our chosen regulator at $f$.

\subsection{The BDP--Katz factorization}
\label{subsec:bdp-katz-factorization}

The Heegner regulator gives a BDP function, whereas the elliptic-unit regulators give Katz functions. The splitting of $V_f$ suggests comparing the square of the BDP function with the product of the two Katz functions. Put
\[
 \mathscr L_i
   :=L_{\mathfrak p}^{\mathrm{Katz}}(\Psi_i)
        \big|_{\iota_{\mathrm{ac}}(\mathcal W^{\mathrm{ac}})},
 \qquad i=1,2,
\]
and
\[
 \mathscr L_{\mathrm{BDP}}
   :=L_{\mathfrak p}^{\mathrm{BDP}}(f)
        \big|_{\mathcal W^{\mathrm{ac}}}.
\]

The norm shift in our choice of translate matches the complex arguments. At the relevant arithmetic characters, Artin formalism gives
\[
 L(f/K,\lambda_j^{-1},1)
   =\prod_{i=1}^2L(\Psi_i\lambda_j^{-1},1)
   =\prod_{i=1}^2L(\Psi_i\xi_j^{-1},0),
\]
since $\lambda_j^{-1}=\mathrm{Nm}\,\xi_j^{-1}$ and the Hecke norm shifts $s$ to $s-1$. Thus the complex values match those in Katz's interpolation formula. The remaining comparison concerns the interpolation factors and periods.

We keep Hypothesis~\ref{hyp:heegner-specialization}. For the comparison of analytic functions we need, in addition to the family reciprocity law, compatible interpolation formulas at weight one; compare \cite[\S C3.5]{LR1} for the critical-slope counterpart.

\begin{hyp}[Interpolation compatibility]\label{hyp:heegner-interpolation}
The specialized BDP function and the two Katz functions admit compatible interpolation formulas on a common set $\Sigma$ of arithmetic characters which determines meromorphic functions on the chosen component. More precisely, the formulas refer to the Rankin motive above, with compatible algebraic periods, and, writing $C_i(\nu)$ for the resulting algebraic degree-one values, they take the form
\[
 \mathscr L_{\mathrm{BDP}}(\nu)^2
     =b(\nu)C_1(\nu)C_2(\nu),\qquad
 \mathscr L_i(\nu)=k_i(\nu)C_i(\nu)
 \quad(\nu\in\Sigma).
\]
Here $b$ and $k_i$ record the local interpolation factors and the period normalizations; we assume that $k_1(\nu)k_2(\nu)\ne0$ for $\nu\in\Sigma$, and that the quotient $b/(k_1k_2)$ extends to a specified meromorphic function $\mathscr A(f,K)$ on the chosen component.
\end{hyp}

\begin{proposition}[BDP--Katz comparison]\label{prop:bdp-katz-factorization}
Under Hypothesis~\ref{hyp:heegner-interpolation}, we have
\begin{equation}\label{eq:bdp-katz-factorization}
   \mathscr L_{\mathrm{BDP}}^2
      =\mathscr A(f,K)\,\mathscr L_1\mathscr L_2.
\end{equation}
\end{proposition}

\begin{proof}
At every $\nu\in\Sigma$, the interpolation formulas give $\mathscr L_{\mathrm{BDP}}(\nu)^2=(b/k_1k_2)(\nu)\,\mathscr L_1(\nu)\mathscr L_2(\nu)$, where we use that $k_1(\nu)k_2(\nu)\ne0$; note that no division by the algebraic $L$-values $C_i(\nu)$ is needed. Since $\Sigma$ determines meromorphic functions on the component, the identity follows.
\end{proof}

We keep $\mathscr A(f,K)$ meromorphic unless the local factors and the period comparison in both numerator and denominator are analytic units; in that case it is an analytic unit. In particular, specifying this factor includes the weight-one CM evaluation and its periods.

The two Katz characters are related by the functional equation. Indeed, translating \eqref{eq:heegner-self-duality} back to Hecke characters gives
\[
 \lambda_j\lambda_j^c=\Psi_1\Psi_2\,\mathrm{Nm}^{-1},
 \qquad
 \xi_j\xi_j^c=\Psi_1\Psi_2\,\mathrm{Nm}.
\]
Since $\Psi_i^c=\Psi_i$, the characters $A_i=\Psi_i\lambda_j^{-1}$ and $\Theta_i=\Psi_i\xi_j^{-1}$ satisfy
\begin{equation}\label{eq:katz-functional-equation-characters}
 A_2=A_1^{-c}\mathrm{Nm},
 \qquad
 \Theta_2=\Theta_1^{-c}\mathrm{Nm}^{-1}.
\end{equation}
Hence the values at $s=1$, respectively at $s=0$, form a pair related by the functional equation. This relation uses the finite part of the self-dual twist, and it holds without imposing $\Psi_2=\Psi_1^{-1}$. However, in order to replace $\mathscr L_1\mathscr L_2$ by a square in a fixed $p$-adic normalization, one must also compare the CM periods and the local factors of the Katz functional equation on this translate, and these factors may have zeros or poles. For this reason we keep the product in \eqref{eq:bdp-katz-factorization}, which is the form directly compatible with the two elliptic-unit reciprocity laws.

\subsection{The weighted elliptic-unit class}
\label{subsec:weighted-elliptic-unit}

Let
\[
 \mathbb V_i^{\mathrm{ac}}
   :=\iota_{\mathrm{ac}}^*
        \bigl(V_p(\Psi_i)^*(1-\mathbf j_K)\bigr),
 \qquad i=1,2,
\]
and denote by
\[
 \mathfrak u_i
   :=\iota_{\mathrm{ac}}^*\kappa(\Psi_i,K)
      \in H^1(K,\mathbb V_i^{\mathrm{ac}})^{\sim}
\]
the corresponding anticyclotomic elliptic-unit classes. After fixing the specialization isomorphism $\gamma$ at the weight-one point, there is a decomposition
\begin{equation}\label{eq:eisenstein-ac-decomposition}
 \mathbb V_{f,K}^{\mathrm{ac}}
   \simeq
 \bigl(Lv_{f,\chi_1}\otimes\mathbb V_1^{\mathrm{ac}}\bigr)
 \oplus
 \bigl(Lv_{f,\chi_2}\otimes\mathbb V_2^{\mathrm{ac}}\bigr).
\end{equation}
In \eqref{eq:eisenstein-ac-decomposition}, the vectors $v_{f,\chi_i}$ just label the two summands, since the character action is already included in $\mathbb V_i^{\mathrm{ac}}$. We normalize these vectors compatibly with $\gamma$; changing this normalization multiplies both sides of the final comparison by the same nonzero scalar.

\begin{defi}\label{def:weighted-elliptic-unit}
The \emph{weighted elliptic-unit class} associated with $(f,K)$ is
\begin{equation}\label{eq:weighted-elliptic-unit}
 \mathfrak u(f,K)
   :=v_{f,\chi_1}\otimes\mathscr L_2\mathfrak u_1
      +v_{f,\chi_2}\otimes\mathscr L_1\mathfrak u_2
 \in H^1\bigl(K,\mathbb V_{f,K}^{\mathrm{ac}}\bigr)^{\sim}.
\end{equation}
\end{defi}

By the elliptic-unit reciprocity laws, both weighted summands have regulator $\mathscr L_1\mathscr L_2$. In order to turn this equality into the Heegner local condition, we must identify the positive line selected by the family with the line of equal coordinates in these Coleman normalizations. This compatibility is not determined by the equality of the two local characters, so we retain it as an explicit assumption for the rest of the section.

\begin{ass}[Compatibility with the ordinary local condition]\label{ass:elliptic-local-injectivity}
The specialized positive line has nonzero projection to each of the two character lines in \eqref{eq:eisenstein-ac-decomposition}. Let $Q_i$ be the local Iwasawa cohomology module of $\mathbb V_i^{\mathrm{ac}}$ at $\mathfrak p$. These modules are free of rank one, and there are isomorphisms $a_i:Q_i\to Q$ to a common rank-one module such that the Coleman maps are $\operatorname{Col}_i=\operatorname{Col}\circ a_i$, with $\operatorname{Col}$ injective. Under \eqref{eq:eisenstein-ac-decomposition}, the specialized positive local condition is
\[
 \{(x_1,x_2):a_1(x_1)=a_2(x_2)\},
\]
and its regulator sends such a pair to the common Coleman value $\operatorname{Col}(a_1(x_1))=\operatorname{Col}(a_2(x_2))$. The local conditions at $\bar{\mathfrak p}$ are compatible with this description under conjugation, including the finite twists and the weights of the elliptic-unit classes. Away from $p$, the classes satisfy the chosen unramified or conductor-allowed Selmer conditions.
\end{ass}

\begin{proposition}\label{prop:weighted-unit-local-condition}
Under Hypothesis~\ref{hyp:heegner-specialization} and Assumption~\ref{ass:elliptic-local-injectivity}, the class $\mathfrak u(f,K)$ satisfies the same local conditions as the specialized Heegner class $\kappa(f,K)$. In particular,
\[
 \mathfrak u(f,K)
    \in H^1_{\mathcal F_{\mathrm{Hg}}}
       \bigl(K,\mathbb V_{f,K}^{\mathrm{ac}}\bigr)^{\sim}.
\]
\end{proposition}

\begin{proof}
Away from $p$, the assertion is one of the stated local conditions, consistent with the ramification of the elliptic units at the conductor primes. At $\mathfrak p$, we use \eqref{eq:eisenstein-ac-decomposition} and the specialization of the ordinary filtration selected by the family $\mathbf f$. By the specified identifications $a_i$, applying the injective map $\operatorname{Col}$ to the difference that defines the local quotient gives
\[
 \mathscr L_2\,
 \operatorname{Col}_{\mathfrak p,\Psi_1}
       (\operatorname{res}_{\mathfrak p}\mathfrak u_1)
  -\mathscr L_1\,
 \operatorname{Col}_{\mathfrak p,\Psi_2}
       (\operatorname{res}_{\mathfrak p}\mathfrak u_2).
\]
By \eqref{eq:elliptic-unit-reciprocity}, this is
\[
       \mathscr L_2\mathscr L_1-
       \mathscr L_1\mathscr L_2=0.
\]
By the injectivity of the rank-one Coleman map, the corresponding local cohomology class vanishes. Finally, the condition at $\bar{\mathfrak p}$ follows from the conjugation compatibility included in the assumption, and hence $\mathfrak u(f,K)$ belongs to the Heegner Selmer group.
\end{proof}

\begin{remark}\label{rem:weighted-unit-artin}
Proposition~\ref{prop:weighted-unit-local-condition} is the analogue of the verification of the balanced local condition for the weighted Beilinson--Flach class in \cite[\S 5.3]{AOR}. Once the local identifications of Assumption~\ref{ass:elliptic-local-injectivity} are available, the analytic factorization determines a combination of rank-one Euler systems satisfying the Heegner local condition.
\end{remark}

\begin{example}\label{ex:common-conductor-heegner}
The conductor and Heegner conditions can be satisfied outside the CM case. For instance, let $M=19$, let $\zeta_{18}$ be a primitive eighteenth root of unity, and let $\vartheta$ be the character modulo $19$ with $\vartheta(2)=\zeta_{18}$. Take
\[
 \chi_1=\vartheta^2,\qquad \chi_2=\vartheta^5,\qquad p=83,
 \qquad K=\Q(\sqrt{-67}),\qquad N=19^2.
\]
The ratio $\chi_1\chi_2^{-1}=\vartheta^{-3}$ has order six and is primitive. We have $83\nmid19\cdot18$. Since $83\equiv7\equiv2^6\pmod{19}$, both characters take the value $\zeta_{18}^{12}\ne1$ at $83$. The field $K$ has class number one, and both $19$ and $83$ split in it: the relevant quadratic residues are $-67\equiv9\pmod{19}$ and $83\equiv16\pmod{67}$. Choosing a prime $\mathfrak n\mid19$ and $\mathfrak N=\mathfrak n^2$ gives $\cO_K/\mathfrak N\simeq\Z/19^2\Z$. The odd nebentype $\vartheta^7$ satisfies the finite-type parity condition for $\lambda_0$. These data therefore satisfy the conductor and Heegner conditions. The comparison of classes also requires regular specialization, compatible interpolation and local regulators, and the Selmer conditions in Theorem~\ref{thm:heegner-elliptic-comparison} below.
\end{example}

\subsection{Comparison of classes}
\label{subsec:heegner-elliptic-comparison}

Both classes now satisfy the same local conditions. Write
\[
 \mathcal S(f,K):=
 H^1_{\mathcal F_{\mathrm{Hg}}}
       \bigl(K,\mathbb V_{f,K}^{\mathrm{ac}}\bigr)^{\sim}.
\]
The elliptic-unit reciprocity laws give
\begin{equation}\label{eq:weighted-unit-reciprocity}
 \mathcal L_{f,K}^{+}
   \bigl(\operatorname{res}_{\mathfrak p}\mathfrak u(f,K)\bigr)
     =\mathscr L_1\mathscr L_2,
\end{equation}
where the positive regulator is evaluated on the common coordinate $a_1(x_1)=a_2(x_2)$ from Assumption~\ref{ass:elliptic-local-injectivity}. 

\begin{theorem}\label{thm:heegner-elliptic-comparison}
Assume Hypotheses~\ref{hyp:heegner-specialization} and~\ref{hyp:heegner-interpolation}, and Assumption~\ref{ass:elliptic-local-injectivity}. Suppose that $\mathcal S(f,K)$ is torsion-free of rank one over $\mathcal H_{\mathrm{ac}}$ and that $\mathcal L_{f,K}^{+}\circ\operatorname{res}_{\mathfrak p}$ is nonzero on it. Then, in $\mathcal S(f,K)\otimes_{\mathcal H_{\mathrm{ac}}} \mathcal K_{\mathrm{ac}}$, we have
\begin{equation}\label{eq:heegner-elliptic-comparison}
\mathscr A(f,K)\,\mathfrak u(f,K) =\Omega_{f,\gamma}\, \mathscr L_{\mathrm{BDP}}\,\kappa(f,K).
\end{equation}
If $\mathscr A(f,K)$ is analytic and both classes belong to $\mathcal S(f,K)$, the equality holds in $\mathcal S(f,K)$ itself. Note that this is a statement over the indicated coefficient algebra, after inverting $p$, and not over an integral $\cO_L$-lattice.
\end{theorem}

\begin{proof}
By \eqref{eq:weighted-unit-reciprocity} and Proposition~\ref{prop:bdp-katz-factorization}, the image of the left-hand side of \eqref{eq:heegner-elliptic-comparison} under the regulator is
\[
 \mathscr A(f,K)\mathscr L_1\mathscr L_2
     =\mathscr L_{\mathrm{BDP}}^2.
\]
On the other hand, the explicit reciprocity law \eqref{eq:heegner-reciprocity} gives
\[
 \mathcal L_{f,K}^{+}
  \bigl(\operatorname{res}_{\mathfrak p}
      (\Omega_{f,\gamma}\mathscr L_{\mathrm{BDP}}
          \kappa(f,K))\bigr)
   =\mathscr L_{\mathrm{BDP}}^2.
\]
Thus the two sides have the same image under $\mathcal L_{f,K}^{+}\circ\operatorname{res}_{\mathfrak p}$, which is injective because it is nonzero on a torsion-free module of rank one.
\end{proof}

The extra factor $\mathscr L_{\mathrm{BDP}}$ in \eqref{eq:heegner-elliptic-comparison} is therefore dictated by reciprocity: the Heegner regulator gives one BDP factor, while the weighted elliptic-unit regulator gives its square.

\begin{remark}\label{rem:comparison-dictionary}
Theorem~\ref{thm:heegner-elliptic-comparison} is the rank-one anticyclotomic analogue of the comparison theorem of \cite{AOR}. The dictionary is
\[
\begin{array}{c|c}
 \text{diagonal-cycle setting} & \text{present setting} \\
 \hline
 \text{triple-product }p\text{-adic }L\text{-function}
   & \text{BDP }p\text{-adic }L\text{-function} \\
 \text{two Hida--Rankin functions}
   & \text{two Katz functions} \\
 \text{weighted Beilinson--Flach class}
   & \text{weighted elliptic-unit class} \\
 \text{diagonal-cycle class}
   & \text{Heegner-cycle class}.
\end{array}
\]
In both cases, the factorization of the square of a larger $p$-adic $L$-function determines a weighted combination of two smaller Euler systems, in which each class is multiplied by the $p$-adic $L$-function of the complementary constituent. The resulting class is then compared with the geometric class, after multiplying the latter by one copy of the larger $p$-adic $L$-function.
\end{remark}

\begin{remark}\label{rem:self-dual-special-case}
If the functional equation identifies the two Katz factors, and there is also a compatible identification of the two coefficient lines and of their elliptic-unit classes, the formula reduces to a comparison with a single elliptic-unit system, and it then has the shape of \cite[Theorem C3.9]{LR1}. Observe, however, that the equality of the two scalar Katz functions is not enough to identify the two global cohomology summands.
\end{remark}

\section{Double Eisenstein degeneration}\label{sec:double-degeneration}

The preceding comparisons suggest asking what happens when two modular variables specialize to Eisenstein points. Starting with Beilinson--Flach classes, one would first obtain Kato classes and then seek a comparison with circular units as in Section~\ref{sec:circular-kato}; see Remark~\ref{rem:two-weight-one-limits} for why we do not pursue this route. Starting with diagonal cycles, the corresponding two steps lead to Beilinson--Flach and then Beilinson--Kato classes. We study this second route here, keeping the remaining modular family cuspidal.

The first comparison follows \cite{AOR}, and the second is modelled on \cite{BDV}. The analytic factorization at each step prescribes the weight of each character component. To combine the two comparisons, we need the first identity to specialize regularly at the second Eisenstein point, and the resulting classes to satisfy the local conditions used in the second comparison.

\subsection{The two Eisenstein specializations}\label{subsec:two-Eis-spec}

Let $\mathbf f$, $\mathbf g$ and $\mathbf h$ be three $p$-adic families of modular forms for which the diagonal-cycle class
\[
 \kappa(\mathbf f,\mathbf g,\mathbf h)
 \in
 H^1\bigl(\Q,\mathbb V_{\mathbf f\mathbf g\mathbf h}^{\dagger}\bigr)
\]
is defined.  Here, as usual, the superscript $\dagger$ incorporates the self-dual twist. We assume that the families $\mathbf f$ and $\mathbf h$ specialize at points $x_f$ and $x_h$ to irregular Eisenstein series of weight one, while the family $\mathbf g$ remains cuspidal throughout the two specializations.

For the first Eisenstein point, we use a CM family and write
\[
 f=\mathbf f_{x_f}=E_1(1,\varepsilon_K),
\]
where $K/\Q$ is imaginary quadratic, $\varepsilon_K$ is its quadratic character, and $p$ splits in $K$. Throughout this section, $f$ denotes this CM point; it is not the common-conductor form of Sections~\ref{sec:circular-kato}--\ref{sec:heegner-elliptic}. After enlarging $L$, the chosen CM realization gives
\[
 V_f=L\cdot v_{f,1}\oplus L\cdot v_{f,\varepsilon_K}.
\]
The vectors $v_{f,1}$ and $v_{f,\varepsilon_K}$ are fixed once and for all; all the periods appearing below depend on this choice.

For the second Eisenstein point, write
\[
 h=\mathbf h_{x_h}=E_1(\psi_1,\psi_2),
\]
so that
\[
 V_h=L\cdot v_{h,1}\oplus L\cdot v_{h,2},
 \qquad
 G_{\Q}\text{ acts on }L\cdot v_{h,j}\text{ through }\psi_j.
\]
For $h$ we assume the split realization displayed above and the local compatibility specified in \S\ref{subsec:second-degeneration}. Both character components are retained.

To specify the rings over which we specialize, choose connected affinoid neighbourhoods with algebras $A_{\mathbf f}$, $A_{\mathbf g}$ and $A_{\mathbf h}$ over $L$, and enlarge $L$ so that $x_f,x_h$ are rational. Fix the finite character components and the square root in the self-dual twist, passing to a finite \'etale cover if necessary. If $\mathfrak m_f$ and $\mathfrak m_h$ are the ideals of the two points, the successive algebras are
\[
\begin{aligned}
 R_3&=A_{\mathbf f}\widehat\otimes_L A_{\mathbf g}
                         \widehat\otimes_L A_{\mathbf h},\\
 R_f&=R_3/\mathfrak m_fR_3
       \simeq A_{\mathbf g}\widehat\otimes_L A_{\mathbf h},\\
 R_{f,h}&=R_f/\mathfrak m_hR_f\simeq A_{\mathbf g}.
\end{aligned}
\]
On compatible reduced irreducible components, write these last two domains as $\mathscr R_f$ and $\mathscr R_{f,h}$, with fraction fields $\mathscr K_f$ and $\mathscr K_{f,h}$. When no confusion can arise, we abbreviate them to $\mathscr R$ and $\mathscr K$ at the stage under discussion. Since the Rankin cyclotomic variable is restricted below to the central-character locus, the final identity only varies with $\mathbf g$.

We require the regularity at $x_h$ of the classes and comparison factors that we specialize, after clearing the stated denominators. Observe that the specialization is made in these modules before passing to $\mathscr K_{f,h}$, since a meromorphic identity alone cannot be evaluated at $x_h$; the resulting equality then holds in an unlocalized module whenever both sides belong to the same torsion-free module.

There are two successive decompositions of the specialized representation:
\[
 \left.\mathbb V_{\mathbf f\mathbf g\mathbf h}^{\dagger}
 \right|_{x_f}
 \simeq
 \mathbb V_{1}^{\dagger}\oplus
 \mathbb V_{\varepsilon_K}^{\dagger},
\]
and
\[
 \left.\mathbb V_i^{\dagger}\right|_{x_h}
 \simeq
 \mathbb V_{i,1}^{\dagger}\oplus
 \mathbb V_{i,2}^{\dagger},
 \qquad i\in\{1,\varepsilon_K\}.
\]
The four representations $\mathbb V_{i,j}^{\dagger}$ are rank-two twists of the representation attached to $\mathbf g$, and hence they are the natural coefficient modules for the four Beilinson--Kato classes appearing at the end of the double degeneration.

We work in the non-exceptional range of the regulators and retain their Euler factors. A vanishing factor that destroys injectivity would require a separate argument.

\subsection{From diagonal cycles to Beilinson--Flach classes}
\label{subsec:first-degeneration}

We first specialize the $\mathbf f$-variable at $x_f$. Let
\[
 \kappa_{1}^{\mathrm{BF}}
   =\kappa(\mathbf g,\mathbf h),
 \qquad
 \kappa_{\varepsilon_K}^{\mathrm{BF}}
   =\kappa(\mathbf g,\mathbf h\otimes\varepsilon_K)
\]
denote the two Beilinson--Flach classes obtained from the two character components of $V_f$.  Their precise self-dual twists are understood in the notation, and they take values in
\[
 H^1(\Q,\mathbb V_{1}^{\dagger})
 \quad\text{and}\quad
 H^1(\Q,\mathbb V_{\varepsilon_K}^{\dagger}),
\]
respectively.

Let
\[
 \mathscr L^{\mathrm{HR}}_1
   =L_p^{\mathbf g}\left(\mathbf g,\mathbf h,
        \frac{\ell+m-1}{2}\right),
 \qquad
 \mathscr L^{\mathrm{HR}}_{\varepsilon_K}
   =L_p^{\mathbf g}\left(\mathbf g,\mathbf h\otimes\varepsilon_K,
        \frac{\ell+m-1}{2}\right)
\]
be the two $\mathbf g$-dominant Hida--Rankin $p$-adic $L$-functions occurring in the Artin-formalism factorization. Here $\ell$ and $m$ stand for the weights of $\mathbf g$ and $\mathbf h$. Similarly, we write $\mathscr L_p^g$ for the triple-product $p$-adic $L$-function with dominant variable $\mathbf g$, which is distinguished from a Hida--Rankin function by its three arguments. Up to the explicit Euler and period factor $\mathscr E_{\mathrm{tri}}$, this factorization has the form
\begin{equation}\label{eq:first-factorization}
 \bigl(\mathscr L_p^{g}
        (f,\mathbf g,\mathbf h)\bigr)^2
 =
 \mathscr E_{\mathrm{tri}}\,
 \mathscr L^{\mathrm{HR}}_1\mathscr L^{\mathrm{HR}}_{\varepsilon_K}.
\end{equation}
The factor $\mathscr E_{\mathrm{tri}}$ is the quotient of the local interpolation and period factors in the cited factorization, restricted to $x_f$, and the following comparison uses the same normalizations. More precisely, we use the normalizations of \cite{AOR} for the reciprocity laws, and we assume that in these normalizations the factor $\mathscr E_{\mathrm{tri}}$ is absorbed in the reciprocity law for $\kappa(\mathbf f,\mathbf g,\mathbf h)$, together with the Atkin--Lehner factor and the period introduced below, so that it does not appear in \eqref{eq:diag-BF-comparison}; otherwise, it would appear there as an additional explicit factor.

The factorization \eqref{eq:first-factorization} suggests multiplying each Beilinson--Flach component by the Hida--Rankin function associated with the other component.

\begin{definition}\label{def:weighted-BF-double}
The \emph{weighted Beilinson--Flach class} attached to $(f,\mathbf g,\mathbf h)$ is
\begin{align*}
 \operatorname{BF}^{\mathrm{wt}}(f,\mathbf g,\mathbf h)
 ={}&
 v_{f,1}\otimes
 \mathscr L^{\mathrm{HR}}_{\varepsilon_K}\,
 \kappa_{1}^{\mathrm{BF}}
 \\
 &+
 v_{f,\varepsilon_K}\otimes
 \mathscr L^{\mathrm{HR}}_1\,
 \kappa_{\varepsilon_K}^{\mathrm{BF}}.
\end{align*}
It is regarded as an element of
\[
 H^1\left(
   \Q,
   \left.\mathbb V_{\mathbf f\mathbf g\mathbf h}^{\dagger}
   \right|_{x_f}
 \right)\otimes_{\mathscr R}\mathscr K.
\]
\end{definition}

With the normalizations of \eqref{eq:first-factorization}, both weighted components have regulator $\mathscr L^{\mathrm{HR}}_1\mathscr L^{\mathrm{HR}}_{\varepsilon_K}$.

For the first comparison, one must initially use a local condition which is slightly larger than the balanced one. Put $\mathbb V=\left.\mathbb V_{\mathbf f\mathbf g\mathbf h}^{\dagger}\right|_{x_f}$. With the common self-dual twist understood, define
\[
\begin{aligned}
 F^2\mathbb V
   &=\bigl(V_f\otimes\mathbb V_{\mathbf g}^+\otimes\mathbb V_{\mathbf h}^+
       +V_f^+\otimes\mathbb V_{\mathbf g}\otimes\mathbb V_{\mathbf h}^+
       +V_f^+\otimes\mathbb V_{\mathbf g}^+\otimes\mathbb V_{\mathbf h}\bigr)^{\dagger},\\
 F^{\mathcal G\cup+}\mathbb V
   &=\bigl(V_f\otimes\mathbb V_{\mathbf g}^+\otimes\mathbb V_{\mathbf h}\bigr)^{\dagger}
           +F^2\mathbb V.
\end{aligned}
\]
Here $V_f^+$ is the local line selected by the CM family, which does not agree with either global character summand. The balanced condition at $p$ is the kernel of projection to $H^1(\Q_p,\mathbb V/F^2\mathbb V)$; replacing $F^2$ by $F^{\mathcal G\cup+}$ gives the condition $\mathcal G\cup+$. For global cohomology we use $H^1(\Q_\Sigma/\Q,\mathbb V)$, where $\Sigma$ contains $p$, infinity and the tame level primes, as in \cite[\S5.3]{AOR}.

In the setting of \cite[Proposition~5.14]{AOR}, the local computation places $\operatorname{BF}^{\mathrm{wt}}$ in $H^1_{\mathcal G\cup+}$, and then Corollary~5.19 of loc.\,cit.\ compares it with the diagonal class and proves that it is balanced. Hence, for the statement in families, the rank-one hypothesis concerns this enlarged Selmer group. At a fixed classical specialization, \cite[Corollary~5.25]{AOR} gives a variant; in its boundary case $\ell_0=2$, $m_0=1$, the additional condition is $\alpha_g\beta_h\ne1$.

The integral constructions of \cite[Assumptions~4.1, 4.10 and 5.1]{AOR} require $\mathbf g$ and $\mathbf h$ to be residually irreducible and $p$-distinguished. At the second irregular Eisenstein point, we instead assume that the classes, local maps and reciprocity laws extend to our characteristic-zero realization, with the regularity at $x_h$ specified in \S\ref{subsec:two-Eis-spec}. This extension is the additional input needed to apply the argument of \cite{AOR} to the double specialization.

We use the Atkin--Lehner normalization of \cite[\S4.1]{AOR}: $\lambda_{N_g}(\mathbf g)$ interpolates the scalars defined by $W_{N_g}g_y=\lambda_{N_g}(g_y)g_y^*$ at crystalline cuspidal points, where $g_y^*$ is the conjugate eigenform. Observe that this factor is different from the congruence ideal used to normalize the Perrin--Riou map. In the CM basis of \cite[\S2.4]{AOR}, put $v_f^+=(v_{f,1}+v_{f,\varepsilon_K})/2$. The period is
\[
 \Omega_f=2\langle v_f^+,\omega_f\rangle_f\in L^\times,
\]
where $\omega_f$ is the specialized de Rham differential and $\langle\ ,\ \rangle_f$ is the cohomological pairing used in the reciprocity law; this is \cite[Definition~5.16]{AOR}. Changing the CM basis changes the period compatibly.

\begin{proposition}[First degeneration]\label{thm:first-degeneration}
Assume the characteristic-zero construction, local conditions and reciprocity laws specified above, and suppose that $H^1_{\mathcal G\cup+}(\Q,\mathbb V)\otimes_{\mathscr R}\mathscr K$ has dimension one. Assume moreover that
\[
 \mathscr L_p^{g}(f,\mathbf g,\mathbf h)\neq0.
\]
Then $\operatorname{BF}^{\mathrm{wt}}(f,\mathbf g,\mathbf h)$ satisfies the balanced local conditions over $\mathscr K$, and
\begin{equation}\label{eq:diag-BF-comparison}
 \lambda_{N_g}(\mathbf g)\,
 \operatorname{BF}^{\mathrm{wt}}(f,\mathbf g,\mathbf h)
 =
 \Omega_f\,
 \mathscr L_p^{g}(f,\mathbf g,\mathbf h)\,
 \kappa(f,\mathbf g,\mathbf h).
\end{equation}
Here $\lambda_{N_g}(\mathbf g)$ is the Atkin--Lehner pseudo-eigenvalue just defined, and $\Omega_f$ is the $p$-adic period determined by the choice of $v_{f,1}$ and $v_{f,\varepsilon_K}$.
\end{proposition}

\begin{proof}
The analytic factorization and reciprocity laws identify the two regulator values. Both classes lie in the enlarged Selmer group, where the regulator is nonzero and hence injective over the fraction field. They therefore agree, and the weighted class is balanced because the diagonal class is.
\end{proof}

\subsection{From Beilinson--Flach to Beilinson--Kato classes}
\label{subsec:second-degeneration}

We next specialize $\mathbf h$ at $x_h$, treating the two Beilinson--Flach components separately. For each $i\in\{1,\varepsilon_K\}$, let
\[
 \kappa_{i,1}^{\mathrm{Kato}}
 \in H^1(\Q,\mathbb V_{i,1}^{\dagger}),
 \qquad
 \kappa_{i,2}^{\mathrm{Kato}}
 \in H^1(\Q,\mathbb V_{i,2}^{\dagger})
\]
be the Beilinson--Kato classes attached to the two character components of $V_h$. In concrete terms, these are Kato classes for $\mathbf g$ twisted by the characters determined by $i$ and $\psi_j$.

Let $\mathscr M_{i,1}$ and $\mathscr M_{i,2}$ denote the two Mazur--Kitagawa $p$-adic $L$-functions arising from the Eisenstein specialization of the Hida--Rankin function $\mathscr L^{\mathrm{HR}}_i$.  Thus
\begin{equation}\label{eq:second-factorization}
 \left.\mathscr L^{\mathrm{HR}}_i\right|_{x_h}
 =
 \mathscr E_i\,
 \mathscr M_{i,1}\mathscr M_{i,2},
 \qquad i\in\{1,\varepsilon_K\},
\end{equation}
where $\mathscr E_i$ is the corresponding explicit normalization factor. There are now four Mazur--Kitagawa factors, indexed by $(i,j)$.

\begin{definition}\label{def:weighted-Kato-components}
For $i\in\{1,\varepsilon_K\}$, define the weighted Beilinson--Kato class
\begin{align*}
 \operatorname{BK}^{\mathrm{wt}}_i(\mathbf g,h)
 ={}&
 v_{h,1}\otimes
 \mathscr M_{i,2}\,
 \kappa_{i,1}^{\mathrm{Kato}}
 \\
 &+
 v_{h,2}\otimes
 \mathscr M_{i,1}\,
 \kappa_{i,2}^{\mathrm{Kato}}.
\end{align*}
\end{definition}

The Kato reciprocity laws send both weighted components to $\mathscr M_{i,1}\mathscr M_{i,2}$. In order to compare the classes globally, we first check that the weighted class has zero image in the ordinary local quotient.

For $i\in\{1,\varepsilon_K\}$, let $\mathbb D_i$ be the local $(\varphi,\Gamma)$-module of $\left.\mathbb V_i^{\dagger}\right|_{x_h}$, and let $\mathbb D_{\mathbf g}$ and $D_h$ be those of the two factors. The balanced condition is the kernel of
\[
 H^1(\Q_p,\mathbb D_i)\longrightarrow H^1(\Q_p,\mathbb Q_i),
 \qquad
 \mathbb Q_i=(\mathbb D_{\mathbf g}^-\otimes D_h^-)^{\dagger}_i.
\]
The subscript $i$ includes the character twist from the first specialization, and $\dagger$ includes the fixed self-dual twist. Equivalently, the allowed local submodule is $(\mathbb D_{\mathbf g}^+\otimes D_h+\mathbb D_{\mathbf g}\otimes D_h^+)^{\dagger}_i$. Write $\mathcal S_i$ for the resulting global Selmer module, using the same tame local conditions for both Euler systems. The family Beilinson--Flach class satisfies this balanced condition, and the regular specialization of its local filtration preserves the vanishing in the quotient.

For the weighted Kato class, the mechanism is the one in \cite[proof of Theorem~4.2]{BDV}. In its CM realization, the global character vectors have the form $v_{h,1}=v_h^++v_h^-$ and $v_{h,2}=v_h^+-v_h^-$; see \emph{loc.\,cit.}, equations~(23)--(24). Hence the projection to $D_h^-$ takes the difference of the two coordinates. Let $x_{i,j}\in H^1(\Q_p,\mathbb Q_i)$ denote the ordinary-quotient image of the $j$-th Kato class, after these identifications. With the Kato periods fixed compatibly, its rank-one reciprocity law is
\[
 \operatorname{Col}_i(x_{i,j})=\mathscr M_{i,j}.
\]
For a general $h$, we assume this description of the ordinary quotient, with the two Kato coordinates and their reciprocity laws normalized as above. The following proposition uses this local compatibility in addition to $p$-irregularity.

Write $\mathscr A_i$ and $\Omega_{h,i}$ for the nonzero factors supplied by the reciprocity laws in the normalization of \eqref{eq:second-factorization}, recording the Euler and congruence factors (including the factor $\mathscr E_i$ of \eqref{eq:second-factorization}) and the comparison of regulator periods, respectively.

\begin{proposition}[Beilinson--Flach--Kato comparison]
\label{prop:BF-Kato-components}
Fix $i\in\{1,\varepsilon_K\}$. Assume the regular specialization and local compatibility above, and the tame local conditions for the Kato classes. If $H^1(\Q_p,\mathbb Q_i)$ is torsion-free of rank one and $\operatorname{Col}_i\ne0$, then
\[
 \operatorname{BK}^{\mathrm{wt}}_i(\mathbf g,h)\in\mathcal S_i.
\]

If, moreover, the global regulator on $\mathcal S_i\otimes_{\mathscr R_{f,h}}\mathscr K_{f,h}$ is injective and the two reciprocity laws have the normalization just specified, then
\begin{equation}\label{eq:component-BF-Kato}
 \mathscr A_i\,
 \operatorname{BK}^{\mathrm{wt}}_i(\mathbf g,h)
 =\Omega_{h,i}\,
       \left.\kappa_i^{\mathrm{BF}}\right|_{x_h}
\end{equation}
over $\mathscr K_{f,h}$.
\end{proposition}

\begin{proof}
Projection of the weighted Kato class to the local quotient is
\[
 \mathscr M_{i,2}x_{i,1}-\mathscr M_{i,1}x_{i,2}.
\]
Its Coleman image is $\mathscr M_{i,2}\mathscr M_{i,1}-\mathscr M_{i,1}\mathscr M_{i,2}=0$, and since a nonzero map from a torsion-free rank-one module to a domain is injective, the quotient class itself vanishes. Together with the tame local conditions, this proves that the class belongs to $\mathcal S_i$.

Both classes in \eqref{eq:component-BF-Kato} now lie in the same Selmer module. Their reciprocity laws and \eqref{eq:second-factorization} give the same regulator image after multiplication by $\mathscr A_i$ and $\Omega_{h,i}$, respectively. Global injectivity proves the equality. It also holds before localization whenever both sides lie in the same torsion-free module.
\end{proof}

A sufficient condition for the global injectivity is that $\mathcal S_i$ be torsion-free of rank one and its regulator nonzero. In the CM case, \cite[Proposition~4.1 and Theorem~4.2]{BDV} obtain the comparison from the vanishing of a strict cyclotomic Selmer group. They treat $h=E_1(1,\varepsilon_{K_h})^{\mathrm{ord}}$, where $K_h$ is the imaginary quadratic field of the second CM family, and their period $\Omega_{h,\gamma}=2\langle v_h^+,\omega_h\rangle_h$ is the CM instance of $\Omega_{h,i}$; see \cite[equation~(25)]{BDV}. Their hypotheses cover even-weight cuspidal forms of trivial nebentype, $p$-semistable conductor, and a $p$-regular refinement of slope less than the weight minus one, with $K_h$ of odd discriminant and all primes dividing $p$ and the tame level split in $K_h$, subject to their non-exceptionality condition. To apply this result to both components, these hypotheses must hold for both cuspidal twists. The comparison in families also requires the regularity of the classes and local filtrations at specialization assumed above.

As in \cite[Lemma~5.9 and Propositions~5.10, 5.14]{AOR}, local cancellation concerns the weighted class, whose individual Kato components may have nonzero ordinary-quotient images. For non-CM $h$, the local compatibility and global injectivity in Proposition~\ref{prop:BF-Kato-components} remain assumptions.

We keep the componentwise periods until both comparisons have been normalized. To state the double-degeneration formula with common factors, we assume in addition that the two comparisons have compatible normalizations:
\begin{equation}\label{eq:compatible-normalizations}
 \mathscr A_1=\mathscr A_{\varepsilon_K}=\mathscr A_h,
 \qquad
 \Omega_{h,1}=\Omega_{h,\varepsilon_K}=\Omega_h.
\end{equation}
If \eqref{eq:compatible-normalizations} is not imposed, one may instead absorb $\mathscr A_i/\Omega_{h,i}$ into the definition of the $i$-th component of the doubly weighted class below.

\begin{definition}\label{def:doubly-weighted-Kato}
Assume first that the compatible normalization \eqref{eq:compatible-normalizations} has been chosen. The \emph{doubly weighted Beilinson--Kato class} is
\begin{align}
 \operatorname{BK}^{(2)}(f,\mathbf g,h)
 ={}&
 v_{f,1}\otimes
 \left.\mathscr L^{\mathrm{HR}}_{\varepsilon_K}\right|_{x_h}\,
 \operatorname{BK}^{\mathrm{wt}}_1(\mathbf g,h)
 \nonumber\\
 &+
 v_{f,\varepsilon_K}\otimes
 \left.\mathscr L^{\mathrm{HR}}_1\right|_{x_h}\,
 \operatorname{BK}^{\mathrm{wt}}_{\varepsilon_K}(\mathbf g,h).
 \label{eq:doubly-weighted-Kato}
\end{align}
Without \eqref{eq:compatible-normalizations}, the normalized version is obtained by replacing $\operatorname{BK}^{\mathrm{wt}}_i$ with $(\mathscr A_i/\Omega_{h,i}) \operatorname{BK}^{\mathrm{wt}}_i$ in \eqref{eq:doubly-weighted-Kato}.
\end{definition}

Expanding \eqref{eq:doubly-weighted-Kato} displays the four Kato components created by the two successive Eisenstein specializations:
\begin{align*}
 \operatorname{BK}^{(2)}(f,\mathbf g,h)
 ={}&
 v_{f,1}\otimes v_{h,1}\otimes
 \left.\mathscr L^{\mathrm{HR}}_{\varepsilon_K}\right|_{x_h}
 \mathscr M_{1,2}\,
 \kappa_{1,1}^{\mathrm{Kato}}
 \\
 &+
 v_{f,1}\otimes v_{h,2}\otimes
 \left.\mathscr L^{\mathrm{HR}}_{\varepsilon_K}\right|_{x_h}
 \mathscr M_{1,1}\,
 \kappa_{1,2}^{\mathrm{Kato}}
 \\
 &+
 v_{f,\varepsilon_K}\otimes v_{h,1}\otimes
 \left.\mathscr L^{\mathrm{HR}}_1\right|_{x_h}
 \mathscr M_{\varepsilon_K,2}\,
 \kappa_{\varepsilon_K,1}^{\mathrm{Kato}}
 \\
 &+
 v_{f,\varepsilon_K}\otimes v_{h,2}\otimes
 \left.\mathscr L^{\mathrm{HR}}_1\right|_{x_h}
 \mathscr M_{\varepsilon_K,1}\,
 \kappa_{\varepsilon_K,2}^{\mathrm{Kato}}.
\end{align*}
The Hida--Rankin weights come from the specialization at $f$, and the Mazur--Kitagawa weights from that at $h$.

\begin{proposition}[Second degeneration]\label{thm:second-degeneration}
Under the hypotheses of Proposition~\ref{prop:BF-Kato-components} for both $i=1$ and $i=\varepsilon_K$, and the compatible normalization \eqref{eq:compatible-normalizations}, one has
\begin{equation}\label{eq:BF-doubleKato}
 \mathscr A_h\,
 \operatorname{BK}^{(2)}(f,\mathbf g,h)
 =
 \Omega_h\,
 \left.
 \operatorname{BF}^{\mathrm{wt}}(f,\mathbf g,\mathbf h)
 \right|_{x_h}.
\end{equation}
In particular, the doubly weighted Beilinson--Kato class satisfies the same balanced local conditions as the specialization of the weighted Beilinson--Flach class.
\end{proposition}

\begin{proof}
Apply \eqref{eq:component-BF-Kato} to the two terms in Definition \ref{def:doubly-weighted-Kato}. Since the outer Hida--Rankin weights are unchanged during the componentwise comparison, we obtain
\begin{align*}
 \mathscr A_h\operatorname{BK}^{(2)}(f,\mathbf g,h)
 ={}&
 \Omega_h v_{f,1}\otimes
 \left.\mathscr L^{\mathrm{HR}}_{\varepsilon_K}\right|_{x_h}
 \left.\kappa_1^{\mathrm{BF}}\right|_{x_h}
 \\
 &+
 \Omega_h v_{f,\varepsilon_K}\otimes
 \left.\mathscr L^{\mathrm{HR}}_1\right|_{x_h}
 \left.\kappa_{\varepsilon_K}^{\mathrm{BF}}\right|_{x_h},
\end{align*}
which is precisely the right-hand side of \eqref{eq:BF-doubleKato}. The assertion about local conditions follows from the corresponding assertion for the weighted Beilinson--Flach class.
\end{proof}

\subsection{The double-degeneration formula}\label{subsec:double-formula}

Combining Propositions \ref{thm:first-degeneration} and \ref{thm:second-degeneration} gives the main result of this section.

\begin{theorem}[Double Eisenstein degeneration]\label{thm:double-degeneration}
Assume the hypotheses of Proposition~\ref{thm:first-degeneration} and of Proposition~\ref{prop:BF-Kato-components} for both $i=1$ and $i=\varepsilon_K$. Suppose that the classes and comparison factors specialize regularly at $x_h$ as specified in \S\ref{subsec:two-Eis-spec}, and choose the common normalization \eqref{eq:compatible-normalizations}. Then
\begin{equation}\label{eq:main-double-degeneration}
 \lambda_{N_g}(\mathbf g)\,
 \mathscr A_h\,
 \operatorname{BK}^{(2)}(f,\mathbf g,h)
 =
 \Omega_f\Omega_h\,
 \mathscr L_p^{g}(f,\mathbf g,h)\,
 \kappa(f,\mathbf g,h).
\end{equation}
Here all classes and functions are restricted to the common specialization locus determined by $x_f$ and $x_h$.
\end{theorem}

\begin{proof}
Specializing \eqref{eq:diag-BF-comparison} at $x_h$ gives
\[
 \lambda_{N_g}(\mathbf g)\,
 \left.\operatorname{BF}^{\mathrm{wt}}
 (f,\mathbf g,\mathbf h)\right|_{x_h}
 =
 \Omega_f\,
 \mathscr L_p^{g}(f,\mathbf g,h)\,
 \kappa(f,\mathbf g,h).
\]
Multiplying by $\Omega_h$ and using \eqref{eq:BF-doubleKato} yields \eqref{eq:main-double-degeneration}.
\end{proof}

\begin{remark}\label{rem:period-square}
If the normalizations used in the two comparisons identify $\Omega_f$ and $\Omega_h$ with a single period $\Omega_{f,\gamma}$, then the right-hand side of \eqref{eq:main-double-degeneration} may be written as
\[
 \Omega_{f,\gamma}^{2}\,
 \mathscr L_p^{g}(f,\mathbf g,h)\,
 \kappa(f,\mathbf g,h).
\]
We retain $\Omega_f\Omega_h$ until this compatibility is checked from the interpolation formulas, so that each period records its own Eisenstein specialization.
\end{remark}

\subsection{Comparison with the Heegner-point picture}
\label{subsec:double-Heegner}

The construction of \cite{BSV2} uses two weight-one specializations, at regular theta series, to relate diagonal cycles with Heegner classes. Our irregular Eisenstein specializations instead give reducible fibres, whose character components are weighted by the analytic factorizations.

The link with Heegner points in our comparison comes through the Beilinson--Flach classes of \cite{BDV}: as in Section~\ref{sec:heegner-elliptic}, the weights enforce the local condition and reciprocity gives the global comparison.

\begin{remark}\label{rem:two-weight-one-limits}
The two routes mentioned at the beginning of the section have different parity constraints. The diagonal construction uses the self-dual twist with $2\mathbf t=\mathbf k+\boldsymbol\ell+\mathbf m$ (see \cite[\S4.2]{AOR}). At weights $(1,\ell,1)$, the remaining classical weight $\ell$ is even. A further specialization to $(1,1,1)$ is therefore outside this classical parity range and does not give a diagonal-to-circular-unit comparison in the same construction.

For Beilinson--Flach classes, one can instead ask about specializing both modular variables to weight one while retaining the cyclotomic variable. This is not an application of \cite[Theorem~4.2]{BDV}, where the remaining cuspidal form has even weight at least two. Such a comparison would require regular continuation of the classes and local maps to the second Eisenstein point, followed by a Kato--circular-unit comparison in the appropriate realization. In particular, the CM point $E_1(1,\varepsilon_K)$ has ordinary eigenvalue $1$ and lies outside the non-exceptional hypotheses of Section~\ref{sec:circular-kato}.
\end{remark}

\bibliographystyle{amsalpha}
\raggedbottom
\bibliography{PR-refs}

\end{document}